\documentclass[10pt]{amsart}
\usepackage[cp1251]{inputenc}
\usepackage[english]{babel}
\usepackage{amsmath}
\usepackage{amssymb}
\usepackage{amsfonts}

\begin{document}
\renewcommand{\refname}{References}

\thispagestyle{empty}

\title[Strong Numerical Methods of Orders 2.0, 2.5, and 3.0]
{Strong Numerical Methods of Orders 2.0, 2.5, and 3.0 for
Ito Stochastic Differential Equations Based on the Unified Stochastic Taylor
Expansions and Multiple Fourier--Legendre Series }
\author[D.F. Kuznetsov]{Dmitriy F. Kuznetsov}
\address{Dmitriy Feliksovich Kuznetsov
\newline\hphantom{iii} Peter the Great Saint-Petersburg Polytechnic University,
\newline\hphantom{iii} Polytechnicheskaya ul., 29,
\newline\hphantom{iii} 195251, Saint-Petersburg, Russia}%
\email{sde\_kuznetsov@inbox.ru}

\maketitle {\small
\begin{quote}
\noindent{\sc Abstract.} 
The article is devoted to the construction 
of explicit one-step strong numerical methods
with the orders of convergence 2.0, 2,5, and 3.0 for Ito
stochastic differential equations with
multidimensional non-commutative noise.
We consider the numerical methods based
on the unified Taylor--Ito and Taylor--Stratonovich
expansions.
For numerical modeling of iterated Ito and Stratonovich
stochastic integrals of multiplicities 1 to 6 we apply the method
of multiple Fourier--Legendre series converging in the
sense of norm in Hilbert space $L_2([t, T]^k),$\ $k=1,\ldots,6$.
The article is addressed to engineers who use
numerical modeling in stochastic control
and for solving the non-linear filtering problem. 
The article can be interesting 
for mathematicians who working in the field 
of high-order strong numerical methods for 
Ito stochastic differential equations.
\medskip
\end{quote}
}

\vspace{5mm}


\setlength{\baselineskip}{1.8em}

\tableofcontents

\setlength{\baselineskip}{1.2em}


\section{Introduction}

\vspace{5mm}

Let $(\Omega,$ ${\rm F},$ ${\sf P})$ be a complete probability space, let 
$\{{\rm F}_t, t\in[0,T]\}$ be a nondecreasing 
right-continous family of $\sigma$-algebras of ${\rm F},$
and let ${\bf f}_t$ be a standard $m$-dimensional Wiener 
stochastic process, which is
${\rm F}_t$-measurable for any $t\in[0, T].$ We assume that the components
${\bf f}_{t}^{(i)}$ $(i=1,\ldots,m)$ of this process are independent. Consider
an Ito stochastic differential equation (SDE) in the integral form

\vspace{-1mm}
\begin{equation}
\label{1.5.2}
{\bf x}_t={\bf x}_0+\int\limits_0^t {\bf a}({\bf x}_{\tau},\tau)d\tau+
\int\limits_0^t \Sigma({\bf x}_{\tau},\tau)d{\bf f}_{\tau},\ 
{\bf x}_0={\bf x}(0,\omega).
\end{equation}

\vspace{2mm}
\noindent
Here ${\bf x}_t$ is some $n$-dimensional stochastic process 
satisfying the equation (\ref{1.5.2}). 
The nonrandom functions ${\bf a}: \mathbb{R}^n\times[0, T]\to\mathbb{R}^n$,
$\Sigma: \mathbb{R}^n\times[0, T]\to\mathbb{R}^{n\times m}$
guarantee the existence and uniqueness up to stochastic 
equivalence of a solution
of the equation (\ref{1.5.2}) \cite{1}. 
The second integral on the right-hand side of (\ref{1.5.2}) is 
interpreted as the Ito stochastic integral.
Let ${\bf x}_0$ be an $n$-dimensional random variable, which is 
${\rm F}_0$-measurable and 
${\sf M}\{\left|{\bf x}_0\right|^2\}<\infty$ 
(${\sf M}$ denotes a mathematical expectation).
We assume that
${\bf x}_0$ and ${\bf f}_t-{\bf f}_0$ are independent when $t>0.$

It is well known \cite{KlPl2}-\cite{2006}
that Ito SDEs are 
adequate mathematical models of dynamic systems of different 
physical nature under 
the influence of random disturbances. One of the effective approaches 
to the numerical integration of 
Ito SDEs is an approach based on the Taylor--Ito and 
Taylor--Stratonovich expansions
\cite{KlPl2}-\cite{arxiv-24999}. The most important feature of such 
expansions is a presence in them of the so-called iterated
Ito and Stratonovich stochastic integrals, which play the key 
role for solving the 
problem of numerical integration of Ito SDEs 
and have the 
following form

\vspace{-1mm}
\begin{equation}
\label{ito}
J[\psi^{(k)}]_{T,t}=\int\limits_t^T\psi_k(t_k) \ldots \int\limits_t^{t_{2}}
\psi_1(t_1) d{\bf w}_{t_1}^{(i_1)}\ldots
d{\bf w}_{t_k}^{(i_k)},
\end{equation}

\begin{equation}
\label{str}
J^{*}[\psi^{(k)}]_{T,t}=
\int\limits_t^{*T}\psi_k(t_k) \ldots \int\limits_t^{*t_{2}}
\psi_1(t_1) d{\bf w}_{t_1}^{(i_1)}\ldots
d{\bf w}_{t_k}^{(i_k)},
\end{equation}

\vspace{2mm}
\noindent
where every $\psi_l(\tau)$ $(l=1,\ldots,k)$ is a 
non-random function 
on $[t,T],$ ${\bf w}_{\tau}^{(i)}={\bf f}_{\tau}^{(i)}$
for $i=1,\ldots,m$ and
${\bf w}_{\tau}^{(0)}=\tau,$\
$i_1,\ldots,i_k = 0, 1,\ldots,m,$
$$
\int\limits\ \hbox{\rm and} 
\int\limits^{*}
$$

\vspace{1mm}
\noindent
denote Ito and 
Stratonovich stochastic integrals,
respectively (in this paper, 
we use the definition of the Stratonovich stochastic integral from \cite{KlPl2}).

Note that $\psi_l(\tau)\equiv 1$ $(l=1,\ldots,k)$ and
$i_1,\ldots,i_k = 0, 1,\ldots,m$ in  
\cite{KlPl2}-\cite{Mi3}, \cite{PW1}, \cite{KlPl1}. At the same time
$\psi_l(\tau)\equiv (t-\tau)^{q_l}$ ($l=1,\ldots,k$; 
$q_1,\ldots,q_k=0, 1, 2,\ldots $) and $i_1,\ldots,i_k = 1,\ldots,m$ in
\cite{KK5}-\cite{arxiv-24999}.
 
Effective solution 
of the problem of
combined mean-square approximation of 
the iterated Ito and Stratonovich stochastic integrals
(\ref{ito}) and (\ref{str}) of multiplicities 1 to 6
composes one of the subjects of this article.

We want to mention in short that there are 
two main criteria of numerical methods convergence 
for Ito SDEs \cite{KlPl2}-\cite{Mi3}:  
a strong or mean-square
criterion and a 
weak criterion where the subject of approximation is not the solution 
of Ito SDE, simply stated, but the 
distribution of Ito SDE solution.

Using strong numerical methods, we may build
sample pathes
of Ito SDEs numerically. 
These methods require combined mean-square approximation 
of the iterated Ito and Stratonovich stochastic integrals
(\ref{ito}) and (\ref{str}). 

Also numerical
integration of Ito SDEs based on the strong convergence criterion of 
approximation is widely used
for the numerical solution
of different mathematical problems connected with Ito SDEs. 
Among these problems, we note the
following: signal filtering under influence of random
noises in 
various statements, optimal stochastic control, testing estimation
procedures of parameters of stochastic systems, stochastic stability 
and bifurcations analysis \cite{KlPl2}-\cite{Mi3}. 

The problem of effective jointly numerical modeling 
(in accordance to the mean-square convergence criterion) of the iterated 
Ito and Stratonovich stochastic integrals 
(\ref{ito}) and (\ref{str}) is 
difficult from 
theoretical and computing point of view \cite{KlPl2}-\cite{rr}.

The only exception is connected with a narrow particular case, when 
$i_1=\ldots=i_k\ne 0$ and
$\psi_1(s),\ldots,\psi_k(s)\equiv \psi(s)$.
This case allows 
the investigation with using of the Ito formula 
\cite{KlPl2}-\cite{Mi3}.

Note that even for the mentioned coincidence ($i_1=\ldots=i_k\ne 0$), 
but for different 
functions $\psi_1(s),\ldots,\psi_k(s)$ the mentioned 
difficulties persist, and 
relatively simple families of 
iterated Ito and Stratonovich stochastic integrals, 
which can be often 
met in the applications, cannot be represented effectively in a finite 
form (for the mean-square approximation) using the system of standard 
Gaussian random variables.

Note that for a number of special types of Ito SDEs 
the problem of approximation of iterated
stochastic integrals can be simplified but cannot be solved. The equations
with additive scalar noise, with additive vector noise, with non-additive  
scalar noise, with a small parameter are related to such 
types of equations \cite{KlPl2}-\cite{Mi3}. 
For the mentioned types of equations, simplifications 
are connected with the fact that either some coefficient functions 
from stochastic analogues of the Taylor formula identically equal to zero, 
or scalar noise has a strong effect, or due to the presence 
of a small parameter we may neglect some members from the stochastic 
analogues of the Taylor formula, which include difficult for approximation 
iterated stochastic integrals \cite{KlPl2}-\cite{Mi3}, \cite{Bruti}.
In this article, we consider Ito SDEs 
with multidimensional, non-additive and non-commutative noise.

Seems that iterated stochastic integrals may be approximated by multiple 
integral sums of different types \cite{Mi2}, \cite{Mi3}, \cite{xxxxa}. 
However, this approach implies partitioning of the interval 
of integration $[t, T]$ of iterated stochastic integrals 
(the length $T-t$ of this interval is a small 
value, because it is a step of integration of numerical methods for 
Ito SDEs) and according to numerical 
experiments this additional partitioning leads to significant calculating 
costs \cite{2006}.

In \cite{Mi2} (also see \cite{KlPl2}, \cite{Mi3}, 
\cite{KPW}, \cite{Bruti}), 
Milstein proposed to expand (\ref{ito}) or (\ref{str})
into iterated series in terms of products
of standard Gaussian random variables by representing the Wiener
process as a trigonometric Fourier series with random coefficients 
(the version of the so-called Karhunen--Loeve expansion
of the Brownian bridge process).
To obtain the Milstein expansion of (\ref{str}), the truncated Fourier
expansions of components of the Wiener process ${\bf f}_s$ must be
iteratively substituted in the single integrals, and the integrals
must be calculated, starting from the innermost integral.
This is a complicated procedure that does not lead to a general
expansion of (\ref{str}) valid for an arbitrary multiplicity $k.$
For this reason, only expansions of single, double, and triple
stochastic integrals of the form (\ref{str}) were presented 
in \cite{KlPl2}, \cite{KPS}-\cite{Bruti} ($k=1, 2, 3$)
and in \cite{Mi2}, \cite{Mi3} ($k=1, 2$) 
for the simplest case $\psi_1(s), \psi_2(s), \psi_3(s)\equiv 1;$ 
$i_1, i_2, i_3=0,1,\ldots,m.$

Moreover, 
in \cite{KlPl2}
(Sect.~5.8, pp.~202--204), \cite{KPS} (pp.~82-84),
\cite{KPW} (pp.~438-439),  
\cite{Bruti} (pp.~263-264) 
the authors use (without rigorous proof)
the Wong--Zakai approximation
\cite{W-Z-1}-\cite{Watanabe}
within the frames of the Milstein approach \cite{Mi2}
based on the Karhunen--Loeve expansion of the Brownian bridge
process (see discussions in \cite{2018a} (Sect.~2.18, 6.2),
\cite{2018aa}, \cite{2018aaa} (Sect.~2.6.2, 6.2), \cite{arxiv-4}).

Note that in \cite{rr} the method of expansion
of the iterated Ito stochastic integrals (\ref{ito}) 
($k=2;$ $\psi_1(s), \psi_2(s) \equiv 1;$ $i_1, i_2 =1,\ldots,m$) 
based on expansion
of the Wiener process using Haar functions and 
trigonometric functions has been considered.

It is necessary to note that the Milstein approach \cite{Mi2} 
excelled in several times or even in several orders
the methods of multiple integral sums \cite{Mi2}, \cite{Mi3}, \cite{xxxxa}
considering computational costs in the sense 
of their diminishing.

An alternative strong approximation method was 
proposed for (\ref{str}) in \cite{kuz1997}-\cite{arxiv100} 
(also see \cite{2010-1}-\cite{5-001}),
where $J^{*}[\psi^{(k)}]_{T,t}$ was represented as a multiple stochastic 
integral
from the certain discontinuous non-random function of $k$ variables, and the 
function
was then expressed as an iterated generalized Fourier series in a complete
systems of continuous functions that are orthonormal in $L_2([t, T]).$
In \cite{kuz1997}-\cite{arxiv100} 
(also see \cite{2010-1}-\cite{5-001}) the cases of Legendre polynomials and
trigonometric functions are considered in details.
As a result,
a general iterated series expansion of (\ref{str}) in terms of products
of standard Gaussian random variables was obtained in 
\cite{kuz1997}-\cite{arxiv100} 
(also see \cite{2010-1}-\cite{5-001}) for an arbitrary multiplicity $k.$
Hereinafter, this method is referred to as the method of generalized
iterated Fourier series.

It was shown in \cite{kuz1997}, \cite{kuz1997a} 
(also see \cite{2017}-\cite{2013}) that the method of 
generalized iterated Fourier series leads to the 
Milstein expansion \cite{Mi2} of (\ref{str})
in the case of trigonometric 
functions (at least for 
$k=2;$ $\psi_1(s), \psi_2(s) \equiv 1;$ $i_1, i_2 =1,\ldots,m$) 
and to a substantially simpler expansion of (\ref{str}) in the case
of Legendre polynomials.

Note that the method of generalized 
iterated Fourier series as well as the Milstein approach
\cite{Mi2}
lead to iterated application of the operation of limit transition. 
This problem appears for triple stochastic integrals
($i_1, i_2, i_3=1,\ldots,m$)
or even for some double stochastic integrals 
in the case, when $\psi_1(\tau),$ 
$\psi_2(\tau)\not\equiv 1$ ($i_1, i_2=1,\ldots,m$)
\cite{kuz1997}, \cite{kuz1997a} 
(also see \cite{2017}-\cite{2013}).

The mentioned problem (iterated application of the operation 
of limit transition) not appears 
in the method, which 
is considered for (\ref{ito}) in Theorems 1, 2 (see below)
\cite{2006}, \cite{2017}-\cite{2013}, \cite{5-002}-\cite{2018aaa}.
The idea of this method is as follows: 
the iterated Ito stochastic 
integral (\ref{ito}) of multiplicity $k$ is represented as 
the multiple stochastic 
integral from the certain discontinuous non-random function of $k$ variables 
defined on the hypercube $[t, T]^k$, where $[t, T]$ is the interval of 
integration of the iterated Ito stochastic 
integral (\ref{ito}). Then, 
the indicated 
non-random function of $k$ variables 
is expanded in the hypercube into the generalized 
multiple Fourier series converging 
in the mean-square sense
in the space 
$L_2([t,T]^k)$. After a number of nontrivial transformations we come 
(see Theorems 1, 2 below) to the 
mean-square converging expansion of the
iterated Ito stochastic 
integral (\ref{ito})
into the multiple 
series in terms of products
of standard  Gaussian random 
variables. The coefficients of this 
series are the coefficients of 
generalized multiple Fourier series for the mentioned non-random function 
of $k$ variables, which can be calculated using the explicit formula 
regardless of the multiplicity $k$ of the
iterated Ito stochastic 
integral (\ref{ito}).
Hereinafter, this method is referred to as the method of generalized
multiple Fourier series.

Thus, we obtain the following useful possibilities
of the method of generalized multiple Fourier series.

1. There is the explicit formula (see (\ref{ppppa})) for calculation 
of expansion coefficients 
of the iterated Ito stochastic integral (\ref{ito}) with any
fixed multiplicity $k$. 

2. We have new possibilities for exact calculation of the mean-square 
error of approximation 
of the iterated Ito stochastic integral (\ref{ito})
(see \cite{2017-1}, \cite{2017-1a},
\cite{5-002}-\cite{5-005}, \cite{5-007}, \cite{5-008},
\cite{arxiv-2}, \cite{2018a}-\cite{2018aaa}).

3. Since the used
multiple Fourier series is a generalized in the sense
that it is built using various complete orthonormal
systems of functions in the space $L_2([t, T])$, then we 
have new possibilities 
for approximation --- we may 
use not only the trigonometric functions as in \cite{KlPl2}-\cite{Mi3}
but the Legendre polynomials.

4. As it turned out (see \cite{2006}, \cite{kuz1997},
\cite{kuz1997a}, \cite{2017}-\cite{5-008}, \cite{arxiv-3}, 
\cite{arxiv-8}, \cite{arxiv-44}, \cite{99999}, \cite{arxiv-12}-\cite{5-010},
\cite{2018a}-\cite{2018aaa}),
it is more convenient to work 
with Legendre polynomials for constructing of approximations 
of the iterated Ito stochastic integrals (\ref{ito}). 
Approximations based on the Legendre polynomials essentially simpler 
than their analogues based on the trigonometric functions.
Another advantages of Legendre polynomials 
in the framework of the mentioned problem are considered
in \cite{5-005}, \cite{5-010}, \cite{2018a}-\cite{2018aaa}.

5. The approach based on the Karhunen--Loeve expansion
of the Brownian bridge process (also see \cite{rr})
leads to 
iterated application of the operation of limit
transition (the operation of limit transition 
is implemented only once in Theorems 1, 2 (see below))
starting from  
the second multiplicity (in the general case) 
and third multiplicity (for the case
$\psi_1(s), \psi_2(s), \psi_3(s)\equiv 1;$ 
$i_1, i_2, i_3=1,\ldots,m$)
of iterated Ito stochastic integrals.
Multiple series (the operation of limit transition 
is implemented only once) are more convenient 
for approximation than the iterated ones
(iterated application of the operation of limit
transition), 
since partial sums of multiple series converge for any possible case of  
convergence to infinity of their upper limits of summation 
(let us denote them as $p_1,\ldots, p_k$). 
For example, 
when $p_1=\ldots=p_k=p\to\infty$. For iterated series, 
the condition $p_1=\ldots=p_k=p\to\infty$ obviously 
does not guarantee the convergence of this series.

However, the authors of the works
\cite{KlPl2}
(Sect.~5.8, pp.~202--204), \cite{KPS} (pp.~82-84),
\cite{KPW} (pp.~438-439),  
\cite{Bruti} (pp.~263-264) unreasonably use 
the condition $p_1=p_2=p_3=p\to\infty$
within the frames
of the Milstein approach 
\cite{Mi2} based on the series expansion 
of the Brownian bridge process together
with the Wong--Zakai approximation 
\cite{W-Z-1}-\cite{Watanabe}.

\vspace{5mm}

\section{Explicit One-Step Strong Numerical Schemes of Orders 2.0,
2.5, and 3.0 Based
on the Unified Taylor--Ito expansion}

\vspace{5mm}

Consider the partition $\{\tau_j\}_{j=0}^N$ of the 
interval $[0, \bar T]$ such that

\vspace{-1mm}
$$
0=\tau_0<\ldots <\tau_N=\bar T,\ \ \
\Delta_N=
\hbox{\vtop{\offinterlineskip\halign{
\hfil#\hfil\cr
{\rm max}\cr
$\stackrel{}{{}_{0\le j\le N-1}}$\cr
}} }\Delta\tau_j,\ \ \ \Delta\tau_j=\tau_{j+1}-\tau_j.
$$

\vspace{2mm}

Let ${\bf y}_{\tau_j}\stackrel{\sf def}{=}
{\bf y}_{j};$ $j=0, 1,\ldots,N$ be a time discrete approximation
of the process ${\bf x}_s,$ $s\in[0,\bar T],$ which is a solution of the Ito
SDE (\ref{1.5.2}). 

{\bf Definition 1}\ \cite{KlPl2}.
{\it We will say that a time discrete approximation 
${\bf y}_{j};$ $j=0, 1,\ldots,N,$
corresponding to the maximal step of discretization $\Delta_N,$
converges strongly with order
$\gamma>0$ at time moment 
$\bar T$ to the process ${\bf x}_s,$ $s\in[0,\bar T]$
if there exists a constant $C>0,$ which does not depend on 
$\Delta_N,$ and a $\delta>0$ such that 

\vspace{-1mm}
$$
{\sf M}\{|{\bf x}_{\bar T}-{\bf y}_{\bar T}|\}\le
C(\Delta_N)^{\gamma}
$$

\vspace{2mm}
\noindent
for each $\Delta_N\in(0, \delta).$}

Consider the following
explicit one-step strong numerical scheme of order 3.0
based on the so-called unified Taylor--Ito expansion 
\cite{2006}, \cite{arxiv-24999}, \cite{2017-1}-\cite{2010-1}, 
\cite{2018a}-\cite{2018aaa}

\vspace{1mm}
$$
{\bf y}_{p+1}={\bf y}_p+\sum_{i_{1}=1}^{m}\Sigma_{i_{1}}
\hat I_{0_{\tau_{p+1},\tau_p}}^{(i_{1})}+\Delta{\bf a}
+\sum_{i_{1},i_{2}=1}^{m}G_0^{(i_{2})}
\Sigma_{i_{1}}\hat I_{00_{\tau_{p+1},\tau_p}}^{(i_{2}i_{1})}+
$$

\vspace{2mm}
$$
+
\sum_{i_{1}=1}^{m}\left[G_0^{(i_{1})}{\bf a}\left(
\Delta \hat I_{0_{\tau_{p+1},\tau_p}}^{(i_{1})}+
\hat I_{1_{\tau_{p+1},\tau_p}}^{(i_{1})}\right)
-L\Sigma_{i_{1}}\hat I_{1_{\tau_{p+1},\tau_p}}^{(i_{1})}\right]+
$$

\vspace{2mm}
$$
+\sum_{i_{1},i_{2},i_{3}=1}^{m} G_0^{(i_{3})}G_0^{(i_{2})}
\Sigma_{i_{1}}\hat I_{000_{\tau_{p+1},\tau_p}}^{(i_{3}i_{2}i_{1})}+
\frac{\Delta^2}{2}L{\bf a}+
$$

\vspace{2mm}
$$
+\sum_{i_{1},i_{2}=1}^{m}
\left[G_0^{(i_{2})}L\Sigma_{i_{1}}\left(
\hat I_{10_{\tau_{p+1},\tau_p}}^{(i_{2}i_{1})}-
\hat I_{01_{\tau_{p+1},\tau_p}}^{(i_{2}i_{1})}
\right)
-LG_0^{(i_{2})}\Sigma_{i_{1}}
\hat I_{10_{\tau_{p+1},\tau_p}}^{(i_{2}i_{1})}
+\right.
$$

\vspace{2mm}
$$
\left.+G_0^{(i_{2})}G_0^{(i_{1})}{\bf a}\left(
\hat I_{01_{\tau_{p+1},\tau_p}}
^{(i_{2}i_{1})}+\Delta \hat I_{00_{\tau_{p+1},\tau_p}}^{(i_{2}i_{1})}
\right)\right]+
$$

\vspace{2mm}
\begin{equation}
\label{4.45}
+
\sum_{i_{1},i_2,i_3,i_{4}=1}^{m}G_0^{(i_{4})}G_0^{(i_{3})}G_0^{(i_{2})}
\Sigma_{i_{1}} \hat I_{0000_{\tau_{p+1},\tau_p}}^{(i_{4}i_{3}i_{2}i_{1})}+
{\bf u}_{p+1,p}+{\bf v}_{p+1,p},
\end{equation}

\vspace{10mm}

$$
{\bf u}_{p+1,p}=
\sum_{i_{1}=1}^{m}\Biggl[G_0^{(i_{1})}L{\bf a}\Biggl(\frac{1}{2}
\hat I_{2_{\tau_{p+1},\tau_p}}
^{(i_{1})}+\Delta \hat I_{1_{\tau_{p+1},\tau_p}}^{(i_{1})}+
\frac{\Delta^2}{2}\hat I_{0_{\tau_{p+1},\tau_p}}^{(i_{1})}\Biggr)+\Biggr.
$$

\vspace{2mm}
$$
\Biggl.+\frac{1}{2}LL\Sigma_{i_{1}}\hat I_{2_{\tau_{p+1},\tau_p}}^{(i_{1})}-
LG_0^{(i_{1})}{\bf a}\left(\hat I_{2_{\tau_{p+1},\tau_p}}^{(i_{1})}+
\Delta \hat I_{1{\tau_{p+1},\tau_p}}^{(i_{1})}\right)\Biggr]+
$$

\vspace{2mm}
$$
+
\sum_{i_{1},i_{2},i_{3}=1}^m\left[
G_0^{(i_{3})}LG_0^{(i_{2})}\Sigma_{i_{1}}
\left(\hat I_{100_{\tau_{p+1},\tau_p}}
^{(i_{3}i_{2}i_{1})}-\hat I_{010_{\tau_{p+1},\tau_p}}
^{(i_{3}i_{2}i_{1})}\right)
\right.+
$$

\vspace{2mm}
$$
+G_0^{(i_{3})}G_0^{(i_{2})}L\Sigma_{i_{1}}\left(
\hat I_{010_{\tau_{p+1},\tau_p}}^{(i_{3}i_{2}i_{1})}-
\hat I_{001_{\tau_{p+1},\tau_p}}^{(i_{3}i_{2}i_{1})}\right)+
$$

\vspace{2mm}
$$
+G_0^{(i_{3})}G_0^{(i_{2})}G_0^{(i_{1})}{\bf a}
\left(\Delta \hat I_{000_{\tau_{p+1},\tau_p}}^{(i_{3}i_{2}i_{1})}+
\hat I_{001_{\tau_{p+1},\tau_p}}^{(i_{3}i_{2}i_{1})}\right)-
$$

\vspace{2mm}
$$
\left.-LG_0^{(i_{3})}G_0^{(i_{2})}\Sigma_{i_{1}}
\hat I_{100_{\tau_{p+1},\tau_p}}^{(i_{3}i_{2}i_{1})}\right]+
$$

\vspace{2mm}
$$
+\sum_{i_{1},i_2,i_3,i_4,i_{5}=1}^m
G_0^{(i_{5})}G_0^{(i_{4})}G_0^{(i_{3})}G_0^{(i_{2})}\Sigma_{i_{1}}
\hat I_{00000_{\tau_{p+1},\tau_p}}^{(i_{5}i_{4}i_{3}i_{2}i_{1})}+
$$

\vspace{2mm}
$$
+\frac{\Delta^3}{6}LL{\bf a},
$$

\vspace{10mm}

$$
{\bf v}_{p+1,p}=\sum_{i_{1},i_{2}=1}^{m}
\Biggl[G_0^{(i_{2})}G_0^{(i_{1})}L{\bf a}\left(
\frac{1}{2}\hat I_{02_{\tau_{p+1},\tau_p}}^{(i_{2}i_{1})}+
\Delta \hat I_{01_{\tau_{p+1},\tau_p}}^{(i_{2}i_{1})}+\frac{\Delta^2}{2}
\hat I_{00_{\tau_{p+1},\tau_p}}^{(i_{2}i_{1})}\right)+
$$

\vspace{2mm}
$$
+
\frac{1}{2}LLG_0^{(i_{2})}\Sigma_{i_{1}}
\hat I_{20_{\tau_{p+1},\tau_p}}^{(i_{2}i_{1})}
+
$$

\vspace{2mm}
$$
+G_0^{(i_{2})}LG_0^{(i_{1})}{\bf a}\left(
\hat I_{11_{\tau_{p+1},\tau_p}}
^{(i_{2}i_{1})}-\hat I_{02_{\tau_{p+1},\tau_p}}^{(i_{2}i_{1})}+
\Delta\left(\hat I_{10_{\tau_{p+1},\tau_p}}
^{(i_{2}i_{1})}-\hat I_{01_{\tau_{p+1},\tau_p}}^{(i_{2}i_{1})}
\right)\right)+
$$

\vspace{2mm}
$$
+LG_0^{(i_{2})}L\Sigma_{i_1}\left(
\hat I_{11_{\tau_{p+1},\tau_p}}
^{(i_{2}i_{1})}-\hat I_{20_{\tau_{p+1},\tau_p}}^{(i_{2}i_{1})}\right)+
$$

\vspace{2mm}
$$
+G_0^{(i_{2})}LL\Sigma_{i_1}\left(
\frac{1}{2}\hat I_{02_{\tau_{p+1},\tau_p}}^{(i_{2}i_{1})}+
\frac{1}{2}\hat I_{20_{\tau_{p+1},\tau_p}}^{(i_{2}i_{1})}-
\hat I_{11_{\tau_{p+1},\tau_p}}^{(i_{2}i_{1})}\right)-
$$

\vspace{2mm}
$$
\Biggl.-LG_0^{(i_{2})}G_0^{(i_{1})}{\bf a}\left(
\Delta \hat I_{10_{\tau_{p+1},\tau_p}}
^{(i_{2}i_{1})}+\hat I_{11_{\tau_{p+1},\tau_p}}^{(i_{2}i_{1})}\right)
\Biggr]+
$$

\vspace{2mm}
$$
+
\sum_{i_{1},i_2,i_3,i_{4}=1}^m\Biggl[
G_0^{(i_{4})}G_0^{(i_{3})}G_0^{(i_{2})}G_0^{(i_{1})}{\bf a}
\left(\Delta \hat I_{0000_{\tau_{p+1},\tau_p}}
^{(i_4i_{3}i_{2}i_{1})}+\hat I_{0001_{\tau_{p+1},\tau_p}}
^{(i_4i_{3}i_{2}i_{1})}\right)
+\Biggr.
$$

\vspace{2mm}
$$
+G_0^{(i_{4})}G_0^{(i_{3})}LG_0^{(i_{2})}\Sigma_{i_1}
\left(\hat I_{0100_{\tau_{p+1},\tau_p}}
^{(i_4i_{3}i_{2}i_{1})}-\hat I_{0010_{\tau_{p+1},\tau_p}}
^{(i_4i_{3}i_{2}i_{1})}\right)-
$$

\vspace{2mm}
$$
-LG_0^{(i_{4})}G_0^{(i_{3})}G_0^{(i_{2})}\Sigma_{i_1}
\hat I_{1000_{\tau_{p+1},\tau_p}}
^{(i_4i_{3}i_{2}i_{1})}+
$$

\vspace{2mm}
$$
+G_0^{(i_{4})}LG_0^{(i_{3})}G_0^{(i_{2})}\Sigma_{i_1}
\left(\hat I_{1000_{\tau_{p+1},\tau_p}}
^{(i_4i_{3}i_{2}i_{1})}-\hat I_{0100_{\tau_{p+1},\tau_p}}
^{(i_4i_{3}i_{2}i_{1})}\right)+
$$

\vspace{2mm}
$$
\Biggl.+G_0^{(i_{4})}G_0^{(i_{3})}G_0^{(i_{2})}L\Sigma_{i_1}
\left(\hat I_{0010_{\tau_{p+1},\tau_p}}
^{(i_4i_{3}i_{2}i_{1})}-\hat I_{0001_{\tau_{p+1},\tau_p}}
^{(i_4i_{3}i_{2}i_{1})}\right)\Biggr]+
$$

\vspace{2mm}
$$
+\sum_{i_{1},i_2,i_3,i_4,i_5,i_{6}=1}^m
G_0^{(i_{6})}G_0^{(i_{5})}
G_0^{(i_{4})}G_0^{(i_{3})}G_0^{(i_{2})}\Sigma_{i_{1}}
\hat I_{000000_{\tau_{p+1},\tau_p}}^{(i_6i_{5}i_{4}i_{3}i_{2}i_{1})},
$$

\vspace{6mm}
\noindent
where $\Delta={\bar T}/N$ $(N>1)$ is a constant (for simplicity)
step of integration,
$\tau_p=p\Delta$ $(p=0, 1,\ldots,N)$,
$\hat 
I_{{l_1\ldots\hspace{0.2mm} l_k}_{\hspace{0.2mm}s,t}}^{(i_1\ldots i_k)}$ 
is an 
approximation of the iterated
Ito stochastic integral 

\vspace{-1mm}
\begin{equation}
\label{ll1}
I_{{l_1\ldots\hspace{0.2mm} l_k}_{\hspace{0.2mm}s,t}}^{(i_1\ldots i_k)}=
 \int\limits^ {s} _ {t} (t-\tau _
{k}) ^ {l_ {k}} 
\ldots \int\limits^ {\tau _ {2}} _ {t} (t-\tau _ {1}) ^ {l_ {1}} d
{\bf f} ^ {(i_ {1})} _ {\tau_ {1}} \ldots 
d {\bf f} _ {\tau_ {k}} ^ {(i_ {k})},
\end{equation}

\vspace{2mm}
$$
L= {\partial \over \partial t}
+ \sum^ {n} _ {i=1} {\bf a}_i ({\bf x},  t) 
{\partial  \over  \partial  {\bf  x}_i}
+ {1\over 2} \sum^ {m} _ {j=1} \sum^ {n} _ {l,i=1}
\Sigma_{lj} ({\bf x}, t) \Sigma_{ij} ({\bf x}, t) {\partial
^{2} \over \partial {\bf x}_l \partial {\bf x}_i},
$$

\vspace{3mm}
$$
G_0^{(i)} = \sum^ {n} _ {j=1} \Sigma_{ji} ({\bf x}, t)
{\partial  \over \partial {\bf x}_j},\ \ \
i=1,\ldots,m,
$$

\vspace{4mm}
\noindent
$l_1,\ldots, l_k=0, 1, 2\ldots,$\ \
$i_1,\ldots, i_k=1,\ldots,m,$\ \ $k=1, 2,\ldots$,\ \
$\Sigma_i$ is the $i$th column of the matrix function $\Sigma$ and 
$\Sigma_{ij}$ is the $ij$th component
of the matrix function 
$\Sigma$,\  ${\bf a}_i$ is the $i$th component of the vector function 
${\bf a}$ and ${\bf x}_i$ 
is the
$i$th component of the column ${\bf x}$,\ 
the columns 

\vspace{-2mm}
$$
\Sigma_{i_{1}},\ {\bf a},\ G_0^{(i_{2})}\Sigma_{i_{1}},\
G_0^{(i_{1})}{\bf a},\  L\Sigma_{i_{1}},\ G_0^{(i_{3})}G_0^{(i_{2})}\Sigma_{i_{1}},\
L{\bf a},\ G_0^{(i_{2})}L\Sigma_{i_{1}},\ LG_0^{(i_{2})}\Sigma_{i_{1}},\
G_0^{(i_{2})}G_0^{(i_{1})}{\bf a},\ 
$$

\vspace{-4mm}
$$
G_0^{(i_{4})}G_0^{(i_{3})}G_0^{(i_{2})}\Sigma_{i_{1}},\
G_0^{(i_{1})}L{\bf a},\ LL\Sigma_{i_{1}},\
LG_0^{(i_{1})}{\bf a},\ 
G_0^{(i_{3})}LG_0^{(i_{2})}\Sigma_{i_{1}},\ 
G_0^{(i_{3})}G_0^{(i_{2})}L\Sigma_{i_{1}},\ 
G_0^{(i_{3})}G_0^{(i_{2})}G_0^{(i_{1})}{\bf a},\
$$

\vspace{-4mm}
$$
LG_0^{(i_{3})}G_0^{(i_{2})}\Sigma_{i_{1}},\
G_0^{(i_{5})}G_0^{(i_{4})}G_0^{(i_{3})}G_0^{(i_{2})}\Sigma_{i_{1}},\ LL{\bf a},\
G_0^{(i_{2})}G_0^{(i_{1})}L{\bf a},\
LLG_0^{(i_{2})}\Sigma_{i_{1}},\
G_0^{(i_{2})}LG_0^{(i_{1})}{\bf a},\
LG_0^{(i_{2})}L\Sigma_{i_1},\
$$

\vspace{-3mm}
$$
G_0^{(i_{2})}LL\Sigma_{i_1},\
LG_0^{(i_{2})}G_0^{(i_{1})}{\bf a},\
G_0^{(i_{4})}G_0^{(i_{3})}G_0^{(i_{2})}G_0^{(i_{1})}{\bf a},\
G_0^{(i_{4})}G_0^{(i_{3})}LG_0^{(i_{2})}\Sigma_{i_1},\
LG_0^{(i_{4})}G_0^{(i_{3})}G_0^{(i_{2})}\Sigma_{i_1},\
$$

\vspace{-2mm}
$$
G_0^{(i_{4})}LG_0^{(i_{3})}G_0^{(i_{2})}\Sigma_{i_1},\
G_0^{(i_{4})}G_0^{(i_{3})}G_0^{(i_{2})}L\Sigma_{i_1},\ 
G_0^{(i_{6})}G_0^{(i_{5})}G_0^{(i_{4})}G_0^{(i_{3})}G_0^{(i_{2})}\Sigma_{i_{1}}
$$

\vspace{4mm}
\noindent
are calculated at the point $({\bf y}_p,p).$

It is well known \cite{KlPl2} that under the standard conditions
the numerical scheme (\ref{4.45}) has strong order of convergence 3.0. 
Among these conditions we consider only the condition
for approximations of iterated Ito stochastic integrals from the numerical
scheme (\ref{4.45}) \cite{KlPl2}, \cite{2006}

\vspace{-1mm}
\begin{equation}
\label{4.3}
{\sf M}\left\{\Biggl(I_{{l_{1}\ldots\hspace{0.2mm} l_{k}}_{\hspace{0.2mm}\tau_{p+1},\tau_p}}
^{(i_{1}\ldots i_{k})} 
-\hat 
I_{{l_{1}\ldots\hspace{0.2mm} l_{k}}_{\hspace{0.2mm}
\tau_{p+1},\tau_p}}^{(i_{1}\ldots i_{k})}
\Biggr)^2\right\}\le C\Delta^{7},
\end{equation}

\vspace{2mm}
\noindent
where 
$\hat 
I_{{l_1\ldots\hspace{0.2mm} l_k}_{\hspace{0.2mm}
\tau_{p+1},\tau_p}}^{(i_1\ldots i_k)}$
is an approximation 
of 
$I_{{l_{1}\ldots\hspace{0.2mm} l_{k}}_{\hspace{0.2mm}
\tau_{p+1},\tau_p}}^{(i_{1}\ldots i_{k})},$
constant $C$ does not depend on $\Delta$.

Note that if we exclude ${\bf u}_{p+1,p}+{\bf v}_{p+1,p}$ from the
right-hand side of (\ref{4.45}), then we will have the explicit 
one-step strong numerical scheme of order 2.0.
The right-hand side of (\ref{4.45}) but without the value
${\bf v}_{p+1,p}$ define 
the explicit 
one-step strong numerical scheme of order 2.5.

Note that the 
truncated 
unified 
Taylor--Ito expansion \cite{2006}, \cite{KK5}, \cite{arxiv-24999},
\cite{kuz1997a}, \cite{2017}-\cite{2013}, \cite{2018a}-\cite{2018aaa}
contains 
the less number of various 
types of iterated Ito stochastic integrals (moreover, their major part 
will have 
less multiplicities) in comparison with 
the classical
Taylor--Ito expansion \cite{KlPl2}, \cite{KlPl1}.

Furthermore, some iterated Ito stochastic integrals from 
the Taylor--Ito expansion \cite{KlPl2}, \cite{KlPl1}
are connected by linear relations. 
However, the iterated stochastic integrals from the 
unified 
Taylor--Ito expansion \cite{2006}, \cite{KK5}, \cite{arxiv-24999},
\cite{kuz1997a}, \cite{2017}-\cite{2013}, \cite{2018a}-\cite{2018aaa}
cannot be connected by linear relations. 
Therefore, we call these families of stochastic integrals 
from the unified Taylor--Ito expansion as the stochastic 
bases \cite{2006}, \cite{arxiv-24999}, \cite{2018a}-\cite{2018aaa}.
Note that (\ref{4.45}) contains 20 different types
of iterated Ito stochastic integrals. At the same time,
the analogue of (\ref{4.45}) based on 
the classical
Taylor--Ito expansion \cite{KlPl2}, \cite{KlPl1} contains 29 different types 
of iterated Ito stochastic integrals.

\vspace{5mm}

\section{Approximation of iterated Ito Stochastic
Integrals Based on Multiple Fourier--Legendre Series}

\vspace{5mm}

Suppose that every $\psi_l(\tau)$ $(l=1,\ldots,k)$ is a 
non-random function from the space $L_2([t, T])$.
Define the following function on the hypercube $[t, T]^k$

\vspace{-1mm}
\begin{equation}
\label{ppp}
K(t_1,\ldots,t_k)=
\begin{cases}
\psi_1(t_1)\ldots \psi_k(t_k)\ &\hbox{for}\ \ t_1<\ldots<t_k\\
~\\
~\\
0\ &\hbox{otherwise}
\end{cases},\ \ \ \ t_1,\ldots,t_k\in[t, T],\ \ \ \ k\ge 2,
\end{equation}

\vspace{2mm}
\noindent
and 
$K(t_1)\equiv\psi_1(t_1)$ for $t_1\in[t, T].$

Suppose that $\{\phi_j(x)\}_{j=0}^{\infty}$
is a complete orthonormal system of functions in the space
$L_2([t, T])$. 
The function $K(t_1,\ldots,t_k)$ belongs to the space
$L_2([t, T]^k).$
At this situation it is well known that the generalized 
multiple Fourier series 
of $K(t_1,\ldots,t_k)\in L_2([t, T]^k)$ is converging 
to $K(t_1,\ldots,t_k)$ in the hypercube $[t, T]^k$ in 
the mean-square sense, i.e.

\vspace{-1mm}
$$
\hbox{\vtop{\offinterlineskip\halign{
\hfil#\hfil\cr
{\rm lim}\cr
$\stackrel{}{{}_{p_1,\ldots,p_k\to \infty}}$\cr
}} }\Biggl\Vert
K(t_1,\ldots,t_k)-
\sum_{j_1=0}^{p_1}\ldots \sum_{j_k=0}^{p_k}
C_{j_k\ldots j_1}\prod_{l=1}^{k} \phi_{j_l}(t_l)
\Biggr\Vert_{L_2([t, T]^k)}=0,
$$

\vspace{2mm}
\noindent
where
\begin{equation}
\label{ppppa}
C_{j_k\ldots j_1}=\int\limits_{[t,T]^k}
K(t_1,\ldots,t_k)\prod_{l=1}^{k}\phi_{j_l}(t_l)dt_1\ldots dt_k,
\end{equation}

$$
\left\Vert f\right\Vert=\left(\int\limits_{[t,T]^k}
f^2(t_1,\ldots,t_k)dt_1\ldots dt_k\right)^{1/2}.
$$

\vspace{3mm}

Consider the partition $\{\tau_j\}_{j=0}^N$ of $[t,T]$ such that

\begin{equation}
\label{1111}
t=\tau_0<\ldots <\tau_N=T,\ \ \
\Delta_N=
\hbox{\vtop{\offinterlineskip\halign{
\hfil#\hfil\cr
{\rm max}\cr
$\stackrel{}{{}_{0\le j\le N-1}}$\cr
}} }\Delta\tau_j\to 0\ \ \hbox{if}\ \ N\to \infty,\ \ \
\Delta\tau_j=\tau_{j+1}-\tau_j.
\end{equation}

\vspace{2mm}

{\bf Theorem 1}\ \cite{2006} (2006), \cite{2017}-\cite{2013}, 
\cite{5-002}-\cite{2018aaa}.
{\it Suppose that
every $\psi_l(\tau)$ $(l=1,\ldots, k)$ is a continuous non-random function on 
$[t, T]$ and
$\{\phi_j(x)\}_{j=0}^{\infty}$ is a complete orthonormal system  
of continuous func\-ti\-ons in the space $L_2([t,T]).$ Then

$$
J[\psi^{(k)}]_{T,t}\  =\ 
\hbox{\vtop{\offinterlineskip\halign{
\hfil#\hfil\cr
{\rm l.i.m.}\cr
$\stackrel{}{{}_{p_1,\ldots,p_k\to \infty}}$\cr
}} }\sum_{j_1=0}^{p_1}\ldots\sum_{j_k=0}^{p_k}
C_{j_k\ldots j_1}\Biggl(
\prod_{l=1}^k\zeta_{j_l}^{(i_l)}\ -
\Biggr.
$$

\vspace{1mm}
\begin{equation}
\label{tyyy}
-\ \Biggl.
\hbox{\vtop{\offinterlineskip\halign{
\hfil#\hfil\cr
{\rm l.i.m.}\cr
$\stackrel{}{{}_{N\to \infty}}$\cr
}} }\sum_{(l_1,\ldots,l_k)\in {\rm G}_k}
\phi_{j_{1}}(\tau_{l_1})
\Delta{\bf w}_{\tau_{l_1}}^{(i_1)}\ldots
\phi_{j_{k}}(\tau_{l_k})
\Delta{\bf w}_{\tau_{l_k}}^{(i_k)}\Biggr),
\end{equation}

\vspace{4mm}
\noindent
where $J[\psi^{(k)}]_{T,t}$ is defined by {\rm (\ref{ito}),}

\vspace{-1mm}
$$
{\rm G}_k={\rm H}_k\backslash{\rm L}_k,\ \ \
{\rm H}_k=\{(l_1,\ldots,l_k):\ l_1,\ldots,l_k=0,\ 1,\ldots,N-1\},
$$

\vspace{-1mm}
$$
{\rm L}_k=\{(l_1,\ldots,l_k):\ l_1,\ldots,l_k=0,\ 1,\ldots,N-1;\
l_g\ne l_r\ (g\ne r);\ g, r=1,\ldots,k\},
$$

\vspace{2mm}
\noindent
${\rm l.i.m.}$ is a limit in the mean-square sense$,$
$i_1,\ldots,i_k=0,1,\ldots,m,$

\vspace{-1mm}
\begin{equation}
\label{rr23}
\zeta_{j}^{(i)}=
\int\limits_t^T \phi_{j}(s) d{\bf w}_s^{(i)}
\end{equation} 

\vspace{2mm}
\noindent
are independent standard Gaussian random variables
for various
$i$ or $j$ {\rm(}if $i\ne 0${\rm),}
$C_{j_k\ldots j_1}$ is the Fourier coefficient {\rm(\ref{ppppa}),}
$\Delta{\bf w}_{\tau_{j}}^{(i)}=
{\bf w}_{\tau_{j+1}}^{(i)}-{\bf w}_{\tau_{j}}^{(i)}$
$(i=0, 1,\ldots,m),$
$\left\{\tau_{j}\right\}_{j=0}^{N}$ is a partition of
the interval $[t, T],$ which satisfies the condition {\rm (\ref{1111})}.
}

\vspace{2mm}

The convergence in the mean of degree 
$2n$ $(n\in \mathbb{N})$ \cite{2017}-\cite{2017-1a}, 
\cite{2007-2}-\cite{2013},
\cite{2018a}-\cite{2018aaa}
as well as the convergence
with probability 1 \cite{arxiv-1},
\cite{arxiv-3}, \cite{2018a}-\cite{2018aaa} 
of approximations from Theorem 1
are proved.
Moreover, the complete orthonormal systems of Haar and 
Rademacher--Walsh functions in $L_2([t,T])$ 
can also be applied in Theorem 1
\cite{2006}, \cite{2017}-\cite{2013}, \cite{2018a}-\cite{2018aaa}.
The modification of Theorem 1 for 
complete orthonormal with weigth $r(x)\ge 0$ systems
of functions in the space $L_2([t,T])$ can be found in 
\cite{2017-1a}, \cite{5-011}-\cite{2018aaa}.
Application of Theorem 1 and Theorem 2 (see below) for the mean-square
approximation of iterated stochastic integrals 
with respect to the 
infinite-dimensional $Q$-Wiener process can be found
in the monographs \cite{2018a}-\cite{2018aaa} 
(Chapter 7) and in \cite{5-007}, \cite{5-008}, \cite{Kuzh-1}, 
\cite{new-new-1}-\cite{arxiv-21}.

In order to evaluate the significance of Theorem 1 for practice we will
demonstrate its transformed particular cases for 
$k=1,\ldots,6$ 
\cite{2006}, \cite{2017}-\cite{2013}, 
\cite{5-002}-\cite{2018aaa}

\begin{equation}
\label{a1}
J[\psi^{(1)}]_{T,t}
=\hbox{\vtop{\offinterlineskip\halign{
\hfil#\hfil\cr
{\rm l.i.m.}\cr
$\stackrel{}{{}_{p_1\to \infty}}$\cr
}} }\sum_{j_1=0}^{p_1}
C_{j_1}\zeta_{j_1}^{(i_1)},
\end{equation}

\vspace{2mm}
\begin{equation}
\label{a2}
J[\psi^{(2)}]_{T,t}
=\hbox{\vtop{\offinterlineskip\halign{
\hfil#\hfil\cr
{\rm l.i.m.}\cr
$\stackrel{}{{}_{p_1,p_2\to \infty}}$\cr
}} }\sum_{j_1=0}^{p_1}\sum_{j_2=0}^{p_2}
C_{j_2j_1}\Biggl(\zeta_{j_1}^{(i_1)}\zeta_{j_2}^{(i_2)}
-{\bf 1}_{\{i_1=i_2\ne 0\}}
{\bf 1}_{\{j_1=j_2\}}\Biggr),
\end{equation}

\vspace{6mm}

$$
J[\psi^{(3)}]_{T,t}=
\hbox{\vtop{\offinterlineskip\halign{
\hfil#\hfil\cr
{\rm l.i.m.}\cr
$\stackrel{}{{}_{p_1,\ldots,p_3\to \infty}}$\cr
}} }\sum_{j_1=0}^{p_1}\sum_{j_2=0}^{p_2}\sum_{j_3=0}^{p_3}
C_{j_3j_2j_1}\Biggl(
\zeta_{j_1}^{(i_1)}\zeta_{j_2}^{(i_2)}\zeta_{j_3}^{(i_3)}
-\Biggr.
$$
\begin{equation}
\label{a3}
-\Biggl.
{\bf 1}_{\{i_1=i_2\ne 0\}}
{\bf 1}_{\{j_1=j_2\}}
\zeta_{j_3}^{(i_3)}
-{\bf 1}_{\{i_2=i_3\ne 0\}}
{\bf 1}_{\{j_2=j_3\}}
\zeta_{j_1}^{(i_1)}-
{\bf 1}_{\{i_1=i_3\ne 0\}}
{\bf 1}_{\{j_1=j_3\}}
\zeta_{j_2}^{(i_2)}\Biggr),
\end{equation}

\vspace{6mm}

$$
J[\psi^{(4)}]_{T,t}
=
\hbox{\vtop{\offinterlineskip\halign{
\hfil#\hfil\cr
{\rm l.i.m.}\cr
$\stackrel{}{{}_{p_1,\ldots,p_4\to \infty}}$\cr
}} }\sum_{j_1=0}^{p_1}\ldots\sum_{j_4=0}^{p_4}
C_{j_4\ldots j_1}\Biggl(
\prod_{l=1}^4\zeta_{j_l}^{(i_l)}
\Biggr.
-
$$
$$
-
{\bf 1}_{\{i_1=i_2\ne 0\}}
{\bf 1}_{\{j_1=j_2\}}
\zeta_{j_3}^{(i_3)}
\zeta_{j_4}^{(i_4)}
-
{\bf 1}_{\{i_1=i_3\ne 0\}}
{\bf 1}_{\{j_1=j_3\}}
\zeta_{j_2}^{(i_2)}
\zeta_{j_4}^{(i_4)}-
$$
$$
-
{\bf 1}_{\{i_1=i_4\ne 0\}}
{\bf 1}_{\{j_1=j_4\}}
\zeta_{j_2}^{(i_2)}
\zeta_{j_3}^{(i_3)}
-
{\bf 1}_{\{i_2=i_3\ne 0\}}
{\bf 1}_{\{j_2=j_3\}}
\zeta_{j_1}^{(i_1)}
\zeta_{j_4}^{(i_4)}-
$$
$$
-
{\bf 1}_{\{i_2=i_4\ne 0\}}
{\bf 1}_{\{j_2=j_4\}}
\zeta_{j_1}^{(i_1)}
\zeta_{j_3}^{(i_3)}
-
{\bf 1}_{\{i_3=i_4\ne 0\}}
{\bf 1}_{\{j_3=j_4\}}
\zeta_{j_1}^{(i_1)}
\zeta_{j_2}^{(i_2)}+
$$
$$
+
{\bf 1}_{\{i_1=i_2\ne 0\}}
{\bf 1}_{\{j_1=j_2\}}
{\bf 1}_{\{i_3=i_4\ne 0\}}
{\bf 1}_{\{j_3=j_4\}}
+
{\bf 1}_{\{i_1=i_3\ne 0\}}
{\bf 1}_{\{j_1=j_3\}}
{\bf 1}_{\{i_2=i_4\ne 0\}}
{\bf 1}_{\{j_2=j_4\}}+
$$
\begin{equation}
\label{a4}
+\Biggl.
{\bf 1}_{\{i_1=i_4\ne 0\}}
{\bf 1}_{\{j_1=j_4\}}
{\bf 1}_{\{i_2=i_3\ne 0\}}
{\bf 1}_{\{j_2=j_3\}}\Biggr),
\end{equation}

\vspace{8mm}

$$
J[\psi^{(5)}]_{T,t}
=\hbox{\vtop{\offinterlineskip\halign{
\hfil#\hfil\cr
{\rm l.i.m.}\cr
$\stackrel{}{{}_{p_1,\ldots,p_5\to \infty}}$\cr
}} }\sum_{j_1=0}^{p_1}\ldots\sum_{j_5=0}^{p_5}
C_{j_5\ldots j_1}\Biggl(
\prod_{l=1}^5\zeta_{j_l}^{(i_l)}
-\Biggr.
$$
$$
-
{\bf 1}_{\{i_1=i_2\ne 0\}}
{\bf 1}_{\{j_1=j_2\}}
\zeta_{j_3}^{(i_3)}
\zeta_{j_4}^{(i_4)}
\zeta_{j_5}^{(i_5)}-
{\bf 1}_{\{i_1=i_3\ne 0\}}
{\bf 1}_{\{j_1=j_3\}}
\zeta_{j_2}^{(i_2)}
\zeta_{j_4}^{(i_4)}
\zeta_{j_5}^{(i_5)}-
$$
$$
-
{\bf 1}_{\{i_1=i_4\ne 0\}}
{\bf 1}_{\{j_1=j_4\}}
\zeta_{j_2}^{(i_2)}
\zeta_{j_3}^{(i_3)}
\zeta_{j_5}^{(i_5)}-
{\bf 1}_{\{i_1=i_5\ne 0\}}
{\bf 1}_{\{j_1=j_5\}}
\zeta_{j_2}^{(i_2)}
\zeta_{j_3}^{(i_3)}
\zeta_{j_4}^{(i_4)}-
$$
$$
-
{\bf 1}_{\{i_2=i_3\ne 0\}}
{\bf 1}_{\{j_2=j_3\}}
\zeta_{j_1}^{(i_1)}
\zeta_{j_4}^{(i_4)}
\zeta_{j_5}^{(i_5)}-
{\bf 1}_{\{i_2=i_4\ne 0\}}
{\bf 1}_{\{j_2=j_4\}}
\zeta_{j_1}^{(i_1)}
\zeta_{j_3}^{(i_3)}
\zeta_{j_5}^{(i_5)}-
$$
$$
-
{\bf 1}_{\{i_2=i_5\ne 0\}}
{\bf 1}_{\{j_2=j_5\}}
\zeta_{j_1}^{(i_1)}
\zeta_{j_3}^{(i_3)}
\zeta_{j_4}^{(i_4)}
-{\bf 1}_{\{i_3=i_4\ne 0\}}
{\bf 1}_{\{j_3=j_4\}}
\zeta_{j_1}^{(i_1)}
\zeta_{j_2}^{(i_2)}
\zeta_{j_5}^{(i_5)}-
$$
$$
-
{\bf 1}_{\{i_3=i_5\ne 0\}}
{\bf 1}_{\{j_3=j_5\}}
\zeta_{j_1}^{(i_1)}
\zeta_{j_2}^{(i_2)}
\zeta_{j_4}^{(i_4)}
-{\bf 1}_{\{i_4=i_5\ne 0\}}
{\bf 1}_{\{j_4=j_5\}}
\zeta_{j_1}^{(i_1)}
\zeta_{j_2}^{(i_2)}
\zeta_{j_3}^{(i_3)}+
$$
$$
+
{\bf 1}_{\{i_1=i_2\ne 0\}}
{\bf 1}_{\{j_1=j_2\}}
{\bf 1}_{\{i_3=i_4\ne 0\}}
{\bf 1}_{\{j_3=j_4\}}\zeta_{j_5}^{(i_5)}+
{\bf 1}_{\{i_1=i_2\ne 0\}}
{\bf 1}_{\{j_1=j_2\}}
{\bf 1}_{\{i_3=i_5\ne 0\}}
{\bf 1}_{\{j_3=j_5\}}\zeta_{j_4}^{(i_4)}+
$$
$$
+
{\bf 1}_{\{i_1=i_2\ne 0\}}
{\bf 1}_{\{j_1=j_2\}}
{\bf 1}_{\{i_4=i_5\ne 0\}}
{\bf 1}_{\{j_4=j_5\}}\zeta_{j_3}^{(i_3)}+
{\bf 1}_{\{i_1=i_3\ne 0\}}
{\bf 1}_{\{j_1=j_3\}}
{\bf 1}_{\{i_2=i_4\ne 0\}}
{\bf 1}_{\{j_2=j_4\}}\zeta_{j_5}^{(i_5)}+
$$
$$
+
{\bf 1}_{\{i_1=i_3\ne 0\}}
{\bf 1}_{\{j_1=j_3\}}
{\bf 1}_{\{i_2=i_5\ne 0\}}
{\bf 1}_{\{j_2=j_5\}}\zeta_{j_4}^{(i_4)}+
{\bf 1}_{\{i_1=i_3\ne 0\}}
{\bf 1}_{\{j_1=j_3\}}
{\bf 1}_{\{i_4=i_5\ne 0\}}
{\bf 1}_{\{j_4=j_5\}}\zeta_{j_2}^{(i_2)}+
$$
$$
+
{\bf 1}_{\{i_1=i_4\ne 0\}}
{\bf 1}_{\{j_1=j_4\}}
{\bf 1}_{\{i_2=i_3\ne 0\}}
{\bf 1}_{\{j_2=j_3\}}\zeta_{j_5}^{(i_5)}+
{\bf 1}_{\{i_1=i_4\ne 0\}}
{\bf 1}_{\{j_1=j_4\}}
{\bf 1}_{\{i_2=i_5\ne 0\}}
{\bf 1}_{\{j_2=j_5\}}\zeta_{j_3}^{(i_3)}+
$$
$$
+
{\bf 1}_{\{i_1=i_4\ne 0\}}
{\bf 1}_{\{j_1=j_4\}}
{\bf 1}_{\{i_3=i_5\ne 0\}}
{\bf 1}_{\{j_3=j_5\}}\zeta_{j_2}^{(i_2)}+
{\bf 1}_{\{i_1=i_5\ne 0\}}
{\bf 1}_{\{j_1=j_5\}}
{\bf 1}_{\{i_2=i_3\ne 0\}}
{\bf 1}_{\{j_2=j_3\}}\zeta_{j_4}^{(i_4)}+
$$
$$
+
{\bf 1}_{\{i_1=i_5\ne 0\}}
{\bf 1}_{\{j_1=j_5\}}
{\bf 1}_{\{i_2=i_4\ne 0\}}
{\bf 1}_{\{j_2=j_4\}}\zeta_{j_3}^{(i_3)}+
{\bf 1}_{\{i_1=i_5\ne 0\}}
{\bf 1}_{\{j_1=j_5\}}
{\bf 1}_{\{i_3=i_4\ne 0\}}
{\bf 1}_{\{j_3=j_4\}}\zeta_{j_2}^{(i_2)}+
$$
$$
+
{\bf 1}_{\{i_2=i_3\ne 0\}}
{\bf 1}_{\{j_2=j_3\}}
{\bf 1}_{\{i_4=i_5\ne 0\}}
{\bf 1}_{\{j_4=j_5\}}\zeta_{j_1}^{(i_1)}+
{\bf 1}_{\{i_2=i_4\ne 0\}}
{\bf 1}_{\{j_2=j_4\}}
{\bf 1}_{\{i_3=i_5\ne 0\}}
{\bf 1}_{\{j_3=j_5\}}\zeta_{j_1}^{(i_1)}+
$$
\begin{equation}
\label{a5}
+\Biggl.
{\bf 1}_{\{i_2=i_5\ne 0\}}
{\bf 1}_{\{j_2=j_5\}}
{\bf 1}_{\{i_3=i_4\ne 0\}}
{\bf 1}_{\{j_3=j_4\}}\zeta_{j_1}^{(i_1)}\Biggr),
\end{equation}

\vspace{9mm}

$$
J[\psi^{(6)}]_{T,t}
=\hbox{\vtop{\offinterlineskip\halign{
\hfil#\hfil\cr
{\rm l.i.m.}\cr
$\stackrel{}{{}_{p_1,\ldots,p_6\to \infty}}$\cr
}} }\sum_{j_1=0}^{p_1}\ldots\sum_{j_6=0}^{p_6}
C_{j_6\ldots j_1}\Biggl(
\prod_{l=1}^6
\zeta_{j_l}^{(i_l)}
-\Biggr.
$$
$$
-
{\bf 1}_{\{i_1=i_6\ne 0\}}
{\bf 1}_{\{j_1=j_6\}}
\zeta_{j_2}^{(i_2)}
\zeta_{j_3}^{(i_3)}
\zeta_{j_4}^{(i_4)}
\zeta_{j_5}^{(i_5)}-
{\bf 1}_{\{i_2=i_6\ne 0\}}
{\bf 1}_{\{j_2=j_6\}}
\zeta_{j_1}^{(i_1)}
\zeta_{j_3}^{(i_3)}
\zeta_{j_4}^{(i_4)}
\zeta_{j_5}^{(i_5)}-
$$
$$
-
{\bf 1}_{\{i_3=i_6\ne 0\}}
{\bf 1}_{\{j_3=j_6\}}
\zeta_{j_1}^{(i_1)}
\zeta_{j_2}^{(i_2)}
\zeta_{j_4}^{(i_4)}
\zeta_{j_5}^{(i_5)}-
{\bf 1}_{\{i_4=i_6\ne 0\}}
{\bf 1}_{\{j_4=j_6\}}
\zeta_{j_1}^{(i_1)}
\zeta_{j_2}^{(i_2)}
\zeta_{j_3}^{(i_3)}
\zeta_{j_5}^{(i_5)}-
$$
$$
-
{\bf 1}_{\{i_5=i_6\ne 0\}}
{\bf 1}_{\{j_5=j_6\}}
\zeta_{j_1}^{(i_1)}
\zeta_{j_2}^{(i_2)}
\zeta_{j_3}^{(i_3)}
\zeta_{j_4}^{(i_4)}-
{\bf 1}_{\{i_1=i_2\ne 0\}}
{\bf 1}_{\{j_1=j_2\}}
\zeta_{j_3}^{(i_3)}
\zeta_{j_4}^{(i_4)}
\zeta_{j_5}^{(i_5)}
\zeta_{j_6}^{(i_6)}-
$$
$$
-
{\bf 1}_{\{i_1=i_3\ne 0\}}
{\bf 1}_{\{j_1=j_3\}}
\zeta_{j_2}^{(i_2)}
\zeta_{j_4}^{(i_4)}
\zeta_{j_5}^{(i_5)}
\zeta_{j_6}^{(i_6)}-
{\bf 1}_{\{i_1=i_4\ne 0\}}
{\bf 1}_{\{j_1=j_4\}}
\zeta_{j_2}^{(i_2)}
\zeta_{j_3}^{(i_3)}
\zeta_{j_5}^{(i_5)}
\zeta_{j_6}^{(i_6)}-
$$
$$
-
{\bf 1}_{\{i_1=i_5\ne 0\}}
{\bf 1}_{\{j_1=j_5\}}
\zeta_{j_2}^{(i_2)}
\zeta_{j_3}^{(i_3)}
\zeta_{j_4}^{(i_4)}
\zeta_{j_6}^{(i_6)}-
{\bf 1}_{\{i_2=i_3\ne 0\}}
{\bf 1}_{\{j_2=j_3\}}
\zeta_{j_1}^{(i_1)}
\zeta_{j_4}^{(i_4)}
\zeta_{j_5}^{(i_5)}
\zeta_{j_6}^{(i_6)}-
$$
$$
-
{\bf 1}_{\{i_2=i_4\ne 0\}}
{\bf 1}_{\{j_2=j_4\}}
\zeta_{j_1}^{(i_1)}
\zeta_{j_3}^{(i_3)}
\zeta_{j_5}^{(i_5)}
\zeta_{j_6}^{(i_6)}-
{\bf 1}_{\{i_2=i_5\ne 0\}}
{\bf 1}_{\{j_2=j_5\}}
\zeta_{j_1}^{(i_1)}
\zeta_{j_3}^{(i_3)}
\zeta_{j_4}^{(i_4)}
\zeta_{j_6}^{(i_6)}-
$$
$$
-
{\bf 1}_{\{i_3=i_4\ne 0\}}
{\bf 1}_{\{j_3=j_4\}}
\zeta_{j_1}^{(i_1)}
\zeta_{j_2}^{(i_2)}
\zeta_{j_5}^{(i_5)}
\zeta_{j_6}^{(i_6)}-
{\bf 1}_{\{i_3=i_5\ne 0\}}
{\bf 1}_{\{j_3=j_5\}}
\zeta_{j_1}^{(i_1)}
\zeta_{j_2}^{(i_2)}
\zeta_{j_4}^{(i_4)}
\zeta_{j_6}^{(i_6)}-
$$
$$
-
{\bf 1}_{\{i_4=i_5\ne 0\}}
{\bf 1}_{\{j_4=j_5\}}
\zeta_{j_1}^{(i_1)}
\zeta_{j_2}^{(i_2)}
\zeta_{j_3}^{(i_3)}
\zeta_{j_6}^{(i_6)}+
$$
$$
+
{\bf 1}_{\{i_1=i_2\ne 0\}}
{\bf 1}_{\{j_1=j_2\}}
{\bf 1}_{\{i_3=i_4\ne 0\}}
{\bf 1}_{\{j_3=j_4\}}
\zeta_{j_5}^{(i_5)}
\zeta_{j_6}^{(i_6)}+
{\bf 1}_{\{i_1=i_2\ne 0\}}
{\bf 1}_{\{j_1=j_2\}}
{\bf 1}_{\{i_3=i_5\ne 0\}}
{\bf 1}_{\{j_3=j_5\}}
\zeta_{j_4}^{(i_4)}
\zeta_{j_6}^{(i_6)}+
$$
$$
+
{\bf 1}_{\{i_1=i_2\ne 0\}}
{\bf 1}_{\{j_1=j_2\}}
{\bf 1}_{\{i_4=i_5\ne 0\}}
{\bf 1}_{\{j_4=j_5\}}
\zeta_{j_3}^{(i_3)}
\zeta_{j_6}^{(i_6)}
+
{\bf 1}_{\{i_1=i_3\ne 0\}}
{\bf 1}_{\{j_1=j_3\}}
{\bf 1}_{\{i_2=i_4\ne 0\}}
{\bf 1}_{\{j_2=j_4\}}
\zeta_{j_5}^{(i_5)}
\zeta_{j_6}^{(i_6)}+
$$
$$
+
{\bf 1}_{\{i_1=i_3\ne 0\}}
{\bf 1}_{\{j_1=j_3\}}
{\bf 1}_{\{i_2=i_5\ne 0\}}
{\bf 1}_{\{j_2=j_5\}}
\zeta_{j_4}^{(i_4)}
\zeta_{j_6}^{(i_6)}
+{\bf 1}_{\{i_1=i_3\ne 0\}}
{\bf 1}_{\{j_1=j_3\}}
{\bf 1}_{\{i_4=i_5\ne 0\}}
{\bf 1}_{\{j_4=j_5\}}
\zeta_{j_2}^{(i_2)}
\zeta_{j_6}^{(i_6)}+
$$
$$
+
{\bf 1}_{\{i_1=i_4\ne 0\}}
{\bf 1}_{\{j_1=j_4\}}
{\bf 1}_{\{i_2=i_3\ne 0\}}
{\bf 1}_{\{j_2=j_3\}}
\zeta_{j_5}^{(i_5)}
\zeta_{j_6}^{(i_6)}
+
{\bf 1}_{\{i_1=i_4\ne 0\}}
{\bf 1}_{\{j_1=j_4\}}
{\bf 1}_{\{i_2=i_5\ne 0\}}
{\bf 1}_{\{j_2=j_5\}}
\zeta_{j_3}^{(i_3)}
\zeta_{j_6}^{(i_6)}+
$$
$$
+
{\bf 1}_{\{i_1=i_4\ne 0\}}
{\bf 1}_{\{j_1=j_4\}}
{\bf 1}_{\{i_3=i_5\ne 0\}}
{\bf 1}_{\{j_3=j_5\}}
\zeta_{j_2}^{(i_2)}
\zeta_{j_6}^{(i_6)}
+
{\bf 1}_{\{i_1=i_5\ne 0\}}
{\bf 1}_{\{j_1=j_5\}}
{\bf 1}_{\{i_2=i_3\ne 0\}}
{\bf 1}_{\{j_2=j_3\}}
\zeta_{j_4}^{(i_4)}
\zeta_{j_6}^{(i_6)}+
$$
$$
+
{\bf 1}_{\{i_1=i_5\ne 0\}}
{\bf 1}_{\{j_1=j_5\}}
{\bf 1}_{\{i_2=i_4\ne 0\}}
{\bf 1}_{\{j_2=j_4\}}
\zeta_{j_3}^{(i_3)}
\zeta_{j_6}^{(i_6)}
+
{\bf 1}_{\{i_1=i_5\ne 0\}}
{\bf 1}_{\{j_1=j_5\}}
{\bf 1}_{\{i_3=i_4\ne 0\}}
{\bf 1}_{\{j_3=j_4\}}
\zeta_{j_2}^{(i_2)}
\zeta_{j_6}^{(i_6)}+
$$
$$
+
{\bf 1}_{\{i_2=i_3\ne 0\}}
{\bf 1}_{\{j_2=j_3\}}
{\bf 1}_{\{i_4=i_5\ne 0\}}
{\bf 1}_{\{j_4=j_5\}}
\zeta_{j_1}^{(i_1)}
\zeta_{j_6}^{(i_6)}
+
{\bf 1}_{\{i_2=i_4\ne 0\}}
{\bf 1}_{\{j_2=j_4\}}
{\bf 1}_{\{i_3=i_5\ne 0\}}
{\bf 1}_{\{j_3=j_5\}}
\zeta_{j_1}^{(i_1)}
\zeta_{j_6}^{(i_6)}+
$$
$$
+
{\bf 1}_{\{i_2=i_5\ne 0\}}
{\bf 1}_{\{j_2=j_5\}}
{\bf 1}_{\{i_3=i_4\ne 0\}}
{\bf 1}_{\{j_3=j_4\}}
\zeta_{j_1}^{(i_1)}
\zeta_{j_6}^{(i_6)}
+
{\bf 1}_{\{i_6=i_1\ne 0\}}
{\bf 1}_{\{j_6=j_1\}}
{\bf 1}_{\{i_3=i_4\ne 0\}}
{\bf 1}_{\{j_3=j_4\}}
\zeta_{j_2}^{(i_2)}
\zeta_{j_5}^{(i_5)}+
$$
$$
+
{\bf 1}_{\{i_6=i_1\ne 0\}}
{\bf 1}_{\{j_6=j_1\}}
{\bf 1}_{\{i_3=i_5\ne 0\}}
{\bf 1}_{\{j_3=j_5\}}
\zeta_{j_2}^{(i_2)}
\zeta_{j_4}^{(i_4)}
+
{\bf 1}_{\{i_6=i_1\ne 0\}}
{\bf 1}_{\{j_6=j_1\}}
{\bf 1}_{\{i_2=i_5\ne 0\}}
{\bf 1}_{\{j_2=j_5\}}
\zeta_{j_3}^{(i_3)}
\zeta_{j_4}^{(i_4)}+
$$
$$
+
{\bf 1}_{\{i_6=i_1\ne 0\}}
{\bf 1}_{\{j_6=j_1\}}
{\bf 1}_{\{i_2=i_4\ne 0\}}
{\bf 1}_{\{j_2=j_4\}}
\zeta_{j_3}^{(i_3)}
\zeta_{j_5}^{(i_5)}
+
{\bf 1}_{\{i_6=i_1\ne 0\}}
{\bf 1}_{\{j_6=j_1\}}
{\bf 1}_{\{i_4=i_5\ne 0\}}
{\bf 1}_{\{j_4=j_5\}}
\zeta_{j_2}^{(i_2)}
\zeta_{j_3}^{(i_3)}+
$$
$$
+
{\bf 1}_{\{i_6=i_1\ne 0\}}
{\bf 1}_{\{j_6=j_1\}}
{\bf 1}_{\{i_2=i_3\ne 0\}}
{\bf 1}_{\{j_2=j_3\}}
\zeta_{j_4}^{(i_4)}
\zeta_{j_5}^{(i_5)}
+
{\bf 1}_{\{i_6=i_2\ne 0\}}
{\bf 1}_{\{j_6=j_2\}}
{\bf 1}_{\{i_3=i_5\ne 0\}}
{\bf 1}_{\{j_3=j_5\}}
\zeta_{j_1}^{(i_1)}
\zeta_{j_4}^{(i_4)}+
$$
$$
+
{\bf 1}_{\{i_6=i_2\ne 0\}}
{\bf 1}_{\{j_6=j_2\}}
{\bf 1}_{\{i_4=i_5\ne 0\}}
{\bf 1}_{\{j_4=j_5\}}
\zeta_{j_1}^{(i_1)}
\zeta_{j_3}^{(i_3)}
+
{\bf 1}_{\{i_6=i_2\ne 0\}}
{\bf 1}_{\{j_6=j_2\}}
{\bf 1}_{\{i_3=i_4\ne 0\}}
{\bf 1}_{\{j_3=j_4\}}
\zeta_{j_1}^{(i_1)}
\zeta_{j_5}^{(i_5)}+
$$
$$
+
{\bf 1}_{\{i_6=i_2\ne 0\}}
{\bf 1}_{\{j_6=j_2\}}
{\bf 1}_{\{i_1=i_5\ne 0\}}
{\bf 1}_{\{j_1=j_5\}}
\zeta_{j_3}^{(i_3)}
\zeta_{j_4}^{(i_4)}
+
{\bf 1}_{\{i_6=i_2\ne 0\}}
{\bf 1}_{\{j_6=j_2\}}
{\bf 1}_{\{i_1=i_4\ne 0\}}
{\bf 1}_{\{j_1=j_4\}}
\zeta_{j_3}^{(i_3)}
\zeta_{j_5}^{(i_5)}+
$$
$$
+
{\bf 1}_{\{i_6=i_2\ne 0\}}
{\bf 1}_{\{j_6=j_2\}}
{\bf 1}_{\{i_1=i_3\ne 0\}}
{\bf 1}_{\{j_1=j_3\}}
\zeta_{j_4}^{(i_4)}
\zeta_{j_5}^{(i_5)}
+
{\bf 1}_{\{i_6=i_3\ne 0\}}
{\bf 1}_{\{j_6=j_3\}}
{\bf 1}_{\{i_2=i_5\ne 0\}}
{\bf 1}_{\{j_2=j_5\}}
\zeta_{j_1}^{(i_1)}
\zeta_{j_4}^{(i_4)}+
$$
$$
+
{\bf 1}_{\{i_6=i_3\ne 0\}}
{\bf 1}_{\{j_6=j_3\}}
{\bf 1}_{\{i_4=i_5\ne 0\}}
{\bf 1}_{\{j_4=j_5\}}
\zeta_{j_1}^{(i_1)}
\zeta_{j_2}^{(i_2)}
+
{\bf 1}_{\{i_6=i_3\ne 0\}}
{\bf 1}_{\{j_6=j_3\}}
{\bf 1}_{\{i_2=i_4\ne 0\}}
{\bf 1}_{\{j_2=j_4\}}
\zeta_{j_1}^{(i_1)}
\zeta_{j_5}^{(i_5)}+
$$
$$
+
{\bf 1}_{\{i_6=i_3\ne 0\}}
{\bf 1}_{\{j_6=j_3\}}
{\bf 1}_{\{i_1=i_5\ne 0\}}
{\bf 1}_{\{j_1=j_5\}}
\zeta_{j_2}^{(i_2)}
\zeta_{j_4}^{(i_4)}
+
{\bf 1}_{\{i_6=i_3\ne 0\}}
{\bf 1}_{\{j_6=j_3\}}
{\bf 1}_{\{i_1=i_4\ne 0\}}
{\bf 1}_{\{j_1=j_4\}}
\zeta_{j_2}^{(i_2)}
\zeta_{j_5}^{(i_5)}+
$$
$$
+
{\bf 1}_{\{i_6=i_3\ne 0\}}
{\bf 1}_{\{j_6=j_3\}}
{\bf 1}_{\{i_1=i_2\ne 0\}}
{\bf 1}_{\{j_1=j_2\}}
\zeta_{j_4}^{(i_4)}
\zeta_{j_5}^{(i_5)}
+
{\bf 1}_{\{i_6=i_4\ne 0\}}
{\bf 1}_{\{j_6=j_4\}}
{\bf 1}_{\{i_3=i_5\ne 0\}}
{\bf 1}_{\{j_3=j_5\}}
\zeta_{j_1}^{(i_1)}
\zeta_{j_2}^{(i_2)}+
$$
$$
+
{\bf 1}_{\{i_6=i_4\ne 0\}}
{\bf 1}_{\{j_6=j_4\}}
{\bf 1}_{\{i_2=i_5\ne 0\}}
{\bf 1}_{\{j_2=j_5\}}
\zeta_{j_1}^{(i_1)}
\zeta_{j_3}^{(i_3)}
+
{\bf 1}_{\{i_6=i_4\ne 0\}}
{\bf 1}_{\{j_6=j_4\}}
{\bf 1}_{\{i_2=i_3\ne 0\}}
{\bf 1}_{\{j_2=j_3\}}
\zeta_{j_1}^{(i_1)}
\zeta_{j_5}^{(i_5)}+
$$
$$
+
{\bf 1}_{\{i_6=i_4\ne 0\}}
{\bf 1}_{\{j_6=j_4\}}
{\bf 1}_{\{i_1=i_5\ne 0\}}
{\bf 1}_{\{j_1=j_5\}}
\zeta_{j_2}^{(i_2)}
\zeta_{j_3}^{(i_3)}
+
{\bf 1}_{\{i_6=i_4\ne 0\}}
{\bf 1}_{\{j_6=j_4\}}
{\bf 1}_{\{i_1=i_3\ne 0\}}
{\bf 1}_{\{j_1=j_3\}}
\zeta_{j_2}^{(i_2)}
\zeta_{j_5}^{(i_5)}+
$$
$$
+
{\bf 1}_{\{i_6=i_4\ne 0\}}
{\bf 1}_{\{j_6=j_4\}}
{\bf 1}_{\{i_1=i_2\ne 0\}}
{\bf 1}_{\{j_1=j_2\}}
\zeta_{j_3}^{(i_3)}
\zeta_{j_5}^{(i_5)}
+
{\bf 1}_{\{i_6=i_5\ne 0\}}
{\bf 1}_{\{j_6=j_5\}}
{\bf 1}_{\{i_3=i_4\ne 0\}}
{\bf 1}_{\{j_3=j_4\}}
\zeta_{j_1}^{(i_1)}
\zeta_{j_2}^{(i_2)}+
$$
$$
+
{\bf 1}_{\{i_6=i_5\ne 0\}}
{\bf 1}_{\{j_6=j_5\}}
{\bf 1}_{\{i_2=i_4\ne 0\}}
{\bf 1}_{\{j_2=j_4\}}
\zeta_{j_1}^{(i_1)}
\zeta_{j_3}^{(i_3)}
+
{\bf 1}_{\{i_6=i_5\ne 0\}}
{\bf 1}_{\{j_6=j_5\}}
{\bf 1}_{\{i_2=i_3\ne 0\}}
{\bf 1}_{\{j_2=j_3\}}
\zeta_{j_1}^{(i_1)}
\zeta_{j_4}^{(i_4)}+
$$
$$
+
{\bf 1}_{\{i_6=i_5\ne 0\}}
{\bf 1}_{\{j_6=j_5\}}
{\bf 1}_{\{i_1=i_4\ne 0\}}
{\bf 1}_{\{j_1=j_4\}}
\zeta_{j_2}^{(i_2)}
\zeta_{j_3}^{(i_3)}
+
{\bf 1}_{\{i_6=i_5\ne 0\}}
{\bf 1}_{\{j_6=j_5\}}
{\bf 1}_{\{i_1=i_3\ne 0\}}
{\bf 1}_{\{j_1=j_3\}}
\zeta_{j_2}^{(i_2)}
\zeta_{j_4}^{(i_4)}+
$$
$$
+
{\bf 1}_{\{i_6=i_5\ne 0\}}
{\bf 1}_{\{j_6=j_5\}}
{\bf 1}_{\{i_1=i_2\ne 0\}}
{\bf 1}_{\{j_1=j_2\}}
\zeta_{j_3}^{(i_3)}
\zeta_{j_4}^{(i_4)}-
$$
$$
-
{\bf 1}_{\{i_6=i_1\ne 0\}}
{\bf 1}_{\{j_6=j_1\}}
{\bf 1}_{\{i_2=i_5\ne 0\}}
{\bf 1}_{\{j_2=j_5\}}
{\bf 1}_{\{i_3=i_4\ne 0\}}
{\bf 1}_{\{j_3=j_4\}}-
$$
$$
-
{\bf 1}_{\{i_6=i_1\ne 0\}}
{\bf 1}_{\{j_6=j_1\}}
{\bf 1}_{\{i_2=i_4\ne 0\}}
{\bf 1}_{\{j_2=j_4\}}
{\bf 1}_{\{i_3=i_5\ne 0\}}
{\bf 1}_{\{j_3=j_5\}}-
$$
$$
-
{\bf 1}_{\{i_6=i_1\ne 0\}}
{\bf 1}_{\{j_6=j_1\}}
{\bf 1}_{\{i_2=i_3\ne 0\}}
{\bf 1}_{\{j_2=j_3\}}
{\bf 1}_{\{i_4=i_5\ne 0\}}
{\bf 1}_{\{j_4=j_5\}}-
$$
$$
-
{\bf 1}_{\{i_6=i_2\ne 0\}}
{\bf 1}_{\{j_6=j_2\}}
{\bf 1}_{\{i_1=i_5\ne 0\}}
{\bf 1}_{\{j_1=j_5\}}
{\bf 1}_{\{i_3=i_4\ne 0\}}
{\bf 1}_{\{j_3=j_4\}}-
$$
$$
-
{\bf 1}_{\{i_6=i_2\ne 0\}}
{\bf 1}_{\{j_6=j_2\}}
{\bf 1}_{\{i_1=i_4\ne 0\}}
{\bf 1}_{\{j_1=j_4\}}
{\bf 1}_{\{i_3=i_5\ne 0\}}
{\bf 1}_{\{j_3=j_5\}}-
$$
$$
-
{\bf 1}_{\{i_6=i_2\ne 0\}}
{\bf 1}_{\{j_6=j_2\}}
{\bf 1}_{\{i_1=i_3\ne 0\}}
{\bf 1}_{\{j_1=j_3\}}
{\bf 1}_{\{i_4=i_5\ne 0\}}
{\bf 1}_{\{j_4=j_5\}}-
$$
$$
-
{\bf 1}_{\{i_6=i_3\ne 0\}}
{\bf 1}_{\{j_6=j_3\}}
{\bf 1}_{\{i_1=i_5\ne 0\}}
{\bf 1}_{\{j_1=j_5\}}
{\bf 1}_{\{i_2=i_4\ne 0\}}
{\bf 1}_{\{j_2=j_4\}}-
$$
$$
-
{\bf 1}_{\{i_6=i_3\ne 0\}}
{\bf 1}_{\{j_6=j_3\}}
{\bf 1}_{\{i_1=i_4\ne 0\}}
{\bf 1}_{\{j_1=j_4\}}
{\bf 1}_{\{i_2=i_5\ne 0\}}
{\bf 1}_{\{j_2=j_5\}}-
$$
$$
-
{\bf 1}_{\{i_3=i_6\ne 0\}}
{\bf 1}_{\{j_3=j_6\}}
{\bf 1}_{\{i_1=i_2\ne 0\}}
{\bf 1}_{\{j_1=j_2\}}
{\bf 1}_{\{i_4=i_5\ne 0\}}
{\bf 1}_{\{j_4=j_5\}}-
$$
$$
-
{\bf 1}_{\{i_6=i_4\ne 0\}}
{\bf 1}_{\{j_6=j_4\}}
{\bf 1}_{\{i_1=i_5\ne 0\}}
{\bf 1}_{\{j_1=j_5\}}
{\bf 1}_{\{i_2=i_3\ne 0\}}
{\bf 1}_{\{j_2=j_3\}}-
$$
$$
-
{\bf 1}_{\{i_6=i_4\ne 0\}}
{\bf 1}_{\{j_6=j_4\}}
{\bf 1}_{\{i_1=i_3\ne 0\}}
{\bf 1}_{\{j_1=j_3\}}
{\bf 1}_{\{i_2=i_5\ne 0\}}
{\bf 1}_{\{j_2=j_5\}}-
$$
$$
-
{\bf 1}_{\{i_6=i_4\ne 0\}}
{\bf 1}_{\{j_6=j_4\}}
{\bf 1}_{\{i_1=i_2\ne 0\}}
{\bf 1}_{\{j_1=j_2\}}
{\bf 1}_{\{i_3=i_5\ne 0\}}
{\bf 1}_{\{j_3=j_5\}}-
$$
$$
-
{\bf 1}_{\{i_6=i_5\ne 0\}}
{\bf 1}_{\{j_6=j_5\}}
{\bf 1}_{\{i_1=i_4\ne 0\}}
{\bf 1}_{\{j_1=j_4\}}
{\bf 1}_{\{i_2=i_3\ne 0\}}
{\bf 1}_{\{j_2=j_3\}}-
$$
$$
-
{\bf 1}_{\{i_6=i_5\ne 0\}}
{\bf 1}_{\{j_6=j_5\}}
{\bf 1}_{\{i_1=i_2\ne 0\}}
{\bf 1}_{\{j_1=j_2\}}
{\bf 1}_{\{i_3=i_4\ne 0\}}
{\bf 1}_{\{j_3=j_4\}}-
$$
\begin{equation}
\label{a6}
\Biggl.-
{\bf 1}_{\{i_6=i_5\ne 0\}}
{\bf 1}_{\{j_6=j_5\}}
{\bf 1}_{\{i_1=i_3\ne 0\}}
{\bf 1}_{\{j_1=j_3\}}
{\bf 1}_{\{i_2=i_4\ne 0\}}
{\bf 1}_{\{j_2=j_4\}}\Biggr),
\end{equation}

\vspace{5mm}
\noindent
where ${\bf 1}_A$ is the indicator of the set $A$.

Note that we will consider the case $i_1,\ldots,i_6=1,\ldots,m$.
This case corresponds to the numerical method (\ref{4.45}).

For further consideration, let us 
consider the generalization of formulas (\ref{a1})--(\ref{a6})                 
for the case of an arbitrary multiplicity $k$ $(k\in\mathbb{N})$ of 
the iterated Ito stochastic integral $J[\psi^{(k)}]_{T,t}$ defined by (\ref{ito}).
In order to do this, let us
introduce some notations. 
Consider the unordered
set $\{1, 2, \ldots, k\}$ 
and separate it into two parts:
the first part consists of $r$ unordered 
pairs (sequence order of these pairs is also unimportant) and the 
second one consists of the 
remaining $k-2r$ numbers.
So, we have

\begin{equation}
\label{leto5007}
(\{
\underbrace{\{g_1, g_2\}, \ldots, 
\{g_{2r-1}, g_{2r}\}}_{\small{\hbox{part 1}}}
\},
\{\underbrace{q_1, \ldots, q_{k-2r}}_{\small{\hbox{part 2}}}
\}),
\end{equation}

\vspace{4mm}
\noindent
where 
$\{g_1, g_2, \ldots, 
g_{2r-1}, g_{2r}, q_1, \ldots, q_{k-2r}\}=\{1, 2, \ldots, k\},$
braces   
mean an unordered 
set, and pa\-ren\-the\-ses mean an ordered set.

We will say that (\ref{leto5007}) is a partition 
and consider the sum with respect to all possible
partitions

\begin{equation}
\label{leto5008}
\sum_{\stackrel{(\{\{g_1, g_2\}, \ldots, 
\{g_{2r-1}, g_{2r}\}\}, \{q_1, \ldots, q_{k-2r}\})}
{{}_{\{g_1, g_2, \ldots, 
g_{2r-1}, g_{2r}, q_1, \ldots, q_{k-2r}\}=\{1, 2, \ldots, k\}}}}
a_{g_1 g_2, \ldots, 
g_{2r-1} g_{2r}, q_1 \ldots q_{k-2r}}.
\end{equation}

\vspace{4mm}

Below there are several examples of sums in the form (\ref{leto5008})

\vspace{2mm}
$$
\sum_{\stackrel{(\{g_1, g_2\})}{{}_{\{g_1, g_2\}=\{1, 2\}}}}
a_{g_1 g_2}=a_{12},
$$

\vspace{3mm}
$$
\sum_{\stackrel{(\{\{g_1, g_2\}, \{g_3, g_4\}\})}
{{}_{\{g_1, g_2, g_3, g_4\}=\{1, 2, 3, 4\}}}}
a_{g_1 g_2 g_3 g_4}=a_{1234} + a_{1324} + a_{2314},
$$

\vspace{3mm}
$$
\sum_{\stackrel{(\{g_1, g_2\}, \{q_1, q_{2}\})}
{{}_{\{g_1, g_2, q_1, q_{2}\}=\{1, 2, 3, 4\}}}}
a_{g_1 g_2, q_1 q_{2}}=
$$

$$
=a_{12,34}+a_{13,24}+a_{14,23}
+a_{23,14}+a_{24,13}+a_{34,12},
$$

\vspace{3mm}
$$
\sum_{\stackrel{(\{g_1, g_2\}, \{q_1, q_{2}, q_3\})}
{{}_{\{g_1, g_2, q_1, q_{2}, q_3\}=\{1, 2, 3, 4, 5\}}}}
a_{g_1 g_2, q_1 q_{2}q_3}
=
$$

$$
=a_{12,345}+a_{13,245}+a_{14,235}
+a_{15,234}+a_{23,145}+a_{24,135}+
$$
$$
+a_{25,134}+a_{34,125}+a_{35,124}+a_{45,123},
$$

\vspace{4mm}
$$
\sum_{\stackrel{(\{\{g_1, g_2\}, \{g_3, g_{4}\}\}, \{q_1\})}
{{}_{\{g_1, g_2, g_3, g_{4}, q_1\}=\{1, 2, 3, 4, 5\}}}}
a_{g_1 g_2, g_3 g_{4},q_1}
=
$$

$$
=
a_{12,34,5}+a_{13,24,5}+a_{14,23,5}+
a_{12,35,4}+a_{13,25,4}+a_{15,23,4}+
a_{12,54,3}+a_{15,24,3}+a_{14,25,3}+
$$
$$
+a_{15,34,2}+a_{13,54,2}+a_{14,53,2}+
a_{52,34,1}+a_{53,24,1}+a_{54,23,1}.
$$

\vspace{6mm}

Now we can write (\ref{tyyy}) as

\vspace{1mm}

$$
J[\psi^{(k)}]_{T,t}=
\hbox{\vtop{\offinterlineskip\halign{
\hfil#\hfil\cr
{\rm l.i.m.}\cr
$\stackrel{}{{}_{p_1,\ldots,p_k\to \infty}}$\cr
}} }
\sum\limits_{j_1=0}^{p_1}\ldots
\sum\limits_{j_k=0}^{p_k}
C_{j_k\ldots j_1}\Biggl(
\prod_{l=1}^k\zeta_{j_l}^{(i_l)}+\sum\limits_{r=1}^{[k/2]}
(-1)^r \times
\Biggr.
$$

\vspace{3mm}
\begin{equation}
\label{leto6000hh}
\times
\sum_{\stackrel{(\{\{g_1, g_2\}, \ldots, 
\{g_{2r-1}, g_{2r}\}\}, \{q_1, \ldots, q_{k-2r}\})}
{{}_{\{g_1, g_2, \ldots, 
g_{2r-1}, g_{2r}, q_1, \ldots, q_{k-2r}\}=\{1, 2, \ldots, k\}}}}
\prod\limits_{s=1}^r
{\bf 1}_{\{i_{g_{{}_{2s-1}}}=~i_{g_{{}_{2s}}}\ne 0\}}
\Biggl.{\bf 1}_{\{j_{g_{{}_{2s-1}}}=~j_{g_{{}_{2s}}}\}}
\prod_{l=1}^{k-2r}\zeta_{j_{q_l}}^{(i_{q_l})}\Biggr),
\end{equation}

\vspace{5mm}
\noindent
where $[x]$ is an integer part of a real number $x;$
another notations are the same as in Theorem {\bf 1}.

\vspace{2mm}

In particular, from (\ref{leto6000hh}) for $k=5$ we obtain

\vspace{3mm}

$$
J[\psi^{(5)}]_{T,t}=
\hbox{\vtop{\offinterlineskip\halign{
\hfil#\hfil\cr
{\rm l.i.m.}\cr
$\stackrel{}{{}_{p_1,\ldots,p_5\to \infty}}$\cr
}} }\sum_{j_1=0}^{p_1}\ldots\sum_{j_5=0}^{p_5}
C_{j_5\ldots j_1}\Biggl(
\prod_{l=1}^5\zeta_{j_l}^{(i_l)}-\Biggr.
$$

\vspace{2mm}
$$
-
\sum\limits_{\stackrel{(\{g_1, g_2\}, \{q_1, q_{2}, q_3\})}
{{}_{\{g_1, g_2, q_{1}, q_{2}, q_3\}=\{1, 2, 3, 4, 5\}}}}
{\bf 1}_{\{i_{g_{{}_{1}}}=~i_{g_{{}_{2}}}\ne 0\}}
{\bf 1}_{\{j_{g_{{}_{1}}}=~j_{g_{{}_{2}}}\}}
\prod_{l=1}^{3}\zeta_{j_{q_l}}^{(i_{q_l})}+
$$

\vspace{2mm}
$$
+
\sum_{\stackrel{(\{\{g_1, g_2\}, 
\{g_{3}, g_{4}\}\}, \{q_1\})}
{{}_{\{g_1, g_2, g_{3}, g_{4}, q_1\}=\{1, 2, 3, 4, 5\}}}}
{\bf 1}_{\{i_{g_{{}_{1}}}=~i_{g_{{}_{2}}}\ne 0\}}
{\bf 1}_{\{j_{g_{{}_{1}}}=~j_{g_{{}_{2}}}\}}
\Biggl.{\bf 1}_{\{i_{g_{{}_{3}}}=~i_{g_{{}_{4}}}\ne 0\}}
{\bf 1}_{\{j_{g_{{}_{3}}}=~j_{g_{{}_{4}}}\}}
\zeta_{j_{q_1}}^{(i_{q_1})}\Biggr).
$$

\vspace{7mm}
\noindent
The last equality obviously agrees with
(\ref{a5}).

Let us consider the generalization of Theorem 1 for the case
of an arbitrary complete orthonormal systems  
of functions in the space $L_2([t,T])$ 
and $\psi_1(\tau),\ldots,\psi_k(\tau)\in L_2([t, T]).$

\vspace{2mm}

{\bf Theorem~2}\ \cite{2018a} (Sect.~1.11), \cite{arxiv-1} (Sect.~15).
{\it Suppose that
$\psi_1(\tau),\ldots,\psi_k(\tau)\in L_2([t, T])$ and
$\{\phi_j(x)\}_{j=0}^{\infty}$ is an arbitrary complete orthonormal system  
of functions in the space $L_2([t,T]).$
Then the following expansion

\vspace{1mm}
$$
J[\psi^{(k)}]_{T,t}=
\hbox{\vtop{\offinterlineskip\halign{
\hfil#\hfil\cr
{\rm l.i.m.}\cr
$\stackrel{}{{}_{p_1,\ldots,p_k\to \infty}}$\cr
}} }
\sum\limits_{j_1=0}^{p_1}\ldots
\sum\limits_{j_k=0}^{p_k}
C_{j_k\ldots j_1}\Biggl(
\prod_{l=1}^k\zeta_{j_l}^{(i_l)}+\sum\limits_{r=1}^{[k/2]}
(-1)^r \times
\Biggr.
$$

\vspace{2mm}
\begin{equation}
\label{leto6000}
\times
\sum_{\stackrel{(\{\{g_1, g_2\}, \ldots, 
\{g_{2r-1}, g_{2r}\}\}, \{q_1, \ldots, q_{k-2r}\})}
{{}_{\{g_1, g_2, \ldots, 
g_{2r-1}, g_{2r}, q_1, \ldots, q_{k-2r}\}=\{1, 2, \ldots, k\}}}}
\prod\limits_{s=1}^r
{\bf 1}_{\{i_{g_{{}_{2s-1}}}=~i_{g_{{}_{2s}}}\ne 0\}}
\Biggl.{\bf 1}_{\{j_{g_{{}_{2s-1}}}=~j_{g_{{}_{2s}}}\}}
\prod_{l=1}^{k-2r}\zeta_{j_{q_l}}^{(i_{q_l})}\Biggr)
\end{equation}

\vspace{6mm}
\noindent
con\-verg\-ing in the mean-square sense is valid,
where $[x]$ is an integer part of a real number $x;$
another notations are the same as in Theorem~{\rm 1}.}

\vspace{2mm}

It should be noted that an analogue of Theorem 2 was considered 
in \cite{Rybakov1000}. 
Note that we use another notations 
\cite{2018a} (Sect.~1.11), \cite{arxiv-1} (Sect.~15)
in comparison with \cite{Rybakov1000}.
Moreover, the proof of an analogue of Theorem 2
from \cite{Rybakov1000} is somewhat different from the proof given in 
\cite{2018a} (Sect.~1.11), \cite{arxiv-1} (Sect.~15).

Let us consider the exact calculation and 
effective estimation
of the mean-square error of appro\-xi\-ma\-ti\-on
$J[\psi^{(k)}]_{T,t}^{q}$. Here $J[\psi^{(k)}]_{T,t}^{q}$ is the expression 
on the right-hand side of (\ref{leto6000}) before passing  to the limit 
$\hbox{\vtop{\offinterlineskip\halign{
\hfil#\hfil\cr
{\rm l.i.m.}\cr
$\stackrel{}{{}_{p_1,\ldots,p_k\to \infty}}$\cr
}} }$
for the case
$p_1=\ldots=p_k=q$

$$
J[\psi^{(k)}]_{T,t}^q=
\sum\limits_{j_1,\ldots,j_k=0}^{q}
C_{j_k\ldots j_1}\Biggl(
\prod_{l=1}^k\zeta_{j_l}^{(i_l)}+\sum\limits_{r=1}^{[k/2]}
(-1)^r \times
\Biggr.
$$

\vspace{2mm}
$$
\times
\sum_{\stackrel{(\{\{g_1, g_2\}, \ldots, 
\{g_{2r-1}, g_{2r}\}\}, \{q_1, \ldots, q_{k-2r}\})}
{{}_{\{g_1, g_2, \ldots, 
g_{2r-1}, g_{2r}, q_1, \ldots, q_{k-2r}\}=\{1, 2, \ldots, k\}}}}
\prod\limits_{s=1}^r
{\bf 1}_{\{i_{g_{{}_{2s-1}}}=~i_{g_{{}_{2s}}}\ne 0\}}
\Biggl.{\bf 1}_{\{j_{g_{{}_{2s-1}}}=~j_{g_{{}_{2s}}}\}}
\prod_{l=1}^{k-2r}\zeta_{j_{q_l}}^{(i_{q_l})}\Biggr),
$$

\vspace{5mm}
\noindent
where $[x]$ is an integer part of a real number $x;$
another notations are the same as in Theorems~{\rm 1, 2}.

\vspace{5mm}

Let us denote

$$
{\sf M}\left\{\left(J[\psi^{(k)}]_{T,t}-
J[\psi^{(k)}]_{T,t}^{q}\right)^2\right\}\stackrel{{\rm def}}
{=}E_k^{q},\ \ \ 
\int\limits_{[t,T]^k}
K^2(t_1,\ldots,t_k)dt_1\ldots dt_k
\stackrel{{\rm def}}{=}I_k.
$$

\vspace{4mm}

In \cite{2017-1}, \cite{2017-1a}, \cite{5-002}, 
\cite{arxiv-2}, \cite{2018a}-\cite{2018aaa} it was shown that 

\vspace{-1mm}
\begin{equation}
\label{qq4}
E_k^{q}\le k!\left(I_k-\sum_{j_1,\ldots,j_k=0}^{q}C^2_{j_k\ldots j_1}\right),
\end{equation}

\vspace{2mm}
\noindent
where $i_1,\ldots,i_k=1,\ldots,m$ $(T-t\in (0, +\infty))$
or $i_1,\ldots,i_k=0, 1,\ldots,m$ $(T-t\in (0, 1))$.

\vspace{2mm}

The value $E_k^{q}$
can be calculated exactly.

\vspace{2mm}

{\bf Theorem 3} \cite{2018a} (Sect.~1.12), \cite{arxiv-2} (Sect.~6).
{\it Suppose that $\{\phi_j(x)\}_{j=0}^{\infty}$ 
is an arbitrary complete orthonormal system  
of functions in the space $L_2([t,T])$ and
$\psi_1(\tau),\ldots,\psi_k(\tau)\in L_2([t, T]),$  $i_1,\ldots, i_k=1,\ldots,m$.
Then

\begin{equation}
\label{tttr11}
E_k^q=I_k- \sum_{j_1,\ldots, j_k=0}^{q}
C_{j_k\ldots j_1}
{\sf M}\left\{J[\psi^{(k)}]_{T,t}
\sum\limits_{(j_1,\ldots,j_k)}
\int\limits_t^T \phi_{j_k}(t_k)
\ldots
\int\limits_t^{t_{2}}\phi_{j_{1}}(t_{1})
d{\bf f}_{t_1}^{(i_1)}\ldots
d{\bf f}_{t_k}^{(i_k)}\right\},
\end{equation}

\vspace{5mm}
\noindent
where $i_1,\ldots,i_k = 1,\ldots,m;$
the expression 

\vspace{-1mm}
$$
\sum\limits_{(j_1,\ldots,j_k)}
$$ 

\vspace{3mm}
\noindent
means the sum with respect to all
possible permutations 
$(j_1,\ldots,j_k)$. At the same time if 
$j_r$ swapped with $j_q$ in the permutation $(j_1,\ldots,j_k),$
then $i_r$ swapped with $i_q$ in the permutation
$(i_1,\ldots,i_k);$
another notations are the same as in Theorems {\rm 1, 2.}
}

\vspace{2mm}

Note that 

\vspace{-2mm}
$$
{\sf M}\left\{J[\psi^{(k)}]_{T,t}
\int\limits_t^T \phi_{j_k}(t_k)
\ldots
\int\limits_t^{t_{2}}\phi_{j_{1}}(t_{1})
d{\bf f}_{t_1}^{(i_1)}\ldots
d{\bf f}_{t_k}^{(i_k)}\right\}=C_{j_k\ldots j_1}.
$$

\vspace{4mm}

Then from Theorem 3 for pairwise different $i_1,\ldots,i_k$ 
and for $i_1=\ldots=i_k$
we obtain \cite{2017-1a}, \cite{5-002}, 
\cite{arxiv-2}, \cite{2018a}-\cite{2018aaa}

\vspace{-4mm}
\begin{equation}
\label{qq1}
E_k^q= I_k- \sum_{j_1,\ldots,j_k=0}^{q}
C_{j_k\ldots j_1}^2,
\end{equation}

$$
E_k^q= I_k - \sum_{j_1,\ldots,j_k=0}^{q}
C_{j_k\ldots j_1}\Biggl(\sum\limits_{(j_1,\ldots,j_k)}
C_{j_k\ldots j_1}\Biggr),
$$

\vspace{3mm}
\noindent
where 
$$
\sum\limits_{(j_1,\ldots,j_k)}
$$ 

\vspace{2mm}
\noindent
is a sum with respect to all 
possible permutations
$(j_1,\ldots,j_k)$.

Consider some examples \cite{2017-1a}, \cite{5-002}, 
\cite{arxiv-2}, \cite{2018a}-\cite{2018aaa} of application of Theorem 3
$(i_1,\ldots,i_5=1,\ldots,m)$

\vspace{1mm}
\begin{equation}
\label{qq2}
E_2^q     
=I_2
-\sum_{j_1,j_2=0}^q
C_{j_2j_1}^2-
\sum_{j_1,j_2=0}^q
C_{j_2j_1}C_{j_1j_2}\ \ \ (i_1=i_2),
\end{equation}

\vspace{2mm}
\begin{equation}
\label{qq3}
E_3^q=I_3
-\sum_{j_3,j_2,j_1=0}^q C_{j_3j_2j_1}^2-
\sum_{j_3,j_2,j_1=0}^q C_{j_3j_1j_2}C_{j_3j_2j_1}\ \ \ (i_1=i_2\ne i_3),
\end{equation}

\vspace{2mm}
\begin{equation}
\label{882}
E_3^q=I_3-
\sum_{j_3,j_2,j_1=0}^q C_{j_3j_2j_1}^2-
\sum_{j_3,j_2,j_1=0}^q C_{j_2j_3j_1}C_{j_3j_2j_1}\ \ \ (i_1\ne i_2=i_3),
\end{equation}

\vspace{2mm}
\begin{equation}
\label{883}
E_3^q=I_3
-\sum_{j_3,j_2,j_1=0}^q C_{j_3j_2j_1}^2-
\sum_{j_3,j_2,j_1=0}^q C_{j_3j_2j_1}C_{j_1j_2j_3}\ \ \ (i_1=i_3\ne i_2),
\end{equation}

\vspace{3mm}
$$
E^q_4 = I_4 - \sum_{j_1,\ldots,j_4=0}^{q}
C_{j_4\ldots j_1}\Biggl(\sum\limits_{(j_1,j_2)}
C_{j_4\ldots j_1}\Biggr)\ \ \ (i_1=i_2\ne i_3, i_4;\ i_3\ne i_4),
$$

\vspace{3mm}
$$
E^q_4 = I_4 - \sum_{j_1,\ldots,j_4=0}^{q}
C_{j_4\ldots j_1}\Biggl(\sum\limits_{(j_1,j_3)}
C_{j_4\ldots j_1}\Biggr)\ \ \ (i_1=i_3\ne i_2, i_4;\ i_2\ne i_4),
$$

\vspace{3mm}
$$
E_4^q = I_4 -
\sum_{j_1,\ldots,j_4=0}^{q}
C_{j_4\ldots j_1}\Biggl(\sum\limits_{(j_1,j_2,j_3)}
C_{j_4\ldots j_1}\Biggr)\ \ \ (i_1=i_2=i_3\ne i_4),
$$

\vspace{3mm}
$$
E_4^q = I_4 -
 \sum_{j_1,\ldots,j_4=0}^{q}
C_{j_4\ldots j_1}\Biggl(\sum\limits_{(j_2,j_3,j_4)}
C_{j_4\ldots j_1}\Biggr)\ \ \ (i_2=i_3=i_4\ne i_1),
$$

\vspace{3mm}
$$
E^q_4 = I_4 - \sum_{j_1,\ldots,j_4=0}^{q}
C_{j_4\ldots j_1}\Biggl(\sum\limits_{(j_1,j_2)}\Biggl(
\sum\limits_{(j_3,j_4)}
C_{j_4\ldots j_1}\Biggr)\Biggr)\ \ \ (i_1=i_2\ne i_3=i_4),
$$

\vspace{3mm}
$$
E^q_4 = I_4 - \sum_{j_1,\ldots,j_4=0}^{q}
C_{j_4\ldots j_1}\Biggl(\sum\limits_{(j_1,j_3)}\Biggl(
\sum\limits_{(j_2,j_4)}
C_{j_4\ldots j_1}\Biggr)\Biggr)\ \ \ (i_1=i_3\ne i_2=i_4),
$$

\vspace{4mm}
$$
E_5^q = I_5 - \sum_{j_1,\ldots,j_5=0}^{q}
C_{j_5\ldots j_1}\Biggl(\sum\limits_{(j_2,j_4)}\Biggl(
\sum\limits_{(j_3,j_5)}
C_{j_5\ldots j_1}\Biggr)\Biggr)\ \ \ (i_1\ne i_2=i_4\ne i_3=i_5\ne i_1),
$$

\vspace{4mm}
$$
E^q_5 = I_5 - \sum_{j_1,\ldots,j_5=0}^{q}
C_{j_5\ldots j_1}\Biggl(\sum\limits_{(j_4,j_5)}\Biggl(
\sum\limits_{(j_1,j_2,j_3)}
C_{j_5\ldots j_1}\Biggr)\Biggr)\ \ \ (i_1=i_2=i_3\ne i_4=i_5),
$$

\vspace{4mm}
$$
E^q_5 = I_5 - \sum_{j_1,\ldots,j_5=0}^{q}
C_{j_5\ldots j_1}\Biggl(\sum\limits_{(j_1,j_3,j_4,j_5)}
C_{j_5\ldots j_1}\Biggr)\ \ \ (i_1=i_3=i_4=i_5\ne i_2).
$$

\vspace{6mm}

The values $E_4^q$ and $E_5^q$ were calculated exaclty for all possible 
combinations of $i_1,\ldots,i_5=1,\ldots,m$ in 
\cite{2017-1}, \cite{2017-1a}, \cite{arxiv-2},
\cite{2018a}-\cite{2018aaa}.

Let us consider the approximations of iterated Ito stochastic integrals
from (\ref{4.45}) using (\ref{a1})--(\ref{a6}) 
and complete orthonormal system of Legendre 
polynomials in the space $L_2([\tau_p,\tau_{p+1}])$\
($\tau_p=p\Delta,$ $N\Delta=\bar T,$ 
$p=0,1,\ldots,N$) \cite{2017-1a}
(also see \cite{kuz1997}-\cite{2017-1}, \cite{2007-1}-\cite{2018aaa})

\vspace{2mm}
$$
I_{{0}_{\tau_{p+1},\tau_p}}^{(i_1)}=\sqrt{\Delta}\zeta_0^{(i_1)},
$$

\vspace{2mm}
\begin{equation}
\label{qqqq1}
I_{{00}_{\tau_{p+1},\tau_p}}^{(i_1 i_2)q}=
\frac{\Delta}{2}\left(\zeta_0^{(i_1)}\zeta_0^{(i_2)}+\sum_{i=1}^{q}
\frac{1}{\sqrt{4i^2-1}}\left(
\zeta_{i-1}^{(i_1)}\zeta_{i}^{(i_2)}-
\zeta_i^{(i_1)}\zeta_{i-1}^{(i_2)}\right) - {\bf 1}_{\{i_1=i_2\}}\right),
\end{equation}

\vspace{5mm}

$$
I_{{1}_{\tau_{p+1},\tau_p}}^{(i_1)}=
-\frac{{\Delta}^{3/2}}{2}\left(\zeta_0^{(i_1)}+
\frac{1}{\sqrt{3}}\zeta_1^{(i_1)}\right),
$$

\vspace{5mm}

$$
I_{{000}_{\tau_{p+1},\tau_p}}^{(i_1i_2i_3)q}
=\sum_{j_1,j_2,j_3=0}^{q}
C_{j_3j_2j_1}\Biggl(
\zeta_{j_1}^{(i_1)}\zeta_{j_2}^{(i_2)}\zeta_{j_3}^{(i_3)}
-{\bf 1}_{\{i_1=i_2\}}
{\bf 1}_{\{j_1=j_2\}}
\zeta_{j_3}^{(i_3)}-
\Biggr.
$$
$$
\Biggl.
-{\bf 1}_{\{i_2=i_3\}}
{\bf 1}_{\{j_2=j_3\}}
\zeta_{j_1}^{(i_1)}-
{\bf 1}_{\{i_1=i_3\}}
{\bf 1}_{\{j_1=j_3\}}
\zeta_{j_2}^{(i_2)}\Biggr),
$$

\vspace{5mm}

$$
I_{{0000}_{\tau_{p+1},\tau_p}}^{(i_1 i_2 i_3 i_4)q}
=\sum_{j_1,j_2,j_3,j_4=0}^{q}
C_{j_4 j_3 j_2 j_1}\Biggl(
\prod_{l=1}^4\zeta_{j_l}^{(i_l)}
-\Biggr.
$$
$$
-
{\bf 1}_{\{i_1=i_2\}}
{\bf 1}_{\{j_1=j_2\}}
\zeta_{j_3}^{(i_3)}
\zeta_{j_4}^{(i_4)}
-
{\bf 1}_{\{i_1=i_3\}}
{\bf 1}_{\{j_1=j_3\}}
\zeta_{j_2}^{(i_2)}
\zeta_{j_4}^{(i_4)}-
$$
$$
-
{\bf 1}_{\{i_1=i_4\}}
{\bf 1}_{\{j_1=j_4\}}
\zeta_{j_2}^{(i_2)}
\zeta_{j_3}^{(i_3)}
-
{\bf 1}_{\{i_2=i_3\}}
{\bf 1}_{\{j_2=j_3\}}
\zeta_{j_1}^{(i_1)}
\zeta_{j_4}^{(i_4)}-
$$
$$
-
{\bf 1}_{\{i_2=i_4\}}
{\bf 1}_{\{j_2=j_4\}}
\zeta_{j_1}^{(i_1)}
\zeta_{j_3}^{(i_3)}
-
{\bf 1}_{\{i_3=i_4\}}
{\bf 1}_{\{j_3=j_4\}}
\zeta_{j_1}^{(i_1)}
\zeta_{j_2}^{(i_2)}+
$$
$$
+
{\bf 1}_{\{i_1=i_2\}}
{\bf 1}_{\{j_1=j_2\}}
{\bf 1}_{\{i_3=i_4\}}
{\bf 1}_{\{j_3=j_4\}}+
{\bf 1}_{\{i_1=i_3\}}
{\bf 1}_{\{j_1=j_3\}}
{\bf 1}_{\{i_2=i_4\}}
{\bf 1}_{\{j_2=j_4\}}+
$$
$$
+\Biggl.
{\bf 1}_{\{i_1=i_4\}}
{\bf 1}_{\{j_1=j_4\}}
{\bf 1}_{\{i_2=i_3\}}
{\bf 1}_{\{j_2=j_3\}}\Biggr),
$$

\vspace{6mm}

$$
I_{{01}_{\tau_{p+1},\tau_p}}^{(i_1 i_2)q}=
-\frac{\Delta}{2}
I_{{00}_{\tau_{p+1},\tau_p}}^{(i_1 i_2)q}
-\frac{{\Delta}^2}{4}\Biggl(
\frac{1}{\sqrt{3}}\zeta_0^{(i_1)}\zeta_1^{(i_2)}+\Biggr.
$$

\vspace{1mm}
$$
+\Biggl.\sum_{i=0}^{q}\Biggl(
\frac{(i+2)\zeta_i^{(i_1)}\zeta_{i+2}^{(i_2)}
-(i+1)\zeta_{i+2}^{(i_1)}\zeta_{i}^{(i_2)}}
{\sqrt{(2i+1)(2i+5)}(2i+3)}-
\frac{\zeta_i^{(i_1)}\zeta_{i}^{(i_2)}}{(2i-1)(2i+3)}\Biggr)\Biggr),
$$

\vspace{6mm}

$$
I_{{10}_{\tau_{p+1},\tau_p}}^{(i_1 i_2)q}=
-\frac{\Delta}{2}I_{{00}_{\tau_{p+1},\tau_p}}^{(i_1 i_2)q}
-\frac{\Delta^2}{4}\Biggl(
\frac{1}{\sqrt{3}}\zeta_0^{(i_2)}\zeta_1^{(i_1)}+\Biggr.
$$

\vspace{1mm}
$$
+\Biggl.\sum_{i=0}^{q}\Biggl(
\frac{(i+1)\zeta_{i+2}^{(i_2)}\zeta_{i}^{(i_1)}
-(i+2)\zeta_{i}^{(i_2)}\zeta_{i+2}^{(i_1)}}
{\sqrt{(2i+1)(2i+5)}(2i+3)}+
\frac{\zeta_i^{(i_1)}\zeta_{i}^{(i_2)}}{(2i-1)(2i+3)}\Biggr)\Biggr),
$$

\vspace{8mm}

$$
{I}_{2_{\tau_{p+1},\tau_p}}^{(i_1)}=
\frac{\Delta^{5/2}}{3}\left(
\zeta_0^{(i_1)}+\frac{\sqrt{3}}{2}\zeta_1^{(i_1)}+
\frac{1}{2\sqrt{5}}\zeta_2^{(i_1)}\right),
$$

\vspace{6mm}

$$
I_{{001}_{\tau_{p+1},\tau_p}}^{(i_1i_2i_3)q}
=\sum_{j_1,j_2,j_3=0}^{q}
C_{j_3j_2j_1}^{001}\Biggl(
\zeta_{j_1}^{(i_1)}\zeta_{j_2}^{(i_2)}\zeta_{j_3}^{(i_3)}
-{\bf 1}_{\{i_1=i_2\}}
{\bf 1}_{\{j_1=j_2\}}
\zeta_{j_3}^{(i_3)}-
\Biggr.
$$
$$
\Biggl.
-{\bf 1}_{\{i_2=i_3\}}
{\bf 1}_{\{j_2=j_3\}}
\zeta_{j_1}^{(i_1)}-
{\bf 1}_{\{i_1=i_3\}}
{\bf 1}_{\{j_1=j_3\}}
\zeta_{j_2}^{(i_2)}\Biggr),
$$

\vspace{6mm}

$$
I_{{010}_{\tau_{p+1},\tau_p}}^{(i_1i_2i_3)q}
=\sum_{j_1,j_2,j_3=0}^{q}
C_{j_3j_2j_1}^{010}\Biggl(
\zeta_{j_1}^{(i_1)}\zeta_{j_2}^{(i_2)}\zeta_{j_3}^{(i_3)}
-{\bf 1}_{\{i_1=i_2\}}
{\bf 1}_{\{j_1=j_2\}}
\zeta_{j_3}^{(i_3)}-
\Biggr.
$$
$$
\Biggl.
-{\bf 1}_{\{i_2=i_3\}}
{\bf 1}_{\{j_2=j_3\}}
\zeta_{j_1}^{(i_1)}-
{\bf 1}_{\{i_1=i_3\}}
{\bf 1}_{\{j_1=j_3\}}
\zeta_{j_2}^{(i_2)}\Biggr),
$$

\vspace{6mm}

$$
I_{{100}_{\tau_{p+1},\tau_p}}^{(i_1i_2i_3)q}
=\sum_{j_1,j_2,j_3=0}^{q}
C_{j_3j_2j_1}^{100}\Biggl(
\zeta_{j_1}^{(i_1)}\zeta_{j_2}^{(i_2)}\zeta_{j_3}^{(i_3)}
-{\bf 1}_{\{i_1=i_2\}}
{\bf 1}_{\{j_1=j_2\}}
\zeta_{j_3}^{(i_3)}-
\Biggr.
$$
$$
\Biggl.
-{\bf 1}_{\{i_2=i_3\}}
{\bf 1}_{\{j_2=j_3\}}
\zeta_{j_1}^{(i_1)}-
{\bf 1}_{\{i_1=i_3\}}
{\bf 1}_{\{j_1=j_3\}}
\zeta_{j_2}^{(i_2)}\Biggr),
$$

\vspace{7mm}

$$
I_{{00000}_{\tau_{p+1},\tau_p}}^{(i_1 i_2 i_3 i_4 i_5)q}
=\sum_{j_1,j_2,j_3,j_4,j_5=0}^q
C_{j_5 j_4 j_3 j_2 j_1}\Biggl(
\prod_{l=1}^5\zeta_{j_l}^{(i_l)}
-\Biggr.
$$
$$
-
{\bf 1}_{\{j_1=j_2\}}
{\bf 1}_{\{i_1=i_2\}}
\zeta_{j_3}^{(i_3)}
\zeta_{j_4}^{(i_4)}
\zeta_{j_5}^{(i_5)}-
{\bf 1}_{\{j_1=j_3\}}
{\bf 1}_{\{i_1=i_3\}}
\zeta_{j_2}^{(i_2)}
\zeta_{j_4}^{(i_4)}
\zeta_{j_5}^{(i_5)}-
$$
$$
-
{\bf 1}_{\{j_1=j_4\}}
{\bf 1}_{\{i_1=i_4\}}
\zeta_{j_2}^{(i_2)}
\zeta_{j_3}^{(i_3)}
\zeta_{j_5}^{(i_5)}-
{\bf 1}_{\{j_1=j_5\}}
{\bf 1}_{\{i_1=i_5\}}
\zeta_{j_2}^{(i_2)}
\zeta_{j_3}^{(i_3)}
\zeta_{j_4}^{(i_4)}-
$$
$$
-
{\bf 1}_{\{j_2=j_3\}}
{\bf 1}_{\{i_2=i_3\}}
\zeta_{j_1}^{(i_1)}
\zeta_{j_4}^{(i_4)}
\zeta_{j_5}^{(i_5)}-
{\bf 1}_{\{j_2=j_4\}}
{\bf 1}_{\{i_2=i_4\}}
\zeta_{j_1}^{(i_1)}
\zeta_{j_3}^{(i_3)}
\zeta_{j_5}^{(i_5)}-
$$
$$
-
{\bf 1}_{\{j_2=j_5\}}
{\bf 1}_{\{i_2=i_5\}}
\zeta_{j_1}^{(i_1)}
\zeta_{j_3}^{(i_3)}
\zeta_{j_4}^{(i_4)}-{\bf 1}_{\{j_3=j_4\}}
{\bf 1}_{\{i_3=i_4\}}
\zeta_{j_1}^{(i_1)}
\zeta_{j_2}^{(i_2)}
\zeta_{j_5}^{(i_5)}-
$$
$$
-
{\bf 1}_{\{j_3=j_5\}}
{\bf 1}_{\{i_3=i_5\}}
\zeta_{j_1}^{(i_1)}
\zeta_{j_2}^{(i_2)}
\zeta_{j_4}^{(i_4)}-{\bf 1}_{\{j_4=j_5\}}
{\bf 1}_{\{i_4=i_5\}}
\zeta_{j_1}^{(i_1)}
\zeta_{j_2}^{(i_2)}
\zeta_{j_3}^{(i_3)}+
$$
$$
+
{\bf 1}_{\{j_1=j_2\}}
{\bf 1}_{\{i_1=i_2\}}
{\bf 1}_{\{j_3=j_4\}}
{\bf 1}_{\{i_3=i_4\}}\zeta_{j_5}^{(i_5)}+
{\bf 1}_{\{j_1=j_2\}}
{\bf 1}_{\{i_1=i_2\}}
{\bf 1}_{\{j_3=j_5\}}
{\bf 1}_{\{i_3=i_5\}}\zeta_{j_4}^{(i_4)}+
$$
$$
+
{\bf 1}_{\{j_1=j_2\}}
{\bf 1}_{\{i_1=i_2\}}
{\bf 1}_{\{j_4=j_5\}}
{\bf 1}_{\{i_4=i_5\}}\zeta_{j_3}^{(i_3)}+
{\bf 1}_{\{j_1=j_3\}}
{\bf 1}_{\{i_1=i_3\}}
{\bf 1}_{\{j_2=j_4\}}
{\bf 1}_{\{i_2=i_4\}}\zeta_{j_5}^{(i_5)}+
$$
$$
+
{\bf 1}_{\{j_1=j_3\}}
{\bf 1}_{\{i_1=i_3\}}
{\bf 1}_{\{j_2=j_5\}}
{\bf 1}_{\{i_2=i_5\}}\zeta_{j_4}^{(i_4)}+
{\bf 1}_{\{j_1=j_3\}}
{\bf 1}_{\{i_1=i_3\}}
{\bf 1}_{\{j_4=j_5\}}
{\bf 1}_{\{i_4=i_5\}}\zeta_{j_2}^{(i_2)}+
$$
$$
+
{\bf 1}_{\{j_1=j_4\}}
{\bf 1}_{\{i_1=i_4\}}
{\bf 1}_{\{j_2=j_3\}}
{\bf 1}_{\{i_2=i_3\}}\zeta_{j_5}^{(i_5)}+
{\bf 1}_{\{j_1=j_4\}}
{\bf 1}_{\{i_1=i_4\}}
{\bf 1}_{\{j_2=j_5\}}
{\bf 1}_{\{i_2=i_5\}}\zeta_{j_3}^{(i_3)}+
$$
$$
+
{\bf 1}_{\{j_1=j_4\}}
{\bf 1}_{\{i_1=i_4\}}
{\bf 1}_{\{j_3=j_5\}}
{\bf 1}_{\{i_3=i_5\}}\zeta_{j_2}^{(i_2)}+
{\bf 1}_{\{j_1=j_5\}}
{\bf 1}_{\{i_1=i_5\}}
{\bf 1}_{\{j_2=j_3\}}
{\bf 1}_{\{i_2=i_3\}}\zeta_{j_4}^{(i_4)}+
$$
$$
+
{\bf 1}_{\{j_1=j_5\}}
{\bf 1}_{\{i_1=i_5\}}
{\bf 1}_{\{j_2=j_4\}}
{\bf 1}_{\{i_2=i_4\}}\zeta_{j_3}^{(i_3)}+
{\bf 1}_{\{j_1=j_5\}}
{\bf 1}_{\{i_1=i_5\}}
{\bf 1}_{\{j_3=j_4\}}
{\bf 1}_{\{i_3=i_4\}}\zeta_{j_2}^{(i_2)}+
$$
$$
+
{\bf 1}_{\{j_2=j_3\}}
{\bf 1}_{\{i_2=i_3\}}
{\bf 1}_{\{j_4=j_5\}}
{\bf 1}_{\{i_4=i_5\}}\zeta_{j_1}^{(i_1)}+
{\bf 1}_{\{j_2=j_4\}}
{\bf 1}_{\{i_2=i_4\}}
{\bf 1}_{\{j_3=j_5\}}
{\bf 1}_{\{i_3=i_5\}}\zeta_{j_1}^{(i_1)}+
$$
$$
+\Biggl.
{\bf 1}_{\{j_2=j_5\ne 0\}}
{\bf 1}_{\{i_2=i_5\}}
{\bf 1}_{\{j_3=j_4\ne 0\}}
{\bf 1}_{\{i_3=i_4\}}\zeta_{j_1}^{(i_1)}\Biggr),
$$

\vspace{8mm}

$$
I_{02_{\tau_{p+1},\tau_p}}^{(i_1 i_2)q}
=-\frac{{\Delta}^2}{4}I_{00_{\tau_{p+1},\tau_p}}^{(i_1 i_2)q}
-\Delta I_{01_{\tau_{p+1},\tau_p}}^{(i_1 i_2)q}+
\frac{{\Delta}^3}{8}\Biggl[
\frac{2}{3\sqrt{5}}\zeta_2^{(i_2)}\zeta_0^{(i_1)}+\Biggr.
$$

\vspace{1mm}
$$
+\frac{1}{3}\zeta_0^{(i_1)}\zeta_0^{(i_2)}+
\sum_{i=0}^{q}\Biggl(
\frac{(i+2)(i+3)\zeta_{i+3}^{(i_2)}\zeta_{i}^{(i_1)}
-(i+1)(i+2)\zeta_{i}^{(i_2)}\zeta_{i+3}^{(i_1)}}
{\sqrt{(2i+1)(2i+7)}(2i+3)(2i+5)}+
\Biggr.
$$

\vspace{1mm}
\begin{equation}
\label{777}
\Biggl.\Biggl.+
\frac{(i^2+i-3)\zeta_{i+1}^{(i_2)}\zeta_{i}^{(i_1)}
-(i^2+3i-1)\zeta_{i}^{(i_2)}\zeta_{i+1}^{(i_1)}}
{\sqrt{(2i+1)(2i+3)}(2i-1)(2i+5)}\Biggr)\Biggr] - \frac{1}{24}{\bf 1}_{\{i_1=i_2\}}{\Delta^3},
\end{equation}

\vspace{8mm}

$$
I_{20_{\tau_{p+1},\tau_p}}^{(i_1 i_2)q}=-\frac{{\Delta}^2}{4}
I_{00_{\tau_{p+1},\tau_p}}^{(i_1 i_2)q}
-\Delta I_{10_{\tau_{p+1},\tau_p}}^{(i_1 i_2)q}+
\frac{{\Delta}^3}{8}\Biggl[
\frac{2}{3\sqrt{5}}\zeta_0^{(i_2)}\zeta_2^{(i_1)}+\Biggr.
$$

\vspace{1mm}
$$
+\frac{1}{3}\zeta_0^{(i_1)}\zeta_0^{(i_2)}+
\sum_{i=0}^{q}\Biggl(
\frac{(i+1)(i+2)\zeta_{i+3}^{(i_2)}\zeta_{i}^{(i_1)}
-(i+2)(i+3)\zeta_{i}^{(i_2)}\zeta_{i+3}^{(i_1)}}
{\sqrt{(2i+1)(2i+7)}(2i+3)(2i+5)}+
\Biggr.
$$

\vspace{1mm}
\begin{equation}
\label{778}
\Biggl.\Biggl.+
\frac{(i^2+3i-1)\zeta_{i+1}^{(i_2)}\zeta_{i}^{(i_1)}
-(i^2+i-3)\zeta_{i}^{(i_2)}\zeta_{i+1}^{(i_1)}}
{\sqrt{(2i+1)(2i+3)}(2i-1)(2i+5)}\Biggr)\Biggr] - \frac{1}{24}{\bf 1}_{\{i_1=i_2\}}{\Delta^3},
\end{equation}

\vspace{8mm}

$$
I_{11_{\tau_{p+1},\tau_p}}^{(i_1 i_2)q}
=-\frac{{\Delta}^2}{4}I_{00_{\tau_{p+1},\tau_p}}^{(i_1 i_2)q}
-\frac{\Delta}{2}\left(
I_{10_{\tau_{p+1},\tau_p}}^{(i_1 i_2)q}+
I_{01_{\tau_{p+1},\tau_p}}^{(i_1 i_2)q}\right)+
\frac{{\Delta}^3}{8}\Biggl[
\frac{1}{3}\zeta_1^{(i_1)}\zeta_1^{(i_2)}+\Biggr.
$$

\vspace{1mm}
$$
+
\sum_{i=0}^{q}\Biggl(
\frac{(i+1)(i+3)\left(\zeta_{i+3}^{(i_2)}\zeta_{i}^{(i_1)}
-\zeta_{i}^{(i_2)}\zeta_{i+3}^{(i_1)}\right)}
{\sqrt{(2i+1)(2i+7)}(2i+3)(2i+5)}+
\Biggr.
$$

\vspace{1mm}
\begin{equation}
\label{779}
\Biggl.\Biggl.+
\frac{(i+1)^2\left(\zeta_{i+1}^{(i_2)}\zeta_{i}^{(i_1)}
-\zeta_{i}^{(i_2)}\zeta_{i+1}^{(i_1)}\right)}
{\sqrt{(2i+1)(2i+3)}(2i-1)(2i+5)}\Biggr)\Biggr] - \frac{1}{24}{\bf 1}_{\{i_1=i_2\}}{\Delta^3},
\end{equation}

\vspace{8mm}

$$
I_{{0001}_{\tau_{p+1},\tau_p}}^{(i_1 i_2 i_3 i_4)q}
=\sum_{j_1,j_2,j_3,j_4=0}^{q}
C_{j_4 j_3 j_2 j_1}^{0001}\Biggl(
\prod_{l=1}^4\zeta_{j_l}^{(i_l)}
-\Biggr.
$$
$$
-
{\bf 1}_{\{i_1=i_2\}}
{\bf 1}_{\{j_1=j_2\}}
\zeta_{j_3}^{(i_3)}
\zeta_{j_4}^{(i_4)}
-
{\bf 1}_{\{i_1=i_3\}}
{\bf 1}_{\{j_1=j_3\}}
\zeta_{j_2}^{(i_2)}
\zeta_{j_4}^{(i_4)}-
$$
$$
-
{\bf 1}_{\{i_1=i_4\}}
{\bf 1}_{\{j_1=j_4\}}
\zeta_{j_2}^{(i_2)}
\zeta_{j_3}^{(i_3)}
-
{\bf 1}_{\{i_2=i_3\}}
{\bf 1}_{\{j_2=j_3\}}
\zeta_{j_1}^{(i_1)}
\zeta_{j_4}^{(i_4)}-
$$
$$
-
{\bf 1}_{\{i_2=i_4\}}
{\bf 1}_{\{j_2=j_4\}}
\zeta_{j_1}^{(i_1)}
\zeta_{j_3}^{(i_3)}
-
{\bf 1}_{\{i_3=i_4\}}
{\bf 1}_{\{j_3=j_4\}}
\zeta_{j_1}^{(i_1)}
\zeta_{j_2}^{(i_2)}+
$$
$$
+
{\bf 1}_{\{i_1=i_2\}}
{\bf 1}_{\{j_1=j_2\}}
{\bf 1}_{\{i_3=i_4\}}
{\bf 1}_{\{j_3=j_4\}}+
{\bf 1}_{\{i_1=i_3\}}
{\bf 1}_{\{j_1=j_3\}}
{\bf 1}_{\{i_2=i_4\}}
{\bf 1}_{\{j_2=j_4\}}+
$$
$$
+\Biggl.
{\bf 1}_{\{i_1=i_4\}}
{\bf 1}_{\{j_1=j_4\}}
{\bf 1}_{\{i_2=i_3\}}
{\bf 1}_{\{j_2=j_3\}}\Biggr),
$$

\vspace{7mm}

$$
I_{{0010}_{\tau_{p+1},\tau_p}}^{(i_1 i_2 i_3 i_4)q}
=\sum_{j_1,j_2,j_3,j_4=0}^{q}
C_{j_4 j_3 j_2 j_1}^{0010}\Biggl(
\prod_{l=1}^4\zeta_{j_l}^{(i_l)}
-\Biggr.
$$
$$
-
{\bf 1}_{\{i_1=i_2\}}
{\bf 1}_{\{j_1=j_2\}}
\zeta_{j_3}^{(i_3)}
\zeta_{j_4}^{(i_4)}
-
{\bf 1}_{\{i_1=i_3\}}
{\bf 1}_{\{j_1=j_3\}}
\zeta_{j_2}^{(i_2)}
\zeta_{j_4}^{(i_4)}-
$$
$$
-
{\bf 1}_{\{i_1=i_4\}}
{\bf 1}_{\{j_1=j_4\}}
\zeta_{j_2}^{(i_2)}
\zeta_{j_3}^{(i_3)}
-
{\bf 1}_{\{i_2=i_3\}}
{\bf 1}_{\{j_2=j_3\}}
\zeta_{j_1}^{(i_1)}
\zeta_{j_4}^{(i_4)}-
$$
$$
-
{\bf 1}_{\{i_2=i_4\}}
{\bf 1}_{\{j_2=j_4\}}
\zeta_{j_1}^{(i_1)}
\zeta_{j_3}^{(i_3)}
-
{\bf 1}_{\{i_3=i_4\}}
{\bf 1}_{\{j_3=j_4\}}
\zeta_{j_1}^{(i_1)}
\zeta_{j_2}^{(i_2)}+
$$
$$
+
{\bf 1}_{\{i_1=i_2\}}
{\bf 1}_{\{j_1=j_2\}}
{\bf 1}_{\{i_3=i_4\}}
{\bf 1}_{\{j_3=j_4\}}+
{\bf 1}_{\{i_1=i_3\}}
{\bf 1}_{\{j_1=j_3\}}
{\bf 1}_{\{i_2=i_4\}}
{\bf 1}_{\{j_2=j_4\}}+
$$
$$
+\Biggl.
{\bf 1}_{\{i_1=i_4\}}
{\bf 1}_{\{j_1=j_4\}}
{\bf 1}_{\{i_2=i_3\}}
{\bf 1}_{\{j_2=j_3\}}\Biggr),
$$

\vspace{7mm}

$$
I_{{0100}_{\tau_{p+1},\tau_p}}^{(i_1 i_2 i_3 i_4)q}
=\sum_{j_1,j_2,j_3,j_4=0}^{q}
C_{j_4 j_3 j_2 j_1}^{0100}\Biggl(
\prod_{l=1}^4\zeta_{j_l}^{(i_l)}
-\Biggr.
$$
$$
-
{\bf 1}_{\{i_1=i_2\}}
{\bf 1}_{\{j_1=j_2\}}
\zeta_{j_3}^{(i_3)}
\zeta_{j_4}^{(i_4)}
-
{\bf 1}_{\{i_1=i_3\}}
{\bf 1}_{\{j_1=j_3\}}
\zeta_{j_2}^{(i_2)}
\zeta_{j_4}^{(i_4)}-
$$
$$
-
{\bf 1}_{\{i_1=i_4\}}
{\bf 1}_{\{j_1=j_4\}}
\zeta_{j_2}^{(i_2)}
\zeta_{j_3}^{(i_3)}
-
{\bf 1}_{\{i_2=i_3\}}
{\bf 1}_{\{j_2=j_3\}}
\zeta_{j_1}^{(i_1)}
\zeta_{j_4}^{(i_4)}-
$$
$$
-
{\bf 1}_{\{i_2=i_4\}}
{\bf 1}_{\{j_2=j_4\}}
\zeta_{j_1}^{(i_1)}
\zeta_{j_3}^{(i_3)}
-
{\bf 1}_{\{i_3=i_4\}}
{\bf 1}_{\{j_3=j_4\}}
\zeta_{j_1}^{(i_1)}
\zeta_{j_2}^{(i_2)}+
$$
$$
+
{\bf 1}_{\{i_1=i_2\}}
{\bf 1}_{\{j_1=j_2\}}
{\bf 1}_{\{i_3=i_4\}}
{\bf 1}_{\{j_3=j_4\}}+
{\bf 1}_{\{i_1=i_3\}}
{\bf 1}_{\{j_1=j_3\}}
{\bf 1}_{\{i_2=i_4\}}
{\bf 1}_{\{j_2=j_4\}}+
$$
$$
+\Biggl.
{\bf 1}_{\{i_1=i_4\}}
{\bf 1}_{\{j_1=j_4\}}
{\bf 1}_{\{i_2=i_3\}}
{\bf 1}_{\{j_2=j_3\}}\Biggr),
$$

\vspace{7mm}

$$
I_{{1000}_{\tau_{p+1},\tau_p}}^{(i_1 i_2 i_3 i_4)q}
=\sum_{j_1,j_2,j_3,j_4=0}^{q}
C_{j_4 j_3 j_2 j_1}^{1000}\Biggl(
\prod_{l=1}^4\zeta_{j_l}^{(i_l)}
-\Biggr.
$$
$$
-
{\bf 1}_{\{i_1=i_2\}}
{\bf 1}_{\{j_1=j_2\}}
\zeta_{j_3}^{(i_3)}
\zeta_{j_4}^{(i_4)}
-
{\bf 1}_{\{i_1=i_3\}}
{\bf 1}_{\{j_1=j_3\}}
\zeta_{j_2}^{(i_2)}
\zeta_{j_4}^{(i_4)}-
$$
$$
-
{\bf 1}_{\{i_1=i_4\}}
{\bf 1}_{\{j_1=j_4\}}
\zeta_{j_2}^{(i_2)}
\zeta_{j_3}^{(i_3)}
-
{\bf 1}_{\{i_2=i_3\}}
{\bf 1}_{\{j_2=j_3\}}
\zeta_{j_1}^{(i_1)}
\zeta_{j_4}^{(i_4)}-
$$
$$
-
{\bf 1}_{\{i_2=i_4\}}
{\bf 1}_{\{j_2=j_4\}}
\zeta_{j_1}^{(i_1)}
\zeta_{j_3}^{(i_3)}
-
{\bf 1}_{\{i_3=i_4\}}
{\bf 1}_{\{j_3=j_4\}}
\zeta_{j_1}^{(i_1)}
\zeta_{j_2}^{(i_2)}+
$$
$$
+
{\bf 1}_{\{i_1=i_2\}}
{\bf 1}_{\{j_1=j_2\}}
{\bf 1}_{\{i_3=i_4\}}
{\bf 1}_{\{j_3=j_4\}}+
{\bf 1}_{\{i_1=i_3\}}
{\bf 1}_{\{j_1=j_3\}}
{\bf 1}_{\{i_2=i_4\}}
{\bf 1}_{\{j_2=j_4\}}+
$$
$$
+\Biggl.
{\bf 1}_{\{i_1=i_4\}}
{\bf 1}_{\{j_1=j_4\}}
{\bf 1}_{\{i_2=i_3\}}
{\bf 1}_{\{j_2=j_3\}}\Biggr),
$$

\vspace{9mm}

$$
I_{{000000}_{\tau_{p+1},\tau_p}}^{(i_1 i_2 i_3 i_4 i_5 i_6)q}
=\sum_{j_1,j_2,j_3,j_4,j_5,j_6=0}^{q}
C_{j_6 j_5 j_4 j_3 j_2 j_1}\Biggl(
\prod_{l=1}^6
\zeta_{j_l}^{(i_l)}
-\Biggr.
$$
$$
-
{\bf 1}_{\{j_1=j_6\}}
{\bf 1}_{\{i_1=i_6\}}
\zeta_{j_2}^{(i_2)}
\zeta_{j_3}^{(i_3)}
\zeta_{j_4}^{(i_4)}
\zeta_{j_5}^{(i_5)}-
{\bf 1}_{\{j_2=j_6\}}
{\bf 1}_{\{i_2=i_6\}}
\zeta_{j_1}^{(i_1)}
\zeta_{j_3}^{(i_3)}
\zeta_{j_4}^{(i_4)}
\zeta_{j_5}^{(i_5)}-
$$
$$
-
{\bf 1}_{\{j_3=j_6\}}
{\bf 1}_{\{i_3=i_6\}}
\zeta_{j_1}^{(i_1)}
\zeta_{j_2}^{(i_2)}
\zeta_{j_4}^{(i_4)}
\zeta_{j_5}^{(i_5)}-
{\bf 1}_{\{j_4=j_6\}}
{\bf 1}_{\{i_4=i_6\}}
\zeta_{j_1}^{(i_1)}
\zeta_{j_2}^{(i_2)}
\zeta_{j_3}^{(i_3)}
\zeta_{j_5}^{(i_5)}-
$$
$$
-
{\bf 1}_{\{j_5=j_6\}}
{\bf 1}_{\{i_5=i_6\}}
\zeta_{j_1}^{(i_1)}
\zeta_{j_2}^{(i_2)}
\zeta_{j_3}^{(i_3)}
\zeta_{j_4}^{(i_4)}-
{\bf 1}_{\{j_1=j_2\}}
{\bf 1}_{\{i_1=i_2\}}
\zeta_{j_3}^{(i_3)}
\zeta_{j_4}^{(i_4)}
\zeta_{j_5}^{(i_5)}
\zeta_{j_6}^{(i_6)}-
$$
$$
-
{\bf 1}_{\{j_1=j_3\}}
{\bf 1}_{\{i_1=i_3\}}
\zeta_{j_2}^{(i_2)}
\zeta_{j_4}^{(i_4)}
\zeta_{j_5}^{(i_5)}
\zeta_{j_6}^{(i_6)}-
{\bf 1}_{\{j_1=j_4\}}
{\bf 1}_{\{i_1=i_4\}}
\zeta_{j_2}^{(i_2)}
\zeta_{j_3}^{(i_3)}
\zeta_{j_5}^{(i_5)}
\zeta_{j_6}^{(i_6)}-
$$
$$
-
{\bf 1}_{\{j_1=j_5\}}
{\bf 1}_{\{i_1=i_5\}}
\zeta_{j_2}^{(i_2)}
\zeta_{j_3}^{(i_3)}
\zeta_{j_4}^{(i_4)}
\zeta_{j_6}^{(i_6)}-
{\bf 1}_{\{j_2=j_3\}}
{\bf 1}_{\{i_2=i_3\}}
\zeta_{j_1}^{(i_1)}
\zeta_{j_4}^{(i_4)}
\zeta_{j_5}^{(i_5)}
\zeta_{j_6}^{(i_6)}-
$$
$$
-
{\bf 1}_{\{j_2=j_4\}}
{\bf 1}_{\{i_2=i_4\}}
\zeta_{j_1}^{(i_1)}
\zeta_{j_3}^{(i_3)}
\zeta_{j_5}^{(i_5)}
\zeta_{j_6}^{(i_6)}-
{\bf 1}_{\{j_2=j_5\}}
{\bf 1}_{\{i_2=i_5\}}
\zeta_{j_1}^{(i_1)}
\zeta_{j_3}^{(i_3)}
\zeta_{j_4}^{(i_4)}
\zeta_{j_6}^{(i_6)}-
$$
$$
-
{\bf 1}_{\{j_3=j_4\}}
{\bf 1}_{\{i_3=i_4\}}
\zeta_{j_1}^{(i_1)}
\zeta_{j_2}^{(i_2)}
\zeta_{j_5}^{(i_5)}
\zeta_{j_6}^{(i_6)}-
{\bf 1}_{\{j_3=j_5\}}
{\bf 1}_{\{i_3=i_5\}}
\zeta_{j_1}^{(i_1)}
\zeta_{j_2}^{(i_2)}
\zeta_{j_4}^{(i_4)}
\zeta_{j_6}^{(i_6)}-
$$
$$
-
{\bf 1}_{\{j_4=j_5\}}
{\bf 1}_{\{i_4=i_5\}}
\zeta_{j_1}^{(i_1)}
\zeta_{j_2}^{(i_2)}
\zeta_{j_3}^{(i_3)}
\zeta_{j_6}^{(i_6)}+
$$
$$
+
{\bf 1}_{\{j_1=j_2\}}
{\bf 1}_{\{i_1=i_2\}}
{\bf 1}_{\{j_3=j_4\}}
{\bf 1}_{\{i_3=i_4\}}
\zeta_{j_5}^{(i_5)}
\zeta_{j_6}^{(i_6)}
+
{\bf 1}_{\{j_1=j_2\}}
{\bf 1}_{\{i_1=i_2\}}
{\bf 1}_{\{j_3=j_5\}}
{\bf 1}_{\{i_3=i_5\}}
\zeta_{j_4}^{(i_4)}
\zeta_{j_6}^{(i_6)}+
$$
$$
+
{\bf 1}_{\{j_1=j_2\}}
{\bf 1}_{\{i_1=i_2\}}
{\bf 1}_{\{j_4=j_5\}}
{\bf 1}_{\{i_4=i_5\}}
\zeta_{j_3}^{(i_3)}
\zeta_{j_6}^{(i_6)}
+
{\bf 1}_{\{j_1=j_3\}}
{\bf 1}_{\{i_1=i_3\}}
{\bf 1}_{\{j_2=j_4\}}
{\bf 1}_{\{i_2=i_4\}}
\zeta_{j_5}^{(i_5)}
\zeta_{j_6}^{(i_6)}+
$$
$$
+
{\bf 1}_{\{j_1=j_3\}}
{\bf 1}_{\{i_1=i_3\}}
{\bf 1}_{\{j_2=j_5\}}
{\bf 1}_{\{i_2=i_5\}}
\zeta_{j_4}^{(i_4)}
\zeta_{j_6}^{(i_6)}
+
{\bf 1}_{\{j_1=j_3\}}
{\bf 1}_{\{i_1=i_3\}}
{\bf 1}_{\{j_4=j_5\}}
{\bf 1}_{\{i_4=i_5\}}
\zeta_{j_2}^{(i_2)}
\zeta_{j_6}^{(i_6)}+
$$
$$
+
{\bf 1}_{\{j_1=j_4\}}
{\bf 1}_{\{i_1=i_4\}}
{\bf 1}_{\{j_2=j_3\}}
{\bf 1}_{\{i_2=i_3\}}
\zeta_{j_5}^{(i_5)}
\zeta_{j_6}^{(i_6)}
+
{\bf 1}_{\{j_1=j_4\}}
{\bf 1}_{\{i_1=i_4\}}
{\bf 1}_{\{j_2=j_5\}}
{\bf 1}_{\{i_2=i_5\}}
\zeta_{j_3}^{(i_3)}
\zeta_{j_6}^{(i_6)}+
$$
$$
+
{\bf 1}_{\{j_1=j_4\}}
{\bf 1}_{\{i_1=i_4\}}
{\bf 1}_{\{j_3=j_5\}}
{\bf 1}_{\{i_3=i_5\}}
\zeta_{j_2}^{(i_2)}
\zeta_{j_6}^{(i_6)}
+
{\bf 1}_{\{j_1=j_5\}}
{\bf 1}_{\{i_1=i_5\}}
{\bf 1}_{\{j_2=j_3\}}
{\bf 1}_{\{i_2=i_3\}}
\zeta_{j_4}^{(i_4)}
\zeta_{j_6}^{(i_6)}+
$$
$$
+
{\bf 1}_{\{j_1=j_5\}}
{\bf 1}_{\{i_1=i_5\}}
{\bf 1}_{\{j_2=j_4\}}
{\bf 1}_{\{i_2=i_4\}}
\zeta_{j_3}^{(i_3)}
\zeta_{j_6}^{(i_6)}
+
{\bf 1}_{\{j_1=j_5\}}
{\bf 1}_{\{i_1=i_5\}}
{\bf 1}_{\{j_3=j_4\}}
{\bf 1}_{\{i_3=i_4\}}
\zeta_{j_2}^{(i_2)}
\zeta_{j_6}^{(i_6)}+
$$
$$
+
{\bf 1}_{\{j_2=j_3\}}
{\bf 1}_{\{i_2=i_3\}}
{\bf 1}_{\{j_4=j_5\}}
{\bf 1}_{\{i_4=i_5\}}
\zeta_{j_1}^{(i_1)}
\zeta_{j_6}^{(i_6)}
+
{\bf 1}_{\{j_2=j_4\}}
{\bf 1}_{\{i_2=i_4\}}
{\bf 1}_{\{j_3=j_5\}}
{\bf 1}_{\{i_3=i_5\}}
\zeta_{j_1}^{(i_1)}
\zeta_{j_6}^{(i_6)}+
$$
$$
+
{\bf 1}_{\{j_2=j_5\}}
{\bf 1}_{\{i_2=i_5\}}
{\bf 1}_{\{j_3=j_4\}}
{\bf 1}_{\{i_3=i_4\}}
\zeta_{j_1}^{(i_1)}
\zeta_{j_6}^{(i_6)}
+
{\bf 1}_{\{j_6=j_1\}}
{\bf 1}_{\{i_6=i_1\}}
{\bf 1}_{\{j_3=j_4\}}
{\bf 1}_{\{i_3=i_4\}}
\zeta_{j_2}^{(i_2)}
\zeta_{j_5}^{(i_5)}+
$$
$$
+
{\bf 1}_{\{j_6=j_1\}}
{\bf 1}_{\{i_6=i_1\}}
{\bf 1}_{\{j_3=j_5\}}
{\bf 1}_{\{i_3=i_5\}}
\zeta_{j_2}^{(i_2)}
\zeta_{j_4}^{(i_4)}
+
{\bf 1}_{\{j_6=j_1\}}
{\bf 1}_{\{i_6=i_1\}}
{\bf 1}_{\{j_2=j_5\}}
{\bf 1}_{\{i_2=i_5\}}
\zeta_{j_3}^{(i_3)}
\zeta_{j_4}^{(i_4)}+
$$
$$
+
{\bf 1}_{\{j_6=j_1\}}
{\bf 1}_{\{i_6=i_1\}}
{\bf 1}_{\{j_2=j_4\}}
{\bf 1}_{\{i_2=i_4\}}
\zeta_{j_3}^{(i_3)}
\zeta_{j_5}^{(i_5)}
+
{\bf 1}_{\{j_6=j_1\}}
{\bf 1}_{\{i_6=i_1\}}
{\bf 1}_{\{j_4=j_5\}}
{\bf 1}_{\{i_4=i_5\}}
\zeta_{j_2}^{(i_2)}
\zeta_{j_3}^{(i_3)}+
$$
$$
+
{\bf 1}_{\{j_6=j_1\}}
{\bf 1}_{\{i_6=i_1\}}
{\bf 1}_{\{j_2=j_3\}}
{\bf 1}_{\{i_2=i_3\}}
\zeta_{j_4}^{(i_4)}
\zeta_{j_5}^{(i_5)}
+
{\bf 1}_{\{j_6=j_2\}}
{\bf 1}_{\{i_6=i_2\}}
{\bf 1}_{\{j_3=j_5\}}
{\bf 1}_{\{i_3=i_5\}}
\zeta_{j_1}^{(i_1)}
\zeta_{j_4}^{(i_4)}+
$$
$$
+
{\bf 1}_{\{j_6=j_2\}}
{\bf 1}_{\{i_6=i_2\}}
{\bf 1}_{\{j_4=j_5\}}
{\bf 1}_{\{i_4=i_5\}}
\zeta_{j_1}^{(i_1)}
\zeta_{j_3}^{(i_3)}
+
{\bf 1}_{\{j_6=j_2\}}
{\bf 1}_{\{i_6=i_2\}}
{\bf 1}_{\{j_3=j_4\}}
{\bf 1}_{\{i_3=i_4\}}
\zeta_{j_1}^{(i_1)}
\zeta_{j_5}^{(i_5)}+
$$
$$
+
{\bf 1}_{\{j_6=j_2\}}
{\bf 1}_{\{i_6=i_2\}}
{\bf 1}_{\{j_1=j_5\}}
{\bf 1}_{\{i_1=i_5\}}
\zeta_{j_3}^{(i_3)}
\zeta_{j_4}^{(i_4)}
+
{\bf 1}_{\{j_6=j_2\}}
{\bf 1}_{\{i_6=i_2\}}
{\bf 1}_{\{j_1=j_4\}}
{\bf 1}_{\{i_1=i_4\}}
\zeta_{j_3}^{(i_3)}
\zeta_{j_5}^{(i_5)}+
$$
$$
+
{\bf 1}_{\{j_6=j_2\}}
{\bf 1}_{\{i_6=i_2\}}
{\bf 1}_{\{j_1=j_3\}}
{\bf 1}_{\{i_1=i_3\}}
\zeta_{j_4}^{(i_4)}
\zeta_{j_5}^{(i_5)}
+
{\bf 1}_{\{j_6=j_3\}}
{\bf 1}_{\{i_6=i_3\}}
{\bf 1}_{\{j_2=j_5\}}
{\bf 1}_{\{i_2=i_5\}}
\zeta_{j_1}^{(i_1)}
\zeta_{j_4}^{(i_4)}+
$$
$$
+
{\bf 1}_{\{j_6=j_3\}}
{\bf 1}_{\{i_6=i_3\}}
{\bf 1}_{\{j_4=j_5\}}
{\bf 1}_{\{i_4=i_5\}}
\zeta_{j_1}^{(i_1)}
\zeta_{j_2}^{(i_2)}
+
{\bf 1}_{\{j_6=j_3\}}
{\bf 1}_{\{i_6=i_3\}}
{\bf 1}_{\{j_2=j_4\}}
{\bf 1}_{\{i_2=i_4\}}
\zeta_{j_1}^{(i_1)}
\zeta_{j_5}^{(i_5)}+
$$
$$
+
{\bf 1}_{\{j_6=j_3\}}
{\bf 1}_{\{i_6=i_3\}}
{\bf 1}_{\{j_1=j_5\}}
{\bf 1}_{\{i_1=i_5\}}
\zeta_{j_2}^{(i_2)}
\zeta_{j_4}^{(i_4)}
+
{\bf 1}_{\{j_6=j_3\}}
{\bf 1}_{\{i_6=i_3\}}
{\bf 1}_{\{j_1=j_4\}}
{\bf 1}_{\{i_1=i_4\}}
\zeta_{j_2}^{(i_2)}
\zeta_{j_5}^{(i_5)}+
$$
$$
+
{\bf 1}_{\{j_6=j_3\}}
{\bf 1}_{\{i_6=i_3\}}
{\bf 1}_{\{j_1=j_2\}}
{\bf 1}_{\{i_1=i_2\}}
\zeta_{j_4}^{(i_4)}
\zeta_{j_5}^{(i_5)}
+
{\bf 1}_{\{j_6=j_4\}}
{\bf 1}_{\{i_6=i_4\}}
{\bf 1}_{\{j_3=j_5\}}
{\bf 1}_{\{i_3=i_5\}}
\zeta_{j_1}^{(i_1)}
\zeta_{j_2}^{(i_2)}+
$$
$$
+
{\bf 1}_{\{j_6=j_4\}}
{\bf 1}_{\{i_6=i_4\}}
{\bf 1}_{\{j_2=j_5\}}
{\bf 1}_{\{i_2=i_5\}}
\zeta_{j_1}^{(i_1)}
\zeta_{j_3}^{(i_3)}
+
{\bf 1}_{\{j_6=j_4\}}
{\bf 1}_{\{i_6=i_4\}}
{\bf 1}_{\{j_2=j_3\}}
{\bf 1}_{\{i_2=i_3\}}
\zeta_{j_1}^{(i_1)}
\zeta_{j_5}^{(i_5)}+
$$
$$
+
{\bf 1}_{\{j_6=j_4\}}
{\bf 1}_{\{i_6=i_4\}}
{\bf 1}_{\{j_1=j_5\}}
{\bf 1}_{\{i_1=i_5\}}
\zeta_{j_2}^{(i_2)}
\zeta_{j_3}^{(i_3)}
+
{\bf 1}_{\{j_6=j_4\}}
{\bf 1}_{\{i_6=i_4\}}
{\bf 1}_{\{j_1=j_3\}}
{\bf 1}_{\{i_1=i_3\}}
\zeta_{j_2}^{(i_2)}
\zeta_{j_5}^{(i_5)}+
$$
$$
+
{\bf 1}_{\{j_6=j_4\}}
{\bf 1}_{\{i_6=i_4\}}
{\bf 1}_{\{j_1=j_2\}}
{\bf 1}_{\{i_1=i_2\}}
\zeta_{j_3}^{(i_3)}
\zeta_{j_5}^{(i_5)}
+
{\bf 1}_{\{j_6=j_5\}}
{\bf 1}_{\{i_6=i_5\}}
{\bf 1}_{\{j_3=j_4\}}
{\bf 1}_{\{i_3=i_4\}}
\zeta_{j_1}^{(i_1)}
\zeta_{j_2}^{(i_2)}+
$$
$$
+
{\bf 1}_{\{j_6=j_5\}}
{\bf 1}_{\{i_6=i_5\}}
{\bf 1}_{\{j_2=j_4\}}
{\bf 1}_{\{i_2=i_4\}}
\zeta_{j_1}^{(i_1)}
\zeta_{j_3}^{(i_3)}
+
{\bf 1}_{\{j_6=j_5\}}
{\bf 1}_{\{i_6=i_5\}}
{\bf 1}_{\{j_2=j_3\}}
{\bf 1}_{\{i_2=i_3\}}
\zeta_{j_1}^{(i_1)}
\zeta_{j_4}^{(i_4)}+
$$
$$
+
{\bf 1}_{\{j_6=j_5\}}
{\bf 1}_{\{i_6=i_5\}}
{\bf 1}_{\{j_1=j_4\}}
{\bf 1}_{\{i_1=i_4\}}
\zeta_{j_2}^{(i_2)}
\zeta_{j_3}^{(i_3)}
+
{\bf 1}_{\{j_6=j_5\}}
{\bf 1}_{\{i_6=i_5\}}
{\bf 1}_{\{j_1=j_3\}}
{\bf 1}_{\{i_1=i_3\}}
\zeta_{j_2}^{(i_2)}
\zeta_{j_4}^{(i_4)}+
$$
$$
+
{\bf 1}_{\{j_6=j_5\}}
{\bf 1}_{\{i_6=i_5\}}
{\bf 1}_{\{j_1=j_2\}}
{\bf 1}_{\{i_1=i_2\}}
\zeta_{j_3}^{(i_3)}
\zeta_{j_4}^{(i_4)}-
$$
$$
-
{\bf 1}_{\{j_6=j_1\}}
{\bf 1}_{\{i_6=i_1\}}
{\bf 1}_{\{j_2=j_5\}}
{\bf 1}_{\{i_2=i_5\}}
{\bf 1}_{\{j_3=j_4\}}
{\bf 1}_{\{i_3=i_4\}}-
$$
$$
-
{\bf 1}_{\{j_6=j_1\}}
{\bf 1}_{\{i_6=i_1\}}
{\bf 1}_{\{j_2=j_4\}}
{\bf 1}_{\{i_2=i_4\}}
{\bf 1}_{\{j_3=j_5\}}
{\bf 1}_{\{i_3=i_5\}}-
$$
$$
-
{\bf 1}_{\{j_6=j_1\}}
{\bf 1}_{\{i_6=i_1\}}
{\bf 1}_{\{j_2=j_3\}}
{\bf 1}_{\{i_2=i_3\}}
{\bf 1}_{\{j_4=j_5\}}
{\bf 1}_{\{i_4=i_5\}}-
$$
$$
-               
{\bf 1}_{\{j_6=j_2\}}
{\bf 1}_{\{i_6=i_2\}}
{\bf 1}_{\{j_1=j_5\}}
{\bf 1}_{\{i_1=i_5\}}
{\bf 1}_{\{j_3=j_4\}}
{\bf 1}_{\{i_3=i_4\}}-
$$
$$
-
{\bf 1}_{\{j_6=j_2\}}
{\bf 1}_{\{i_6=i_2\}}
{\bf 1}_{\{j_1=j_4\}}
{\bf 1}_{\{i_1=i_4\}}
{\bf 1}_{\{j_3=j_5\}}
{\bf 1}_{\{i_3=i_5\}}-
$$
$$
-
{\bf 1}_{\{j_6=j_2\}}
{\bf 1}_{\{i_6=i_2\}}
{\bf 1}_{\{j_1=j_3\}}
{\bf 1}_{\{i_1=i_3\}}
{\bf 1}_{\{j_4=j_5\}}
{\bf 1}_{\{i_4=i_5\}}-
$$
$$
-
{\bf 1}_{\{j_6=j_3\}}
{\bf 1}_{\{i_6=i_3\}}
{\bf 1}_{\{j_1=j_5\}}
{\bf 1}_{\{i_1=i_5\}}
{\bf 1}_{\{j_2=j_4\}}
{\bf 1}_{\{i_2=i_4\}}-
$$
$$
-
{\bf 1}_{\{j_6=j_3\}}
{\bf 1}_{\{i_6=i_3\}}
{\bf 1}_{\{j_1=j_4\}}
{\bf 1}_{\{i_1=i_4\}}
{\bf 1}_{\{j_2=j_5\}}
{\bf 1}_{\{i_2=i_5\}}-
$$
$$
-
{\bf 1}_{\{j_3=j_6\}}
{\bf 1}_{\{i_3=i_6\}}
{\bf 1}_{\{j_1=j_2\}}
{\bf 1}_{\{i_1=i_2\}}
{\bf 1}_{\{j_4=j_5\}}
{\bf 1}_{\{i_4=i_5\}}-
$$
$$
-
{\bf 1}_{\{j_6=j_4\}}
{\bf 1}_{\{i_6=i_4\}}
{\bf 1}_{\{j_1=j_5\}}
{\bf 1}_{\{i_1=i_5\}}
{\bf 1}_{\{j_2=j_3\}}
{\bf 1}_{\{i_2=i_3\}}-
$$
$$
-
{\bf 1}_{\{j_6=j_4\}}
{\bf 1}_{\{i_6=i_4\}}
{\bf 1}_{\{j_1=j_3\}}
{\bf 1}_{\{i_1=i_3\}}
{\bf 1}_{\{j_2=j_5\}}
{\bf 1}_{\{i_2=i_5\}}-
$$
$$
-
{\bf 1}_{\{j_6=j_4\}}
{\bf 1}_{\{i_6=i_4\}}
{\bf 1}_{\{j_1=j_2\}}
{\bf 1}_{\{i_1=i_2\}}
{\bf 1}_{\{j_3=j_5\}}
{\bf 1}_{\{i_3=i_5\}}-
$$
$$
-
{\bf 1}_{\{j_6=j_5\}}
{\bf 1}_{\{i_6=i_5\}}
{\bf 1}_{\{j_1=j_4\}}
{\bf 1}_{\{i_1=i_4\}}
{\bf 1}_{\{j_2=j_3\}}
{\bf 1}_{\{i_2=i_3\}}-
$$
$$
-
{\bf 1}_{\{j_6=j_5\}}
{\bf 1}_{\{i_6=i_5\}}
{\bf 1}_{\{j_1=j_2\}}
{\bf 1}_{\{i_1=i_2\}}
{\bf 1}_{\{j_3=j_4\}}
{\bf 1}_{\{i_3=i_4\}}-
$$
$$
\Biggl.-
{\bf 1}_{\{j_6=j_5\}}
{\bf 1}_{\{i_6=i_5\}}
{\bf 1}_{\{j_1=j_3\}}
{\bf 1}_{\{i_1=i_3\}}
{\bf 1}_{\{j_2=j_4\}}
{\bf 1}_{\{i_2=i_4\}}\Biggr),
$$

\vspace{7mm}
\noindent
where

$$
C_{j_3j_2j_1}=\int\limits_{\tau_p}^{\tau_{p+1}}\phi_{j_3}(z)
\int\limits_{\tau_p}^{z}\phi_{j_2}(y)
\int\limits_{\tau_p}^{y}
\phi_{j_1}(x)dx dy dz=
$$

$$
=
\frac{\sqrt{(2j_1+1)(2j_2+1)(2j_3+1)}}{8}\Delta^{3/2}\bar
C_{j_3j_2j_1},
$$

\vspace{4mm}

$$
C_{j_4j_3j_2j_1}=\int\limits_{\tau_p}^{\tau_{p+1}}\phi_{j_4}(u)
\int\limits_{\tau_p}^{u}\phi_{j_3}(z)
\int\limits_{\tau_p}^{z}\phi_{j_2}(y)
\int\limits_{\tau_p}^{y}
\phi_{j_1}(x)dx dy dz du=
$$

$$
=\frac{\sqrt{(2j_1+1)(2j_2+1)(2j_3+1)(2j_4+1)}}{16}\Delta^{2}\bar
C_{j_4j_3j_2j_1},
$$

\vspace{4mm}

$$
C_{j_3j_2j_1}^{001}=\int\limits_{\tau_p}^{\tau_{p+1}}(\tau_p-z)\phi_{j_3}(z)
\int\limits_{\tau_p}^{z}\phi_{j_2}(y)
\int\limits_{\tau_p}^{y}
\phi_{j_1}(x)dx dy dz=
$$

$$
=
\frac{\sqrt{(2j_1+1)(2j_2+1)(2j_3+1)}}{16}\Delta^{5/2}\bar
C_{j_3j_2j_1}^{001},
$$

\vspace{4mm}

$$
C_{j_3j_2j_1}^{010}=\int\limits_{\tau_p}^{\tau_{p+1}}\phi_{j_3}(z)
\int\limits_{\tau_p}^{z}(\tau_p-y)\phi_{j_2}(y)
\int\limits_{\tau_p}^{y}
\phi_{j_1}(x)dx dy dz=
$$

$$
=
\frac{\sqrt{(2j_1+1)(2j_2+1)(2j_3+1)}}{16}\Delta^{5/2}\bar
C_{j_3j_2j_1}^{010},
$$

\vspace{4mm}

$$
C_{j_3j_2j_1}^{100}=\int\limits_{\tau_p}^{\tau_{p+1}}\phi_{j_3}(z)
\int\limits_{\tau_p}^{z}\phi_{j_2}(y)
\int\limits_{\tau_p}^{y}
(\tau_p-x)\phi_{j_1}(x)dx dy dz=
$$

$$
=
\frac{\sqrt{(2j_1+1)(2j_2+1)(2j_3+1)}}{16}\Delta^{5/2}\bar
C_{j_3j_2j_1}^{100},
$$

\vspace{4mm}

$$
C_{j_5j_4 j_3 j_2 j_1}=
\int\limits_{\tau_p}^{\tau_{p+1}}\phi_{j_5}(v)
\int\limits_{\tau_p}^v\phi_{j_4}(u)
\int\limits_{\tau_p}^{u}
\phi_{j_3}(z)
\int\limits_{\tau_p}^{z}\phi_{j_2}(y)\int\limits_{\tau_p}^{y}\phi_{j_1}(x)
dxdydzdudv=
$$

$$
=\frac{\sqrt{(2j_1+1)(2j_2+1)(2j_3+1)(2j_4+1)(2j_5+1)}}{32}\Delta^{5/2}\bar
C_{j_5j_4 j_3 j_2 j_1},
$$

\vspace{4mm}

$$
C_{j_4j_3j_2j_1}^{0001}=\int\limits_{\tau_p}^{\tau_{p+1}}(\tau_p-u)\phi_{j_4}(u)
\int\limits_{\tau_p}^{u}\phi_{j_3}(z)
\int\limits_{\tau_p}^{z}\phi_{j_2}(y)
\int\limits_{\tau_p}^{y}
\phi_{j_1}(x)dx dy dz du=
$$

$$
=\frac{\sqrt{(2j_1+1)(2j_2+1)(2j_3+1)(2j_4+1)}}{32}\Delta^{3}\bar
C_{j_4j_3j_2j_1}^{0001},
$$

\vspace{4mm}

$$
C_{j_3j_2j_1}^{0010}=\int\limits_{\tau_p}^{\tau_{p+1}}\phi_{j_4}(u)
\int\limits_{\tau_p}^{u}(\tau_p-z)\phi_{j_3}(z)
\int\limits_{\tau_p}^{z}\phi_{j_2}(y)
\int\limits_{\tau_p}^{y}
\phi_{j_1}(x)dx dy dz du=
$$

$$
=\frac{\sqrt{(2j_1+1)(2j_2+1)(2j_3+1)(2j_4+1)}}{32}\Delta^{3}\bar
C_{j_4j_3j_2j_1}^{0010},
$$

\vspace{4mm}

$$
C_{j_4j_3j_2j_1}^{0100}=
\int\limits_{\tau_p}^{\tau_{p+1}}\phi_{j_4}(u)
\int\limits_{\tau_p}^{u}\phi_{j_3}(z)
\int\limits_{\tau_p}^{z}(\tau_p-y)\phi_{j_2}(y)
\int\limits_{\tau_p}^{y}
\phi_{j_1}(x)dx dy dz du=
$$

$$
=\frac{\sqrt{(2j_1+1)(2j_2+1)(2j_3+1)(2j_4+1)}}{32}\Delta^{3}\bar
C_{j_3j_2j_1}^{0100},
$$

\vspace{4mm}

$$
C_{j_4j_3j_2j_1}^{1000}=\int\limits_{\tau_p}^{\tau_{p+1}}\phi_{j_4}(u)
\int\limits_{\tau_p}^{u}\phi_{j_3}(z)
\int\limits_{\tau_p}^{z}
\phi_{j_2}(y)
\int\limits_{\tau_p}^{y}
(\tau_p-x)\phi_{j_1}(x)
dx dy dz du=
$$

$$
=\frac{\sqrt{(2j_1+1)(2j_2+1)(2j_3+1)(2j_4+1)}}{32}\Delta^{3}\bar
C_{j_4j_3j_2j_1}^{1000},
$$

\vspace{4mm}

$$
C_{j_6j_5j_4 j_3 j_2 j_1}=
\int\limits_{\tau_p}^{\tau_{p+1}}\phi_{j_6}(w)
\int\limits_{\tau_p}^w\phi_{j_5}(v)
\int\limits_{\tau_p}^v\phi_{j_4}(u)
\int\limits_{\tau_p}^{u}
\phi_{j_3}(z)
\int\limits_{\tau_p}^{z}\phi_{j_2}(y)\int\limits_{\tau_p}^{y}\phi_{j_1}(x)
dxdydzdudvdw=
$$

$$
=\frac{\sqrt{(2j_1+1)(2j_2+1)(2j_3+1)(2j_4+1)(2j_5+1)(2j_6+1)}}{64}\Delta^{3}\bar
C_{j_6j_5j_4 j_3 j_2 j_1},
$$

\vspace{5mm}
\noindent
where

\vspace{-2mm}
$$
\bar C_{j_3j_2j_1}=
\int\limits_{-1}^{1}P_{j_3}(z)
\int\limits_{-1}^{z}P_{j_2}(y)
\int\limits_{-1}^{y}
P_{j_1}(x)dx dy dz,
$$

\vspace{1mm}
$$
\bar C_{j_4j_3j_2j_1}=
\int\limits_{-1}^{1}P_{j_4}(u)
\int\limits_{-1}^{u}P_{j_3}(z)
\int\limits_{-1}^{z}P_{j_2}(y)
\int\limits_{-1}^{y}
P_{j_1}(x)dx dy dz,
$$

\vspace{1mm}

$$
\bar C_{j_3j_2j_1}^{100}=-
\int\limits_{-1}^{1}P_{j_3}(z)
\int\limits_{-1}^{z}P_{j_2}(y)
\int\limits_{-1}^{y}
P_{j_1}(x)(x+1)dx dy dz,
$$

\vspace{1mm}
$$
\bar C_{j_3j_2j_1}^{010}=-
\int\limits_{-1}^{1}P_{j_3}(z)
\int\limits_{-1}^{z}P_{j_2}(y)(y+1)
\int\limits_{-1}^{y}
P_{j_1}(x)dx dy dz,
$$

\vspace{1mm}

$$
\bar C_{j_3j_2j_1}^{001}=-
\int\limits_{-1}^{1}P_{j_3}(z)(z+1)
\int\limits_{-1}^{z}P_{j_2}(y)
\int\limits_{-1}^{y}
P_{j_1}(x)dx dy dz,
$$

\vspace{1mm}

$$
\bar C_{j_5j_4 j_3 j_2 j_1}=
\int\limits_{-1}^{1}P_{j_5}(v)
\int\limits_{-1}^{v}P_{j_4}(u)
\int\limits_{-1}^{u}P_{j_3}(z)
\int\limits_{-1}^{z}P_{j_2}(y)
\int\limits_{-1}^{y}
P_{j_1}(x)dx dy dz du dv,
$$

\vspace{1mm}

$$
\bar C_{j_4j_3j_2j_1}^{1000}=-
\int\limits_{-1}^{1}P_{j_4}(u)
\int\limits_{-1}^{u}P_{j_3}(z)
\int\limits_{-1}^{z}P_{j_2}(y)
\int\limits_{-1}^{y}
P_{j_1}(x)(x+1)dx dy dz du,
$$

\vspace{1mm}

$$
\bar C_{j_4j_3j_2j_1}^{0100}=-
\int\limits_{-1}^{1}P_{j_4}(u)
\int\limits_{-1}^{u}P_{j_3}(z)
\int\limits_{-1}^{z}P_{j_2}(y)(y+1)
\int\limits_{-1}^{y}
P_{j_1}(x)dx dy dz du,
$$

\vspace{1mm}

$$
\bar C_{j_4j_3j_2j_1}^{0010}=-
\int\limits_{-1}^{1}P_{j_4}(u)
\int\limits_{-1}^{u}P_{j_3}(z)(z+1)
\int\limits_{-1}^{z}P_{j_2}(y)
\int\limits_{-1}^{y}
P_{j_1}(x)dx dy dz du,
$$

\vspace{1mm}

$$
\bar C_{j_4j_3j_2j_1}^{0001}=-
\int\limits_{-1}^{1}P_{j_4}(u)(u+1)
\int\limits_{-1}^{u}P_{j_3}(z)
\int\limits_{-1}^{z}P_{j_2}(y)
\int\limits_{-1}^{y}
P_{j_1}(x)dx dy dz du,
$$

\vspace{1mm}

$$
\bar C_{j_6j_5j_4 j_3 j_2 j_1}=
\int\limits_{-1}^{1}P_{j_6}(w)
\int\limits_{-1}^{w}P_{j_5}(v)
\int\limits_{-1}^{v}P_{j_4}(u)
\int\limits_{-1}^{u}P_{j_3}(z)
\int\limits_{-1}^{z}P_{j_2}(y)
\int\limits_{-1}^{y}
P_{j_1}(x)dx dy dz du dv dw,
$$

\vspace{5mm}
\noindent
where $P_i(x)$ $(i=0, 1, 2,\ldots)$ is the Legendre polynomial and

\vspace{2mm}

$$
\phi_i(x)=\sqrt{\frac{2i+1}{\Delta}}P_i\left(\left(x-\tau_p-\frac{\Delta}{2}\right)
\frac{2}{\Delta}\right),\ \ \ i=0, 1, 2,\ldots 
$$

\vspace{4mm}

Let us consider the exact relations and some estimates
for the mean-square errors of appro\-xi\-ma\-ti\-ons of iterated Ito
stochastic integrals.

Using Theorem 3, we get 
\cite{2017}-\cite{2013}, \cite{arxiv-2} (also see \cite{2006}, \cite{kuz1997}, 
\cite{kuz1997a}, \cite{5-000}-\cite{5-008}, \cite{2018a}-\cite{2018aaa})
 
\vspace{2mm}
$$
{\sf M}\left\{\left(I_{{00}_{\tau_{p+1},\tau_p}}^{(i_1 i_2)}-
I_{{00}_{\tau_{p+1},\tau_p}}^{(i_1 i_2)q}
\right)^2\right\}
=\frac{\Delta^2}{2}\Biggl(\frac{1}{2}-\sum_{i=1}^q
\frac{1}{4i^2-1}\Biggr)\ \ \ (i_1\ne i_2),
$$

\vspace{6mm}

$$
{\sf M}\left\{\left(I_{{10}_{\tau_{p+1},\tau_p}}^{(i_1 i_2)}-
I_{{10}_{\tau_{p+1},\tau_p}}^{(i_1 i_2)q}
\right)^2\right\}=
{\sf M}\left\{\left(I_{{01}_{\tau_{p+1},\tau_p}}^{(i_1 i_2)}-
I_{{01}_{\tau_{p+1},\tau_p}}^{(i_1 i_2)q}\right)^2\right\}=
$$

$$
=\frac{\Delta^4}{16}\Biggl(\frac{5}{9}-
2\sum_{i=2}^q\frac{1}{4i^2-1}-
\sum_{i=1}^q
\frac{1}{(2i-1)^2(2i+3)^2}
-\sum_{i=0}^q\frac{(i+2)^2+(i+1)^2}{(2i+1)(2i+5)(2i+3)^2}
\Biggr)\ \ \ (i_1\ne i_2),
$$

\vspace{6mm}

$$
{\sf M}\left\{\left(I_{{10}_{\tau_{p+1},\tau_p}}^{(i_1 i_1)}-
I_{{10}_{\tau_{p+1},\tau_p}}^{(i_1 i_1)q}
\right)^2\right\}=
{\sf M}\left\{\left(I_{{01}_{\tau_{p+1},\tau_p}}^{(i_1 i_1)}-
I_{{01}_{\tau_{p+1},\tau_p}}^{(i_1 i_1)q}\right)^2\right\}=
$$

$$
=\frac{\Delta^4}{16}\Biggl(\frac{1}{9}-
\sum_{i=0}^{q}
\frac{1}{(2i+1)(2i+5)(2i+3)^2}
-2\sum_{i=1}^{q}
\frac{1}{(2i-1)^2(2i+3)^2}\Biggr).
$$

\vspace{6mm}

Applying (\ref{qq1}), (\ref{qq2})--(\ref{883}), we obtain

\vspace{2mm}

$$
{\sf M}\biggl\{\left(I_{{20}_{\tau_{p+1},\tau_p}}^{(i_1 i_2)}-
I_{{20}_{\tau_{p+1},\tau_p}}^{(i_1 i_2)q}
\right)^2\biggr\}=\frac{\Delta^6}{30}
-\sum_{j_1,j_2=0}^q
\left(C_{j_2j_1}^{20}\right)^2-
\sum_{j_1,j_2=0}^q C_{j_2j_1}^{20}C_{j_1j_2}^{20}\ \ \ (i_1=i_2),
$$

\vspace{3mm}

$$
{\sf M}\biggl\{\left(I_{{20}_{\tau_{p+1},\tau_p}}^{(i_1 i_2)}-
I_{{20}_{\tau_{p+1},\tau_p}}^{(i_1 i_2)q}
\right)^2\biggr\}=\frac{\Delta^6}{30}
-\sum_{j_1,j_2=0}^q
\left(C_{j_2j_1}^{20}\right)^2\ \ \ (i_1\ne i_2),
$$

\vspace{3mm}

$$
{\sf M}\biggl\{\left(I_{{11}_{\tau_{p+1},\tau_p}}^{(i_1 i_2)}-
I_{{11}_{\tau_{p+1},\tau_p}}^{(i_1 i_2)q}
\right)^2\biggr\}=\frac{\Delta^6}{18}
-\sum_{j_1,j_2=0}^q
\left(C_{j_2j_1}^{11}\right)^2-
\sum_{j_1,j_2=0}^q C_{j_2j_1}^{11}C_{j_1j_2}^{11}\ \ \ (i_1=i_2),
$$

\vspace{3mm}

$$
{\sf M}\biggl\{\left(I_{{11}_{\tau_{p+1},\tau_p}}^{(i_1 i_2)}-
I_{{11}_{\tau_{p+1},\tau_p}}^{(i_1 i_2)q}
\right)^2\biggr\}=\frac{\Delta^6}{18}
-\sum_{j_1,j_2=0}^q
\left(C_{j_2j_1}^{11}\right)^2\ \ \ (i_1\ne i_2),
$$

\vspace{3mm}

$$
{\sf M}\biggl\{\left(I_{{02}_{\tau_{p+1},\tau_p}}^{(i_1 i_2)}-
I_{{02}_{\tau_{p+1},\tau_p}}^{(i_1 i_2)q}
\right)^2\biggr\}=\frac{\Delta^6}{6}
-\sum_{j_1,j_2=0}^q
\left(C_{j_2j_1}^{02}\right)^2-
\sum_{j_1,j_2=0}^q C_{j_2j_1}^{02}C_{j_1j_2}^{02}\ \ \ (i_1=i_2),
$$

\vspace{3mm}

$$
{\sf M}\biggl\{\left(I_{{02}_{\tau_{p+1},\tau_p}}^{(i_1 i_2)}-
I_{{02}_{\tau_{p+1},\tau_p}}^{(i_1 i_2)q}
\right)^2\biggr\}=\frac{\Delta^6}{6}
-\sum_{j_1,j_2=0}^q
\left(C_{j_2j_1}^{02}\right)^2\ \ \ (i_1\ne i_2),
$$

\vspace{3mm}

$$
{\sf M}\left\{\left(
I_{{000}_{\tau_{p+1},\tau_p}}^{(i_1i_2 i_3)}-
I_{{000}_{\tau_{p+1},\tau_p}}^{(i_1i_2 i_3)q}\right)^2\right\}=
\frac{\Delta^{3}}{6}-\sum_{j_3,j_2,j_1=0}^{q}
C_{j_3j_2j_1}^2\ \ \ (i_1\ne i_2,\  i_1\ne i_3,\ i_2\ne i_3),
$$

\vspace{3mm}

$$
{\sf M}\left\{\left(
I_{{000}_{\tau_{p+1},\tau_p}}^{(i_1i_2 i_3)}-
I_{{000}_{\tau_{p+1},\tau_p}}^{(i_1i_2 i_3)q}\right)^2\right\}=
\frac{\Delta^{3}}{6}-\sum_{j_3,j_2,j_1=0}^{q}
C_{j_3j_2j_1}^2
-\sum_{j_3,j_2,j_1=0}^{q}
C_{j_2j_3j_1}C_{j_3j_2j_1}\ \ \ (i_1\ne i_2=i_3),
$$

\vspace{3mm}

$$
{\sf M}\left\{\left(
I_{{000}_{\tau_{p+1},\tau_p}}^{(i_1i_2 i_3)}-
I_{{000}_{\tau_{p+1},\tau_p}}^{(i_1i_2 i_3)q}\right)^2\right\}=
\frac{\Delta^{3}}{6}-\sum_{j_3,j_2,j_1=0}^{q}
C_{j_3j_2j_1}^2
-\sum_{j_3,j_2,j_1=0}^{q}
C_{j_3j_2j_1}C_{j_1j_2j_3}\ \ \ (i_1=i_3\ne i_2),
$$

\vspace{3mm}

$$
{\sf M}\left\{\left(
I_{{000}_{\tau_{p+1},\tau_p}}^{(i_1i_2 i_3)}-
I_{{000}_{\tau_{p+1},\tau_p}}^{(i_1i_2 i_3)q}\right)^2\right\}=
\frac{\Delta^{3}}{6}-\sum_{j_3,j_2,j_1=0}^{q}
C_{j_3j_2j_1}^2
-\sum_{j_3,j_2,j_1=0}^{q}
C_{j_3j_1j_2}C_{j_3j_2j_1}\ \ \ (i_1=i_2\ne i_3),
$$

\vspace{5mm}
\noindent
where

$$
C_{j_2j_1}^{20}=\int\limits_{\tau_p}^{\tau_{p+1}}\phi_{j_2}(y)
\int\limits_{\tau_p}^{y}
\phi_{j_1}(x)(\tau_p-x)^2 dx dy =
\frac{\sqrt{(2j_1+1)(2j_2+1)}}{16}\Delta^{3}\bar
C_{j_2j_1}^{20},
$$

\vspace{1mm}
$$
C_{j_2j_1}^{02}=\int\limits_{\tau_p}^{\tau_{p+1}}\phi_{j_2}(y)(\tau_p-y)^2
\int\limits_{\tau_p}^{y}
\phi_{j_1}(x)dx dy =
\frac{\sqrt{(2j_1+1)(2j_2+1)}}{16}\Delta^{3}\bar
C_{j_2j_1}^{02},
$$

\vspace{1mm}
$$
C_{j_2j_1}^{11}=\int\limits_{\tau_p}^{\tau_{p+1}}\phi_{j_2}(y)(\tau_p-y)
\int\limits_{\tau_p}^{y}
\phi_{j_1}(x)(\tau_p-x) dx dy =
\frac{\sqrt{(2j_1+1)(2j_2+1)}}{16}\Delta^{3}\bar
C_{j_2j_1}^{11}, 
$$

\vspace{1mm}
$$
\bar C_{j_2j_1}^{20}=
\int\limits_{-1}^{1}P_{j_2}(y)
\int\limits_{-1}^{y}
P_{j_1}(x)(x+1)^2 dx dy,
$$

\vspace{1mm}
$$
\bar C_{j_2j_1}^{02}=
\int\limits_{-1}^{1}P_{j_2}(y)(y+1)^2
\int\limits_{-1}^{y}
P_{j_1}(x)dx dy,
$$

\vspace{1mm}
$$
\bar C_{j_2j_1}^{11}=
\int\limits_{-1}^{1}P_{j_2}(y)(y+1)
\int\limits_{-1}^{y}
P_{j_1}(x)(x+1)dx dy,
$$

\vspace{4mm}
\noindent
where $P_i(x)$ $(i=0, 1, 2,\ldots)$ is the Legendre polynomial and

\vspace{1mm}
$$
\phi_i(x)=\sqrt{\frac{2i+1}{\Delta}}
P_i\left(\left(x-\tau_p-\frac{\Delta}{2}\right)
\frac{2}{\Delta}\right),\ \ \ i=0, 1, 2,\ldots 
$$

\vspace{4mm}

At the same time using the estimate (\ref{qq4}) 
for $i_1,\ldots,i_6=1,\ldots,m$, 
we get

\vspace{1mm}

$$
{\sf M}\left\{\left(
I_{{01}_{\tau_{p+1},\tau_p}}^{(i_1i_2)}-
I_{{01}_{\tau_{p+1},\tau_p}}^{(i_1i_2)q}\right)^2\right\}\le
2\Biggl(\frac{\Delta^{4}}{4}-\sum_{j_1,j_2=0}^{q}
\left(C_{j_2j_1}^{01}\right)^2\Biggr),
$$

\vspace{3mm}
$$
{\sf M}\left\{\left(
I_{{10}_{\tau_{p+1},\tau_p}}^{(i_1i_2)}-
I_{{10}_{\tau_{p+1},\tau_p}}^{(i_1i_2)q}\right)^2\right\}\le
2\Biggl(\frac{\Delta^{4}}{12}-\sum_{j_1,j_2=0}^{q}
\left(C_{j_2j_1}^{10}\right)^2\Biggr),
$$

\vspace{3mm}

$$
{\sf M}\left\{\left(
I_{{000}_{\tau_{p+1},\tau_p}}^{(i_1i_2 i_3)}-
I_{{000}_{\tau_{p+1},\tau_p}}^{(i_1i_2 i_3)q}\right)^2\right\}\le
6\Biggl(\frac{\Delta^{3}}{6}-\sum_{j_3,j_2,j_1=0}^{q}
C_{j_3j_2j_1}^2\Biggr),
$$

\vspace{3mm}

$$
{\sf M}\left\{\left(
I_{{0000}_{\tau_{p+1},\tau_p}}^{(i_1i_2 i_3 i_4)}-
I_{{0000}_{\tau_{p+1},\tau_p}}^{(i_1i_2 i_3 i_4)q}\right)^2\right\}\le
24\Biggl(\frac{\Delta^{4}}{24}-\sum_{j_1,j_2,j_3,j_4=0}^{q}
C_{j_4j_3j_2j_1}^2\Biggr),
$$

\vspace{3mm}

$$
{\sf M}\left\{\left(
I_{{100}_{\tau_{p+1},\tau_p}}^{(i_1i_2 i_3)}-
I_{{100}_{\tau_{p+1},\tau_p}}^{(i_1i_2 i_3)q}\right)^2\right\}\le
6\Biggl(\frac{\Delta^{5}}{60}-\sum_{j_1,j_2,j_3=0}^{q}
\left(C_{j_3j_2j_1}^{100}\right)^2\Biggr),
$$

\vspace{3mm}

$$
{\sf M}\left\{\left(
I_{{010}_{\tau_{p+1},\tau_p}}^{(i_1i_2 i_3)}-
I_{{010}_{\tau_{p+1},\tau_p}}^{(i_1i_2 i_3)q}\right)^2\right\}\le
6\Biggl(\frac{\Delta^{5}}{20}-\sum_{j_1,j_2,j_3=0}^{q}
\left(C_{j_3j_2j_1}^{010}\right)^2\Biggr),
$$
                          
\vspace{3mm}

$$
{\sf M}\left\{\left(
I_{{001}_{\tau_{p+1},\tau_p}}^{(i_1i_2 i_3)}-
I_{{001}_{\tau_{p+1},\tau_p}}^{(i_1i_2 i_3)q}\right)^2\right\}\le
6\Biggl(\frac{\Delta^5}{10}-\sum_{j_1,j_2,j_3=0}^{q}
\left(C_{j_3j_2j_1}^{001}\right)^2\Biggr),
$$

\vspace{3mm}

$$
{\sf M}\left\{\left(
I_{{00000}_{\tau_{p+1},\tau_p}}^{(i_1 i_2 i_3 i_4 i_5)}-
I_{{00000}_{\tau_{p+1},\tau_p}}^{(i_1 i_2 i_3 i_4 i_5)q}\right)^2\right\}\le
120\left(\frac{\Delta^{5}}{120}-\sum_{j_1,j_2,j_3,j_4,j_5=0}^{q}
C_{j_5 i_4 i_3 i_2 j_1}^2\right),
$$

\vspace{3mm}

$$
{\sf M}\left\{\left(
I_{{20}_{\tau_{p+1},\tau_p}}^{(i_1i_2)}-
I_{{20}_{\tau_{p+1},\tau_p}}^{(i_1i_2)q}\right)^2\right\}\le
2\Biggl(\frac{\Delta^6}{30}-\sum_{j_2,j_1=0}^{q}
\left(C_{j_2j_1}^{20}\right)^2\Biggr),
$$

\vspace{3mm}

$$
{\sf M}\left\{\left(
I_{{11}_{\tau_{p+1},\tau_p}}^{(i_1i_2)}-
I_{{11}_{\tau_{p+1},\tau_p}}^{(i_1i_2)q}\right)^2\right\}\le
2\Biggl(\frac{\Delta^6}{18}-\sum_{j_2,j_1=0}^{q}
\left(C_{j_2j_1}^{11}\right)^2\Biggr),
$$

\vspace{3mm}

$$
{\sf M}\left\{\left(
I_{{02}_{\tau_{p+1},\tau_p}}^{(i_1i_2)}-
I_{{02}_{\tau_{p+1},\tau_p}}^{(i_1i_2)q}\right)^2\right\}\le
2\Biggl(\frac{\Delta^6}{6}-\sum_{j_2,j_1=0}^{q}
\left(C_{j_2j_1}^{02}\right)^2\Biggr),
$$

\vspace{3mm}

$$
{\sf M}\left\{\left(
I_{{1000}_{\tau_{p+1},\tau_p}}^{(i_1i_2 i_3i_4)}-
I_{{1000}_{\tau_{p+1},\tau_p}}^{(i_1i_2 i_3i_4)q}\right)^2\right\}\le
24\Biggl(\frac{\Delta^{6}}{360}-\sum_{j_1,j_2,j_3, j_4=0}^{q}
\left(C_{j_4j_3j_2j_1}^{1000}\right)^2\Biggr),
$$

\vspace{3mm}

$$
{\sf M}\left\{\left(
I_{{0100}_{\tau_{p+1},\tau_p}}^{(i_1i_2 i_3i_4)}-
I_{{0100}_{\tau_{p+1},\tau_p}}^{(i_1i_2 i_3i_4)q}\right)^2\right\}\le
24\Biggl(\frac{\Delta^{6}}{120}-\sum_{j_1,j_2,j_3, j_4=0}^{q}
\left(C_{j_4j_3j_2j_1}^{0100}\right)^2\Biggr),
$$

\vspace{3mm}

$$
{\sf M}\left\{\left(
I_{{0010}_{\tau_{p+1},\tau_p}}^{(i_1i_2 i_3i_4)}-
I_{{0010}_{\tau_{p+1},\tau_p}}^{(i_1i_2 i_3 i_4)q}\right)^2\right\}\le
24\Biggl(\frac{\Delta^6}{60}-\sum_{j_1,j_2,j_3, j_4=0}^{q}
\left(C_{j_4j_3j_2j_1}^{0010}\right)^2\Biggr),
$$

\vspace{3mm}

$$
{\sf M}\left\{\left(
I_{{0001}_{\tau_{p+1},\tau_p}}^{(i_1i_2 i_3 i_4)}-
I_{{0001}_{\tau_{p+1},\tau_p}}^{(i_1i_2 i_3 i_4)q}\right)^2\right\}\le
24\Biggl(\frac{\Delta^6}{36}-\sum_{j_1,j_2,j_3, j_4=0}^{q}
\left(C_{j_4j_3j_2j_1}^{0001}\right)^2\Biggr),
$$

\vspace{3mm}

$$
{\sf M}\left\{\left(
I_{{000000}_{\tau_{p+1},\tau_p}}^{(i_1 i_2 i_3 i_4 i_5 i_6)}-
I_{{000000}_{\tau_{p+1},\tau_p}}^{(i_1 i_2 i_3 i_4 i_5 i_6)q}\right)^2\right\}\le
720\left(\frac{\Delta^{6}}{720}-\sum_{j_1,j_2,j_3,j_4,j_5,j_6=0}^{q}
C_{j_6 j_5 j_4 j_3 j_2 j_1}^2\right).
$$

\vspace{5mm}

The Fourier--Legendre coefficients 

$$
\bar C_{j_3 j_2 j_1},\ \bar C_{j_4 j_3 j_2 j_1},\ \bar C_{j_3 j_2 j_1}^{001},\ 
\bar C_{j_3 j_2 j_1}^{010},\ \bar C_{j_3 j_2 j_1}^{100},\
\bar C_{j_5 j_4 j_3 j_2 j_1},\
\bar C_{j_4 j_3 j_2 j_1}^{0001},
$$

$$
\bar C_{j_4 j_3 j_2 j_1}^{0010},\
\bar C_{j_4 j_3 j_2 j_1}^{0100},\
\bar C_{j_4 j_3 j_2 j_1}^{1000},\
\bar C_{j_6 j_5 j_4 j_3 j_2 j_1}
$$

\vspace{3mm}
\noindent
can be calculated exactly before start of the numerical method (\ref{4.45})
using DERIVE or MAPLE (computer
algebra systems). 
In \cite{2006}, \cite{2017}-\cite{2013}, \cite{arxiv-3}, \cite{2018a}-\cite{2018aaa}
several tables with these coefficients can be found.
Note that the mentioned Fourier--Legendre coefficients
are independent of the integration step $\tau_{p+1}-\tau_p$ of the
numerical scheme,
which can be not a constant in a general case.

Note that in \cite{Mikh-1}, \cite{Kuz-Kuz}
the database with 270,000 exactly
calculated Fourier--Legendre coefficients was described.
This database was used in the software package \cite{Mikh-1}, \cite{Kuz-Kuz},
which is written in the Python programming language
for the implementation of explicit one-step numerical schemes 
with strong orders 0.5, 1.0, 1.5, 2.0, 2.5, and 3.0 
of convergence for Ito SDEs.
The optimization of the mean-square approximation 
procedures for iterated Ito stochastic integrals
from these numerical schemes can be found in \cite{Mikh-2}.

On the basis of 
the presented 
approximations of 
iterated Ito stochastic integrals we 
can see that increasing of multiplicities of these integrals 
leads to increasing 
of orders of smallness with respect to $\tau_{p+1}-\tau_p$
($\tau_{p+1}-\tau_p\ll 1$) in the mean-square sense 
for iterated Ito stochastic integrals. This leads to a sharp decrease  
of member 
quantities
in the approximations of iterated Ito stochastic 
integrals (see the number $q$ in Theorem 3), 
which are required for achieving the acceptable accuracy
of approximation.

\vspace{5mm}

\section{Explicit One-Step Strong Numerical Schemes
of Orders 2.0,
2.5, and 3.0 Based
on the Unified Taylor--Stratonovich expansion}

\vspace{5mm}

Consider the following explicit one-step strong numerical scheme of order 3.0
based on the so-called unified 
Taylor--Stratonovich expansion \cite{kuz33}
(also see \cite{2006}, \cite{2017}-\cite{2013}, \cite{2018a}-\cite{2018aaa})

\vspace{1mm}
$$
{\bf y}_{p+1}={\bf y}_{p}+\sum_{i_{1}=1}^{m}\Sigma_{i_{1}}
\hat I_{0_{\tau_{p+1},\tau_p}}^{*(i_{1})}+\Delta\bar{\bf a}
+\sum_{i_{1},i_{2}=1}^{m}G_0^{(i_{2})}
\Sigma_{i_{1}}\hat I_{00_{\tau_{p+1},\tau_p}}^{*(i_{2}i_{1})}+
$$

\vspace{2mm}
$$
+
\sum_{i_{1}=1}^{m}\left[G_0^{(i_{1})}\bar{\bf a}\left(
\Delta \hat I_{0_{\tau_{p+1},\tau_p}}^{*(i_{1})}+
\hat I_{1_{\tau_{p+1},\tau_p}}^{*(i_{1})}\right)
-\bar L\Sigma_{i_{1}}\hat I_{1_{\tau_{p+1},\tau_p}}^{*(i_{1})}\right]+
$$

\vspace{2mm}
$$
+\sum_{i_{1},i_{2},i_{3}=1}^{m} G_0^{(i_{3})}G_0^{(i_{2})}
\Sigma_{i_{1}}\hat I_{000_{\tau_{p+1},\tau_p}}^{*(i_{3}i_{2}i_{1})}+
\frac{\Delta^2}{2}\bar L\bar{\bf a}+
$$

\vspace{2mm}
$$
+\sum_{i_{1},i_{2}=1}^{m}
\left[G_0^{(i_{2})}\bar L\Sigma_{i_{1}}\left(
\hat I_{10_{\tau_{p+1},\tau_p}}^{*(i_{2}i_{1})}-
\hat I_{01_{\tau_{p+1},\tau_p}}^{*(i_{2}i_{1})}
\right)
-\bar LG_0^{(i_{2})}
\Sigma_{i_{1}}\hat I_{10_{\tau_{p+1},\tau_p}}^{*(i_{2}i_{1})}
+\right.
$$

\vspace{2mm}
$$
\left.+G_0^{(i_{2})}G_0^{(i_{1})}\bar{\bf a}\left(
\hat I_{01_{\tau_{p+1},\tau_p}}
^{*(i_{2}i_{1})}+\Delta \hat I_{00_{\tau_{p+1},\tau_p}}^{*(i_{2}i_{1})}
\right)\right]+
$$

\vspace{2mm}
\begin{equation}
\label{4.470}
+
\sum_{i_{1},i_{2},i_{3},i_{4}=1}^{m}G_0^{(i_{4})}G_0^{(i_{3})}G_0^{(i_{2})}
\Sigma_{i_{1}}\hat
I_{0000_{\tau_{p+1},\tau_p}}^{*(i_{4}i_{3}i_{2}i_{1})}+
{\bf q}_{p+1,p}+{\bf r}_{p+1,p},
\end{equation}

\vspace{10mm}
$$
{\bf q}_{p+1,p}=
\sum_{i_{1}=1}^{m}\Biggl[G_0^{(i_{1})}\bar L\bar{\bf a}\left(\frac{1}{2}
\hat I_{2_{\tau_{p+1},\tau_p}}
^{*(i_{1})}+\Delta \hat I_{1_{\tau_{p+1},\tau_p}}^{*(i_{1})}+
\frac{\Delta^2}{2}\hat I_{0_{\tau_{p+1},\tau_p}}^{*(i_{1})}\right)\Biggr.+
$$

\vspace{2mm}
$$
+\frac{1}{2}\bar L\bar L\Sigma_{i_{1}}\hat
I_{2_{\tau_{p+1},\tau_p}}^{*(i_{1})}-
\bar LG_0^{(i_{1})}\bar{\bf a}\Biggl.
\left(\hat I_{2_{\tau_{p+1},\tau_p}}^{*(i_{1})}+
\Delta \hat I_{1{\tau_{p+1},\tau_p}}^{*(i_{1})}\right)\Biggr]+
$$

\vspace{2mm}
$$
+
\sum_{i_{1},i_{2},i_{3}=1}^m\left[
G_0^{(i_{3})}\bar LG_0^{(i_{2})}\Sigma_{i_{1}}
\left(\hat I_{100_{\tau_{p+1},\tau_p}}
^{*(i_{3}i_{2}i_{1})}-\hat I_{010_{\tau_{p+1},\tau_p}}
^{*(i_{3}i_{2}i_{1})}\right)
\right.+
$$

\vspace{2mm}
$$
+G_0^{(i_{3})}G_0^{(i_{2})}\bar L\Sigma_{i_{1}}\left(
\hat I_{010_{\tau_{p+1},\tau_p}}^{*(i_{3}i_{2}i_{1})}-
\hat I_{001_{\tau_{p+1},\tau_p}}^{*(i_{3}i_{2}i_{1})}\right)+
$$

\vspace{2mm}
$$
+G_0^{(i_{3})}G_0^{(i_{2})}G_0^{(i_{1})}\bar {\bf a}
\left(\Delta \hat I_{000_{\tau_{p+1},\tau_p}}^{*(i_{3}i_{2}i_{1})}+
\hat I_{001_{\tau_{p+1},\tau_p}}^{*(i_{3}i_{2}i_{1})}\right)-
$$

\vspace{2mm}

$$
\left.-\bar LG_0^{(i_{3})}G_0^{(i_{2})}\Sigma_{i_{1}}
\hat I_{100_{\tau_{p+1},\tau_p}}^{*(i_{3}i_{2}i_{1})}\right]+
$$

\vspace{2mm}
$$
+\sum_{i_{1},i_{2},i_{3},i_{4},i_{5}=1}^m
G_0^{(i_{5})}G_0^{(i_{4})}G_0^{(i_{3})}G_0^{(i_{2})}\Sigma_{i_{1}}
\hat I_{00000_{\tau_{p+1},\tau_p}}^{*(i_{5}i_{4}i_{3}i_{2}i_{1})}+
$$

\vspace{2mm}
$$
+
\frac{\Delta^3}{6}\bar L\bar L\bar {\bf a},
$$

\vspace{10mm}
$$
{\bf r}_{p+1,p}=\sum_{i_{1},i_{2}=1}^{m}
\Biggl[G_0^{(i_{2})}G_0^{(i_{1})}\bar L\bar {\bf a}\Biggl(
\frac{1}{2}\hat I_{02_{\tau_{p+1},\tau_p}}^{*(i_{2}i_{1})}
+
\Delta \hat I_{01_{\tau_{p+1},\tau_p}}^{*(i_{2}i_{1})}
+
\frac{\Delta^2}{2}
\hat I_{00_{\tau_{p+1},\tau_p}}^{*(i_{2}i_{1})}\Biggr)+\Biggr.
$$

\vspace{2mm}
$$
+
\frac{1}{2}\bar L\bar LG_0^{(i_{2})}\Sigma_{i_{1}}
\hat I_{20_{\tau_{p+1},\tau_p}}^{*(i_{2}i_{1})}
$$

\vspace{2mm}
$$
+G_0^{(i_{2})}\bar LG_0^{(i_{1})}\bar {\bf a}\left(
\hat I_{11_{\tau_{p+1},\tau_p}}
^{*(i_{2}i_{1})}-\hat I_{02_{\tau_{p+1},\tau_p}}^{*(i_{2}i_{1})}+
\Delta\left(\hat I_{10_{\tau_{p+1},\tau_p}}
^{*(i_{2}i_{1})}-\hat I_{01_{\tau_{p+1},\tau_p}}^{*(i_{2}i_{1})}
\right)\right)+
$$

\vspace{2mm}
$$
+\bar LG_0^{(i_{2})}\bar L\Sigma_{i_1}\left(
\hat I_{11_{\tau_{p+1},\tau_p}}
^{*(i_{2}i_{1})}-\hat I_{20_{\tau_{p+1},\tau_p}}^{*(i_{2}i_{1})}\right)+
$$

\vspace{2mm}
$$
+G_0^{(i_{2})}\bar L\bar L\Sigma_{i_1}\Biggl(
\frac{1}{2}\hat I_{02_{\tau_{p+1},\tau_p}}^{*(i_{2}i_{1})}+
\frac{1}{2}\hat I_{20_{\tau_{p+1},\tau_p}}^{*(i_{2}i_{1})}-
\hat I_{11_{\tau_{p+1},\tau_p}}^{*(i_{2}i_{1})}\Biggr)-
$$

\vspace{2mm}
$$
\Biggl.-\bar LG_0^{(i_{2})}G_0^{(i_{1})}\bar{\bf a}\left(
\Delta \hat I_{10_{\tau_{p+1},\tau_p}}
^{*(i_{2}i_{1})}+\hat I_{11_{\tau_{p+1},\tau_p}}^{*(i_{2}i_{1})}\right)
\Biggr]+
$$

\vspace{2mm}
$$
+
\sum_{i_{1},i_2,i_3,i_{4}=1}^m\Biggl[
G_0^{(i_{4})}G_0^{(i_{3})}G_0^{(i_{2})}G_0^{(i_{1})}\bar{\bf a}
\left(\Delta \hat I_{0000_{\tau_{p+1},\tau_p}}
^{*(i_4i_{3}i_{2}i_{1})}+\hat I_{0001_{\tau_{p+1},\tau_p}}
^{*(i_4i_{3}i_{2}i_{1})}\right)
+\Biggr.
$$

\vspace{2mm}
$$
+G_0^{(i_{4})}G_0^{(i_{3})}\bar LG_0^{(i_{2})}\Sigma_{i_1}
\left(\hat I_{0100_{\tau_{p+1},\tau_p}}
^{*(i_4i_{3}i_{2}i_{1})}-\hat I_{0010_{\tau_{p+1},\tau_p}}
^{*(i_4i_{3}i_{2}i_{1})}\right)-
$$

\vspace{2mm}
$$
-\bar LG_0^{(i_{4})}G_0^{(i_{3})}G_0^{(i_{2})}\Sigma_{i_1}
\hat I_{1000_{\tau_{p+1},\tau_p}}
^{*(i_4i_{3}i_{2}i_{1})}+
$$

\vspace{2mm}
$$
+G_0^{(i_{4})}\bar LG_0^{(i_{3})}G_0^{(i_{2})}\Sigma_{i_1}
\left(\hat I_{1000_{\tau_{p+1},\tau_p}}
^{*(i_4i_{3}i_{2}i_{1})}-\hat I_{0100_{\tau_{p+1},\tau_p}}
^{*(i_4i_{3}i_{2}i_{1})}\right)+
$$

\vspace{2mm}
$$
\Biggl.+G_0^{(i_{4})}G_0^{(i_{3})}G_0^{(i_{2})}\bar L\Sigma_{i_1}
\left(\hat I_{0010_{\tau_{p+1},\tau_p}}
^{*(i_4i_{3}i_{2}i_{1})}-\hat I_{0001_{\tau_{p+1},\tau_p}}
^{*(i_4i_{3}i_{2}i_{1})}\right)\Biggr]+
$$

\vspace{2mm}
$$
+\sum_{i_{1},i_2,i_3,i_4,i_5,i_{6}=1}^m
G_0^{(i_{6})}G_0^{(i_{5})}
G_0^{(i_{4})}G_0^{(i_{3})}G_0^{(i_{2})}\Sigma_{i_{1}}
\hat I_{000000_{\tau_{p+1},\tau_p}}^{*(i_6i_{5}i_{4}i_{3}i_{2}i_{1})},
$$

\vspace{6mm}
\noindent
where $\Delta={\bar T}/N$ $(N>1)$ is a constant (for simplicity) 
step of integration,\
$\tau_p=p\Delta$ $(p=0, 1,\ldots,N)$,\
$\hat 
I_{{l_1\ldots\hspace{0.2mm} l_k}_{\hspace{0.2mm}s,t}}^{*(i_1\ldots i_k)}$ 
is an approximation of the iterated
Stratonovich stochastic integral 

\vspace{-1mm}
\begin{equation}
\label{600}
I_{{l_1\ldots\hspace{0.2mm} l_k}_{\hspace{0.2mm}s,t}}^{*(i_1\ldots i_k)}=
  {\int\limits_t^{*}}^s
(t-\tau_{k})^{l_{k}}\ldots 
{\int\limits_t^{*}}^{\tau_2}
(t-\tau _ {1}) ^ {l_ {1}} d
{\bf f} ^ {(i_ {1})} _ {\tau_ {1}} \ldots 
d {\bf f} _ {\tau_ {k}} ^ {(i_ {k})},
\end{equation}

\vspace{2mm}
$$
\bar{\bf a}({\bf x},t)={\bf a}({\bf x},t)-
\frac{1}{2}\sum\limits_{j=1}^m G_0^{(j)}\Sigma_j({\bf x},t),
$$

\vspace{2mm}
$$
\bar L=L-\frac{1}{2}\sum\limits_{j=1}^m G_0^{(j)}G_0^{(j)}=
\frac{\partial }{\partial t}+
\sum\limits_{j=1}^n \bar {\bf a}^{(j)}({\bf x},t)
\frac{\partial }{\partial {\bf x}^{(j)}},
$$

\vspace{2mm}
$$
L= {\partial \over \partial t}
+ \sum^ {n} _ {i=1} {\bf a}_i ({\bf x},t) 
{\partial  \over  \partial  {\bf  x}_i}
+ {1\over 2} \sum^ {m} _ {j=1} \sum^ {n} _ {l,i=1}
\Sigma_{lj} ({\bf x}, t) \Sigma_{ij} ({\bf x}, t) {\partial
^{2} \over \partial {\bf x}_l \partial {\bf x}_i},
$$

\vspace{3mm}
$$
G_0^{(i)} = \sum^ {n} _ {j=1} \Sigma_{ji} ({\bf x}, t)
{\partial  \over \partial {\bf x}_j},\ \ \
i=1,\ldots,m,
$$

\vspace{4mm}
\noindent
$l_1,\ldots, l_k=0, 1, 2\ldots,$\ \
$i_1,\ldots, i_k=1,\ldots,m,$\ \ $k=1, 2,\ldots$,\ \
$\Sigma_i$ is the $i$th column of the matrix function $\Sigma$ and 
$\Sigma_{ij}$ is the $ij$th component
of the matrix function 
$\Sigma$,  ${\bf a}_i$ is the $i$th component of the vector function 
${\bf a}$ and ${\bf x}_i$ 
is the
$i$th component of the column ${\bf x}$,
the columns 
                    
\vspace{-2mm}
$$
\Sigma_{i_{1}},\ \bar{\bf a},\ G_0^{(i_{2})}\Sigma_{i_{1}},\
G_0^{(i_{1})}\bar{\bf a},\  \bar L\Sigma_{i_{1}},\ G_0^{(i_{3})}G_0^{(i_{2})}\Sigma_{i_{1}},\
\bar L\bar{\bf a},\ G_0^{(i_{2})}\bar L\Sigma_{i_{1}},\ \bar LG_0^{(i_{2})}\Sigma_{i_{1}},\
G_0^{(i_{2})}G_0^{(i_{1})}\bar {\bf a},\ 
$$

\vspace{-4mm}
$$
G_0^{(i_{4})}G_0^{(i_{3})}G_0^{(i_{2})}\Sigma_{i_{1}},\
G_0^{(i_{1})}\bar L \bar {\bf a},\ \bar L \bar L\Sigma_{i_{1}},\
\bar LG_0^{(i_{1})}\bar {\bf a},\ 
G_0^{(i_{3})}\bar LG_0^{(i_{2})}\Sigma_{i_{1}},\ 
G_0^{(i_{3})}G_0^{(i_{2})}\bar L\Sigma_{i_{1}},\ 
G_0^{(i_{3})}G_0^{(i_{2})}G_0^{(i_{1})}\bar {\bf a},\
$$

\vspace{-4mm}
$$
\bar LG_0^{(i_{3})}G_0^{(i_{2})}\Sigma_{i_{1}},\
G_0^{(i_{5})}G_0^{(i_{4})}G_0^{(i_{3})}G_0^{(i_{2})}\Sigma_{i_{1}},\ \bar L \bar L \bar{\bf a},\
G_0^{(i_{2})}G_0^{(i_{1})}\bar L \bar{\bf a},\
\bar L \bar LG_0^{(i_{2})}\Sigma_{i_{1}},\
G_0^{(i_{2})}\bar LG_0^{(i_{1})}\bar {\bf a},\
\bar LG_0^{(i_{2})}\bar L\Sigma_{i_1},\
$$

\vspace{-3mm}
$$
G_0^{(i_{2})}\bar L \bar L\Sigma_{i_1},\
\bar LG_0^{(i_{2})}G_0^{(i_{1})}\bar {\bf a},\
G_0^{(i_{4})}G_0^{(i_{3})}G_0^{(i_{2})}G_0^{(i_{1})}\bar {\bf a},\
G_0^{(i_{4})}G_0^{(i_{3})}\bar LG_0^{(i_{2})}\Sigma_{i_1},\
\bar LG_0^{(i_{4})}G_0^{(i_{3})}G_0^{(i_{2})}\Sigma_{i_1},\
$$

\vspace{-2mm}
$$
G_0^{(i_{4})} \bar LG_0^{(i_{3})}G_0^{(i_{2})}\Sigma_{i_1},\
G_0^{(i_{4})}G_0^{(i_{3})}G_0^{(i_{2})}\bar L\Sigma_{i_1},\ 
G_0^{(i_{6})}G_0^{(i_{5})}G_0^{(i_{4})}G_0^{(i_{3})}G_0^{(i_{2})}\Sigma_{i_{1}}
$$

\vspace{4mm}
\noindent
are calculated at the point $({\bf y}_p,p).$

It is well known \cite{KlPl2} that under the standard conditions
the numerical scheme (\ref{4.470}) has strong order of convergence 3.0. 
Among these conditions we consider only the condition
for approximations of iterated Stratonovich stochastic 
integrals from the numerical
scheme (\ref{4.470}) \cite{KlPl2}, \cite{2006}

\vspace{-1mm}
$$
{\sf M}\left\{\Biggl(I_{{l_{1}\ldots\hspace{0.2mm} l_{k}}_{\tau_{p+1},\tau_p}}
^{*(i_{1}\ldots i_{k})} 
-\hat
I_{{l_{1}\ldots\hspace{0.2mm} l_{k}}_{\tau_{p+1},\tau_p}}^{*(i_{1}\ldots i_{k})}
\Biggr)^2\right\}\le C\Delta^{7},
$$

\vspace{2mm}
\noindent
where 
$\hat 
I_{{l_1\ldots\hspace{0.2mm} 
l_k}_{\hspace{0.2mm}\tau_{p+1},\tau_p}}^{*(i_1\ldots i_k)}$
is an approximation 
of 
$I_{{l_{1}\ldots\hspace{0.2mm} 
l_{k}}_{\hspace{0.2mm}\tau_{p+1},\tau_p}}^{*(i_{1}\ldots i_{k})},$
constant $C$ does not depend on $\Delta$.

Note that if we exclude ${\bf q}_{p+1,p}+{\bf r}_{p+1,p}$ from the
right-hand side of (\ref{4.470}), then we will have the explicit 
one-step strong numerical scheme of order 2.0.
The right-hand side of (\ref{4.470}) but without the value
${\bf r}_{p+1,p}$ and with replacing the value 
$\Delta^3\bar L\bar L\bar {\bf a}/6$
by the value $\Delta^3 LL{\bf a}/6$ 
define 
the explicit 
one-step strong numerical scheme of order 2.5.

Note that the 
truncated 
unified 
Taylor--Stratonovich expansion \cite{kuz33}
(also see \cite{2006}, \cite{2017}-\cite{2013}, \cite{2018a}-\cite{2018aaa})
contains 
the less number of various 
types of iterated Stratonovich 
stochastic integrals (moreover, their major part 
will have 
less multiplicities) in comparison with 
the classical
Taylor--Stratonovich expansion \cite{KlPl2}, \cite{KlPl1}.

Furthermore, some iterated stochastic integrals from 
the Taylor--Stratonovich expansion \cite{KlPl2}, \cite{KlPl1}
are connected by linear relations. 
However, the iterated stochastic integrals from the 
unified 
Taylor--Stratonovich expansion \cite{kuz33}
(also see \cite{2006}, \cite{2017}-\cite{2013}, \cite{2018a}-\cite{2018aaa})
cannot be connected by linear relations. 
Therefore, we call these families of stochastic integrals 
as the stochastic 
bases \cite{2006}, \cite{2017}-\cite{2013}, \cite{2018a}-\cite{2018aaa}.
Note that (\ref{4.470}) contains 20 different types
of iterated Stratonovich stochastic integrals. At the same time,
the analogue of (\ref{4.470}) based on 
the classical
Taylor--Stratonovich expansion \cite{KlPl2}, \cite{KlPl1} 
contains 29 different types 
of iterated stochastic integrals.

\vspace{5mm}

\section{Fourier--Legendre Expansions of Iterated Stratonovich 
Stochastic Integrals of Multiplicities 1 to 6}

\vspace{5mm}

As noted above, in a number of works of the author 
Theorems 1, 2 have been adapted for the iterated  Stratonovich stochastic integrals
(\ref{str}) of multiplicities 1 to 6 (the case of multiplicity 1 is given by (\ref{a1})).
Let us first present some old results.

\vspace{2mm}

{\bf Theorem 4}\ \cite{2017}-\cite{2017-1a}, 
\cite{2010-2}-\cite{2013}, \cite{arxiv-5}, \cite{arxiv-7}, 
\cite{2018a}-\cite{2018aaa}.\ {\it Assume that 
the following conditions are fulfilled{\rm :}

{\rm 1}. The function $\psi_2(\tau)$ is continuously 
differentiable at the interval $[t, T]$ and the
function $\psi_1(\tau)$ is two times continuously 
differentiable at the interval $[t, T]$.

{\rm 2}. $\{\phi_j(x)\}_{j=0}^{\infty}$ is a complete orthonormal system 
of Legendre polynomials or system of tri\-go\-no\-met\-ric functions
in the space $L_2([t, T]).$

Then, the iterated Stratonovich stochastic integral 
of multiplicity {\rm 2}

\vspace{-1mm}
$$
J^{*}[\psi^{(2)}]_{T,t}={\int\limits_t^{*}}^T\psi_2(t_2)
{\int\limits_t^{*}}^{t_2}\psi_1(t_1)d{\bf f}_{t_1}^{(i_1)}
d{\bf f}_{t_2}^{(i_2)}\ \ \ (i_1, i_2=1,\ldots,m)
$$

\vspace{2mm}
\noindent
is expanded into the 
converging in the mean-square sense multiple series

$$
J^{*}[\psi^{(2)}]_{T,t}=
\hbox{\vtop{\offinterlineskip\halign{
\hfil#\hfil\cr
{\rm l.i.m.}\cr
$\stackrel{}{{}_{p_1,p_2\to \infty}}$\cr
}} }\sum_{j_1=0}^{p_1}\sum_{j_2=0}^{p_2}
C_{j_2j_1}\zeta_{j_1}^{(i_1)}\zeta_{j_2}^{(i_2)},
$$

\vspace{4mm}
\noindent
where the meaning of notations introduced in the formulation
of Theorem {\rm 1} is 
remained.}

\vspace{2mm}

Proving the theorem 4 \cite{2017}-\cite{2017-1a}, 
\cite{2010-2}-\cite{2013}, \cite{arxiv-5}, \cite{arxiv-7},
\cite{2018a}-\cite{2018aaa}
we used Theorem 1 and double integration by parts. This
procedure leads to the condition of double 
continuous
differentiability of the
function $\psi_1(\tau)$ at the interval $[t, T]$.
The mentioned condition can be weakened \cite{arxiv100},
\cite{5-006}, \cite{arxiv-8}, \cite{arxiv-10}, \cite{2018a}-\cite{2018aaa}
and Theorem 4 will be valid
for continuously 
differentiable functions $\psi_l(\tau)$ $(l=1, 2)$
at the interval $[t, T]$. 

\vspace{2mm}

{\bf Theorem 5}\ \cite{2017}-\cite{2017-1a}, 
\cite{2010-2}-\cite{2013}, \cite{arxiv-5}, \cite{arxiv-7},
\cite{2018a}-\cite{2018aaa}. 
{\it Assume, that
$\{\phi_j(x)\}_{j=0}^{\infty}$ is a complete ortho\-nor\-mal
system of Legendre polynomials or tri\-go\-no\-mer\-tic functions
in the space $L_2([t, T])$. Furthermore,
the function $\psi_2(\tau)$ is continuously
differentiable at the interval $[t, T]$ and
the functions $\psi_1(\tau),$ $\psi_3(\tau)$ are twice continuously
differentiable at the interval $[t, T]$.
Then, for the iterated Stratonovich stochastic integral of multiplicity {\rm 3}

$$
J^{*}[\psi^{(3)}]_{T,t}={\int\limits_t^{*}}^T\psi_3(t_3)
{\int\limits_t^{*}}^{t_3}\psi_2(t_2)
{\int\limits_t^{*}}^{t_2}\psi_1(t_1)
d{\bf f}_{t_1}^{(i_1)}
d{\bf f}_{t_2}^{(i_2)}d{\bf f}_{t_3}^{(i_3)}\ \ \ (i_1, i_2, i_3=1,\ldots,m)
$$

\vspace{4mm}
\noindent
the following 
expansion 

\vspace{-1mm}
\begin{equation}
\label{feto19000a}
J^{*}[\psi^{(3)}]_{T,t}=
\hbox{\vtop{\offinterlineskip\halign{
\hfil#\hfil\cr
{\rm l.i.m.}\cr
$\stackrel{}{{}_{p\to \infty}}$\cr
}} }
\sum\limits_{j_1, j_2, j_3=0}^{p}
C_{j_3 j_2 j_1}\zeta_{j_1}^{(i_1)}\zeta_{j_2}^{(i_2)}\zeta_{j_3}^{(i_3)}
\end{equation}

\vspace{2mm}
\noindent
converging in the mean-square sense is valid, where

\vspace{-1mm}
$$
C_{j_3 j_2 j_1}=\int\limits_t^T\psi_3(t_3)\phi_{j_3}(t_3)
\int\limits_t^{t_3}\psi_2(t_2)\phi_{j_2}(t_2)
\int\limits_t^{t_2}\psi_1(t_1)\phi_{j_1}(t_1)dt_1dt_2dt_3;
$$

\vspace{2mm}
\noindent
another notations are the same as in Theorems {\rm 1, 2}.}

\vspace{2mm}

{\bf Theorem 6}\ \cite{2017}--\cite{2017-1a}, 
\cite{2010-2}-\cite{2013}, \cite{arxiv-5}, \cite{arxiv-6},
\cite{2018a}-\cite{2018aaa}. 
{\it Assume that
$\{\phi_j(x)\}_{j=0}^{\infty}$ is a complete ortho\-nor\-mal
system of Legendre polynomials or trigonometric functions
in the space $L_2([t, T])$.
Then, for the iterated Stratonovich stochastic integral of 
multiplicity {\rm 4}

$$
I_{T,t}^{*(i_1 i_2 i_3 i_4)}=
{\int\limits_t^{*}}^T
{\int\limits_t^{*}}^{t_4}
{\int\limits_t^{*}}^{t_3}
{\int\limits_t^{*}}^{t_2}
d{\bf w}_{t_1}^{(i_1)}
d{\bf w}_{t_2}^{(i_2)}d{\bf w}_{t_3}^{(i_3)}d{\bf w}_{t_4}^{(i_4)},\ 
$$

\vspace{3mm}
\noindent
where $i_1, i_2, i_3, i_4=0, 1,\ldots,m,$
the following 
expansion

$$
I_{T,t}^{*(i_1 i_2 i_3 i_4)}=
\hbox{\vtop{\offinterlineskip\halign{
\hfil#\hfil\cr
{\rm l.i.m.}\cr
$\stackrel{}{{}_{p\to \infty}}$\cr
}} }
\sum\limits_{j_1, j_2, j_3, j_4=0}^{p}
C_{j_4 j_3 j_2 j_1}\zeta_{j_1}^{(i_1)}\zeta_{j_2}^{(i_2)}\zeta_{j_3}^{(i_3)}
\zeta_{j_4}^{(i_4)}
$$

\vspace{3mm}
\noindent
converging in the mean-square sense is valid, where

$$
C_{j_4 j_3 j_2 j_1}=\int\limits_t^T\phi_{j_4}(t_4)\int\limits_t^{t_4}
\phi_{j_3}(t_3)
\int\limits_t^{t_3}\phi_{j_2}(t_2)\int\limits_t^{t_2}\phi_{j_1}(t_1)
dt_1dt_2dt_3dt_4;
$$

\vspace{3mm}
\noindent
${\bf w}_{\tau}^{(i)}={\bf f}_{\tau}^{(i)}$ $(i=1,\ldots,m)$ are independent 
standard Wiener processes and 
${\bf w}_{\tau}^{(0)}=\tau;$
another notations are the same as in Theorems {\rm 1, 2}.}

\vspace{2mm}

Recently, a new approach to the expansion and mean-square 
approximation of iterated Stratonovich stochastic integrals has been obtained
\cite{2018a} (Sect.~2.10--2.16), \cite{arxiv-5} (Sect.~13--19), 
\cite{arxiv-6} (Sect.~5--11), \cite{arxiv-4} (Sect.~7--13), \cite{new-art-1-xxy}
(Sections~4--9).
Let us formulate four theorems that were obtained using this approach.

\vspace{2mm}

{\bf Theorem 7}\ \cite{arxiv-5}, 
\cite{arxiv-6}, \cite{arxiv-4}, \cite{2018a}, \cite{new-art-1-xxy}.\ {\it Suppose 
that $\{\phi_j(x)\}_{j=0}^{\infty}$ is a complete orthonormal system of 
Legendre polynomials or trigonometric functions in the space $L_2([t, T]).$
Furthermore, let $\psi_1(\tau), \psi_2(\tau),$ $\psi_3(\tau)$ are continuously dif\-ferentiable 
nonrandom functions on $[t, T].$ 
Then, for the 
iterated Stra\-to\-no\-vich stochastic integral of third multiplicity

$$
J^{*}[\psi^{(3)}]_{T,t}={\int\limits_t^{*}}^T\psi_3(t_3)
{\int\limits_t^{*}}^{t_3}\psi_2(t_2)
{\int\limits_t^{*}}^{t_2}\psi_1(t_1)
d{\bf w}_{t_1}^{(i_1)}
d{\bf w}_{t_2}^{(i_2)}d{\bf w}_{t_3}^{(i_3)}\ \ \ (i_1,i_2,i_3=0,1,\ldots,m)
$$

\vspace{4mm}
\noindent
the following 
relations

\vspace{-1mm}
\begin{equation}
\label{fin1}
J^{*}[\psi^{(3)}]_{T,t}
=\hbox{\vtop{\offinterlineskip\halign{
\hfil#\hfil\cr
{\rm l.i.m.}\cr
$\stackrel{}{{}_{p\to \infty}}$\cr
}} }
\sum\limits_{j_1, j_2, j_3=0}^{p}
C_{j_3 j_2 j_1}\zeta_{j_1}^{(i_1)}\zeta_{j_2}^{(i_2)}\zeta_{j_3}^{(i_3)},
\end{equation}

\vspace{3mm}
\begin{equation}
\label{fin2}
{\sf M}\left\{\left(
J^{*}[\psi^{(3)}]_{T,t}-
\sum\limits_{j_1, j_2, j_3=0}^{p}
C_{j_3 j_2 j_1}\zeta_{j_1}^{(i_1)}\zeta_{j_2}^{(i_2)}\zeta_{j_3}^{(i_3)}\right)^2\right\}
\le \frac{C}{p}
\end{equation}

\vspace{5mm}
\noindent
are fulfilled, where $i_1, i_2, i_3=0,1,\ldots,m$ in {\rm (\ref{fin1})} and 
$i_1, i_2, i_3=1,\ldots,m$ in {\rm (\ref{fin2})},
constant $C$ is independent of $p,$

$$
C_{j_3 j_2 j_1}=\int\limits_t^T\psi_3(t_3)\phi_{j_3}(t_3)
\int\limits_t^{t_3}\psi_2(t_2)\phi_{j_2}(t_2)
\int\limits_t^{t_2}\psi_1(t_1)\phi_{j_1}(t_1)dt_1dt_2dt_3
$$

\vspace{4mm}
\noindent
and
$$
\zeta_{j}^{(i)}=
\int\limits_t^T \phi_{j}(\tau) d{\bf f}_{\tau}^{(i)}
$$ 

\vspace{2mm}
\noindent
are independent standard Gaussian random variables for various 
$i$ or $j$ {\rm (}in the case when $i\ne 0${\rm );} 
another notations are the same as in Theorems~{\rm 1, 2}.}

\vspace{2mm}

{\bf Theorem 8}\ \cite{arxiv-5}, 
\cite{arxiv-6}, \cite{arxiv-4}, \cite{2018a}, \cite{new-art-1-xxy}.\ {\it Let
$\{\phi_j(x)\}_{j=0}^{\infty}$ be a complete orthonormal system of 
Legendre polynomials or trigonometric functions in the space $L_2([t, T]).$
Furthermore, let $\psi_1(\tau), \ldots,$ $\psi_4(\tau)$ be continuously dif\-ferentiable 
nonrandom functions on $[t, T].$ 
Then, for the 
iterated Stra\-to\-no\-vich stochastic integral of fourth multiplicity

\begin{equation}
\label{fin0}
J^{*}[\psi^{(4)}]_{T,t}={\int\limits_t^{*}}^T\psi_4(t_4)
{\int\limits_t^{*}}^{t_4}\psi_3(t_3)
{\int\limits_t^{*}}^{t_3}\psi_2(t_2)
{\int\limits_t^{*}}^{t_2}\psi_1(t_1)
d{\bf w}_{t_1}^{(i_1)}
d{\bf w}_{t_2}^{(i_2)}d{\bf w}_{t_3}^{(i_3)}d{\bf w}_{t_4}^{(i_4)}
\end{equation}

\vspace{4mm}
\noindent
the following 
relations

\begin{equation}
\label{fin3}
J^{*}[\psi^{(4)}]_{T,t}
=\hbox{\vtop{\offinterlineskip\halign{
\hfil#\hfil\cr
{\rm l.i.m.}\cr
$\stackrel{}{{}_{p\to \infty}}$\cr
}} }
\sum\limits_{j_1, j_2, j_3,j_4=0}^{p}
C_{j_4j_3 j_2 j_1}\zeta_{j_1}^{(i_1)}\zeta_{j_2}^{(i_2)}\zeta_{j_3}^{(i_3)}\zeta_{j_4}^{(i_4)},
\end{equation}

\vspace{3mm}

\begin{equation}
\label{fin4}
{\sf M}\left\{\left(
J^{*}[\psi^{(4)}]_{T,t}-
\sum\limits_{j_1, j_2, j_3, j_4=0}^{p}
C_{j_4 j_3 j_2 j_1}\zeta_{j_1}^{(i_1)}\zeta_{j_2}^{(i_2)}\zeta_{j_3}^{(i_3)}
\zeta_{j_4}^{(i_4)}
\right)^2\right\}
\le \frac{C}{p^{1-\varepsilon}}
\end{equation}

\vspace{5mm}
\noindent
are fulfilled, where $i_1, \ldots , i_4=0,1,\ldots,m$ in {\rm (\ref{fin0}),} {\rm (\ref{fin3})} 
and $i_1, \ldots, i_4=1,\ldots,m$ in {\rm (\ref{fin4}),}
constant $C$ does not depend on $p,$
$\varepsilon$ is an arbitrary
small positive real number 
for the case of complete orthonormal system of 
Legendre polynomials in the space $L_2([t, T])$
and $\varepsilon=0$ for the case of
complete orthonormal system of 
trigonometric functions in the space $L_2([t, T]),$

$$
C_{j_4 j_3 j_2 j_1}=
$$

$$
=
\int\limits_t^T\psi_4(t_4)\phi_{j_4}(t_4)
\int\limits_t^{t_4}\psi_3(t_3)\phi_{j_3}(t_3)
\int\limits_t^{t_3}\psi_2(t_2)\phi_{j_2}(t_2)
\int\limits_t^{t_2}\psi_1(t_1)\phi_{j_1}(t_1)dt_1dt_2dt_3dt_4;
$$

\vspace{4mm}
\noindent
another notations are the same as in Theorem~{\rm 7}.}

\vspace{2mm}

{\bf Theorem 9}\ \cite{arxiv-5}, 
\cite{arxiv-6}, \cite{arxiv-4}, \cite{2018a}, \cite{new-art-1-xxy}.\
{\it Assume 
that $\{\phi_j(x)\}_{j=0}^{\infty}$ is a complete orthonormal system of 
Legendre polynomials or trigonometric functions in the space $L_2([t, T])$
and $\psi_1(\tau), \ldots,$ $\psi_5(\tau)$ are continuously dif\-ferentiable 
nonrandom functions on $[t, T].$ 
Then, for the 
iterated Stra\-to\-no\-vich stochastic integral of fifth multiplicity

\begin{equation}
\label{fin7}
J^{*}[\psi^{(5)}]_{T,t}={\int\limits_t^{*}}^T\psi_5(t_5)
\ldots
{\int\limits_t^{*}}^{t_2}\psi_1(t_1)
d{\bf w}_{t_1}^{(i_1)}
\ldots d{\bf w}_{t_5}^{(i_5)}
\end{equation}

\vspace{4mm}
\noindent
the following 
relations

\begin{equation}
\label{fin8}
J^{*}[\psi^{(5)}]_{T,t}
=\hbox{\vtop{\offinterlineskip\halign{
\hfil#\hfil\cr
{\rm l.i.m.}\cr
$\stackrel{}{{}_{p\to \infty}}$\cr
}} }
\sum\limits_{j_1,\ldots,j_5=0}^{p}
C_{j_5 \ldots j_1}\zeta_{j_1}^{(i_1)}\ldots \zeta_{j_5}^{(i_5)},
\end{equation}

\vspace{3mm}

\begin{equation}
\label{fin9}
{\sf M}\left\{\left(
J^{*}[\psi^{(5)}]_{T,t}-
\sum\limits_{j_1, \ldots, j_5=0}^{p}
C_{j_5 \ldots j_1}\zeta_{j_1}^{(i_1)}\ldots
\zeta_{j_5}^{(i_5)}
\right)^2\right\}
\le \frac{C}{p^{1-\varepsilon}}
\end{equation}

\vspace{5mm}
\noindent
are fulfilled, where $i_1, \ldots , i_5=0,1,\ldots,m$ in {\rm (\ref{fin7}),} {\rm (\ref{fin8})} 
and $i_1, \ldots, i_5=1,\ldots,m$ in {\rm (\ref{fin9}),}
constant $C$ is independent of $p,$
$\varepsilon$ is an arbitrary
small positive real number 
for the case of complete orthonormal system of 
Legendre polynomials in the space $L_2([t, T])$
and $\varepsilon=0$ for the case of
complete orthonormal system of 
trigonometric functions in the space $L_2([t, T]),$

$$
C_{j_5 \ldots j_1}=
\int\limits_t^T\psi_5(t_5)\phi_{j_5}(t_5)\ldots
\int\limits_t^{t_2}\psi_1(t_1)\phi_{j_1}(t_1)dt_1\ldots dt_5;
$$

\vspace{3mm}
\noindent
another notations are the same as in Theorems~{\rm 7, 8}.}

\vspace{2mm}

{\bf Theorem 10}\ \cite{arxiv-5}, 
\cite{arxiv-6}, \cite{arxiv-4}, \cite{2018a}.\
{\it Suppose that 
$\{\phi_j(x)\}_{j=0}^{\infty}$ is a complete orthonormal system of 
Legendre polynomials or trigonometric functions in the space $L_2([t, T]).$
Then, for the 
iterated Stratonovich stochastic integral of sixth multiplicity

\begin{equation}
\label{after10001qu1}
I_{T,t}^{*(i_1\ldots i_6)}={\int\limits_t^{*}}^T
\ldots
{\int\limits_t^{*}}^{t_2}
d{\bf w}_{t_1}^{(i_1)}
\ldots d{\bf w}_{t_6}^{(i_6)}
\end{equation}

\vspace{3mm}
\noindent
the following 
expansion 

\vspace{-1mm}
$$
I_{T,t}^{*(i_1\ldots i_6)}
=\hbox{\vtop{\offinterlineskip\halign{
\hfil#\hfil\cr
{\rm l.i.m.}\cr
$\stackrel{}{{}_{p\to \infty}}$\cr
}} }
\sum\limits_{j_1, \ldots, j_6=0}^{p}
C_{j_6 \ldots j_1}\zeta_{j_1}^{(i_1)}\ldots
\zeta_{j_6}^{(i_6)}
$$

\vspace{4mm}
\noindent
that converges in the mean-square sense is valid, where
$i_1, \ldots, i_6=0, 1,\ldots,m,$

$$
C_{j_6 \ldots j_1}=
\int\limits_t^T\phi_{j_6}(t_6)\ldots
\int\limits_t^{t_2}\phi_{j_1}(t_1)dt_1\ldots dt_6;
$$

\vspace{3mm}
\noindent
another notations are the same as in Theorems~{\rm 7--9}.}

\vspace{2mm}

On the base of Theorems 4--10 the folloing hypothesis was formulated
in \cite{2017}-\cite{2017-1a}, \cite{2010-2}-\cite{2013}, \cite{arxiv-4},
\cite{2018a}-\cite{2018aaa}.

\vspace{2mm}

{\bf Hypothesis 1}\ \cite{2017}-\cite{2017-1a}, \cite{2010-2}-\cite{2013}, 
\cite{arxiv-4}, \cite{2018a}-\cite{2018aaa}. {\it Assume that
$\{\phi_j(x)\}_{j=0}^{\infty}$ is a complete orthonormal
system of Legendre polynomials or trigonometric functions
in the space $L_2([t, T])$.
Moreover,
every $\psi_l(\tau)$ $(l=1,\ldots,k)$ is 
an enough smooth nonrandom function
on $[t,T].$
Then, for the iterated Stratonovich stochastic 
integral 
$J^{*}[\psi^{(k)}]_{T,t}$ defined by {\rm (\ref{str})}
the following expansion 
         
\vspace{-1mm}
\begin{equation}
\label{feto1900otita}
J^{*}[\psi^{(k)}]_{T,t}=
\hbox{\vtop{\offinterlineskip\halign{
\hfil#\hfil\cr
{\rm l.i.m.}\cr
$\stackrel{}{{}_{p\to \infty}}$\cr
}} }
\sum\limits_{j_1, \ldots j_k=0}^{p}
C_{j_k \ldots j_1}\zeta_{j_1}^{(i_1)}
\ldots
\zeta_{j_k}^{(i_k)}
\end{equation}

\vspace{4mm}
\noindent
converging in the mean-square sense 
is valid, where the notations are the same as in Theorems {\rm 1, 2}.}

\vspace{2mm}

Hypothesis 1 allows to approximate the iterated
Stratonovich stochastic integral 
$J^{*}[\psi^{(k)}]_{T,t}$
by the sum
$$
J^{*}[\psi^{(k)}]_{T,t}^p=
\sum\limits_{j_1,\ldots j_k=0}^{p}
C_{j_k \ldots j_1}\prod\limits_{l=1}^k
\zeta_{j_l}^{(i_l)},
$$

\vspace{3mm}
\noindent
where
$$
\lim_{p\to\infty}{\sf M}\left\{\Biggl(
J^{*}[\psi^{(k)}]_{T,t}-
J^{*}[\psi^{(k)}]_{T,t}^p\Biggl)^2\right\}=0.
$$

\vspace{4mm}

Note that Hypothesis~1 is proved in \cite{2018a} (Sect.~2.10)
under the condition of convergence of trace series
(also see \cite{arxiv-5}, 
\cite{arxiv-6}, \cite{arxiv-4}).
In \cite{arxiv-4}, \cite{2018a}-\cite{2018aaa} 
a more general hypothesis
is formulated.

Applying Theorems 4--10,
we obtain the following 
approximations of iterated Stratonovich stochastic integrals 
from (\ref{4.470})

\vspace{2mm}
$$
I_{0_{\tau_{p+1},\tau_p}}^{*(i_1)}=\sqrt{\Delta}\zeta_0^{(i_1)},
$$

\vspace{2mm}

$$
I_{{00}_{\tau_{p+1},\tau_p}}^{*(i_1 i_2)q}=
\frac{\Delta}{2}\left(\zeta_0^{(i_1)}\zeta_0^{(i_2)}+\sum_{i=1}^{q}
\frac{1}{\sqrt{4i^2-1}}\left(
\zeta_{i-1}^{(i_1)}\zeta_{i}^{(i_2)}-
\zeta_i^{(i_1)}\zeta_{i-1}^{(i_2)}\right)\right),
$$

\vspace{4mm}

$$
I_{{1}_{\tau_{p+1},\tau_p}}^{*(i_1)}=
-\frac{{\Delta}^{3/2}}{2}\left(\zeta_0^{(i_1)}+
\frac{1}{\sqrt{3}}\zeta_1^{(i_1)}\right),
$$

\vspace{4mm}

$$
I_{{000}_{\tau_{p+1},\tau_p}}^{*(i_1i_2i_3)q}
=\sum_{j_1,j_2,j_3=0}^{q}
C_{j_3j_2j_1}
\zeta_{j_1}^{(i_1)}\zeta_{j_2}^{(i_2)}\zeta_{j_3}^{(i_3)},
$$

\vspace{7mm}

$$
I_{01_{\tau_{p+1},\tau_p}}^{*(i_1 i_2)q}=-\frac{\Delta}{2}
I_{00_{\tau_{p+1},\tau_p}}^{*(i_1 i_2)q}
-\frac{{\Delta}^2}{4}\Biggl[
\frac{1}{\sqrt{3}}\zeta_0^{(i_1)}\zeta_1^{(i_2)}+\Biggr.
$$

\vspace{1mm}
$$
+\Biggl.\sum_{i=0}^{q}\Biggl(
\frac{(i+2)\zeta_i^{(i_1)}\zeta_{i+2}^{(i_2)}
-(i+1)\zeta_{i+2}^{(i_1)}\zeta_{i}^{(i_2)}}
{\sqrt{(2i+1)(2i+5)}(2i+3)}-
\frac{\zeta_i^{(i_1)}\zeta_{i}^{(i_2)}}{(2i-1)(2i+3)}\Biggr)\Biggr],
$$

\vspace{7mm}

$$
I_{10_{\tau_{p+1},\tau_p}}^{*(i_1 i_2)q}=
-\frac{\Delta}{2}I_{00_{\tau_{p+1},\tau_p}}^{*(i_1 i_2)q}
-\frac{\Delta^2}{4}\Biggl[
\frac{1}{\sqrt{3}}\zeta_0^{(i_2)}\zeta_1^{(i_1)}+\Biggr.
$$

\vspace{1mm}
$$
+\Biggl.\sum_{i=0}^{q}\Biggl(
\frac{(i+1)\zeta_{i+2}^{(i_2)}\zeta_{i}^{(i_1)}
-(i+2)\zeta_{i}^{(i_2)}\zeta_{i+2}^{(i_1)}}
{\sqrt{(2i+1)(2i+5)}(2i+3)}+
\frac{\zeta_i^{(i_1)}\zeta_{i}^{(i_2)}}{(2i-1)(2i+3)}\Biggr)\Biggr],
$$

\vspace{7mm}

$$
I_{{0000}_{\tau_{p+1},\tau_p}}^{*(i_1 i_2 i_3 i_4)q}
=\sum_{j_1, j_2, j_3, j_4=0}^{q}
C_{j_4 j_3 j_2 j_1}
\zeta_{j_1}^{(i_1)}\zeta_{j_2}^{(i_2)}\zeta_{j_3}^{(i_3)}
\zeta_{j_4}^{(i_4)},
$$

\vspace{5mm}

$$
{I}_{2_{\tau_{p+1},\tau_p}}^{*(i_1)}=
\frac{\Delta^{5/2}}{3}\left(
\zeta_0^{(i_1)}+\frac{\sqrt{3}}{2}\zeta_1^{(i_1)}+
\frac{1}{2\sqrt{5}}\zeta_2^{(i_1)}\right),
$$

\vspace{5mm}

$$
I_{{100}_{\tau_{p+1},\tau_p}}^{*(i_1 i_2 i_3)q}
=\sum_{j_1, j_2, j_3=0}^{q}
C_{j_3 j_2 j_1}^{100}
\zeta_{j_1}^{(i_1)}\zeta_{j_2}^{(i_2)}\zeta_{j_3}^{(i_3)},
$$

\vspace{5mm}
$$
I_{{010}_{\tau_{p+1},\tau_p}}^{*(i_1 i_2 i_3)q}
=\sum_{j_1, j_2, j_3=0}^{q}
C_{j_3 j_2 j_1}^{010}
\zeta_{j_1}^{(i_1)}\zeta_{j_2}^{(i_2)}\zeta_{j_3}^{(i_3)},
$$

\vspace{5mm}
$$
I_{{001}_{\tau_{p+1},\tau_p}}^{*(i_1 i_2 i_3)q}
=\sum_{j_1, j_2, j_3=0}^{q}
C_{j_3 j_2 j_1}^{001}
\zeta_{j_1}^{(i_1)}\zeta_{j_2}^{(i_2)}\zeta_{j_3}^{(i_3)},
$$

\vspace{5mm}
$$
I_{{00000}_{\tau_{p+1},\tau_p}}^{*(i_1 i_2 i_3 i_4 i_5)q}=
\sum\limits_{j_1, j_2, j_3, j_4, j_5=0}^{q}
C_{j_5j_4 j_3 j_2 j_1}\zeta_{j_1}^{(i_1)}\zeta_{j_2}^{(i_2)}\zeta_{j_3}^{(i_3)}
\zeta_{j_4}^{(i_4)}\zeta_{j_5}^{(i_5)},
$$

\vspace{9mm}

$$
I_{02_{\tau_{p+1},\tau_p}}^{*(i_1 i_2)q}
=-\frac{{\Delta}^2}{4}I_{00_{\tau_{p+1},\tau_p}}^{*(i_1 i_2)q}
-\Delta I_{01_{\tau_{p+1},\tau_p}}^{*(i_1 i_2)q}+
\frac{{\Delta}^3}{8}\Biggl[
\frac{2}{3\sqrt{5}}\zeta_2^{(i_2)}\zeta_0^{(i_1)}+\Biggr.
$$

\vspace{1mm}
$$
+\frac{1}{3}\zeta_0^{(i_1)}\zeta_0^{(i_2)}+
\sum_{i=0}^{q}\Biggl(
\frac{(i+2)(i+3)\zeta_{i+3}^{(i_2)}\zeta_{i}^{(i_1)}
-(i+1)(i+2)\zeta_{i}^{(i_2)}\zeta_{i+3}^{(i_1)}}
{\sqrt{(2i+1)(2i+7)}(2i+3)(2i+5)}+
\Biggr.
$$

\vspace{1mm}
$$
\Biggl.\Biggl.+
\frac{(i^2+i-3)\zeta_{i+1}^{(i_2)}\zeta_{i}^{(i_1)}
-(i^2+3i-1)\zeta_{i}^{(i_2)}\zeta_{i+1}^{(i_1)}}
{\sqrt{(2i+1)(2i+3)}(2i-1)(2i+5)}\Biggr)\Biggr],
$$

\vspace{9mm}

$$
I_{20_{\tau_{p+1},\tau_p}}^{*(i_1 i_2)q}=-\frac{{\Delta}^2}{4}
I_{00_{\tau_{p+1},\tau_p}}^{*(i_1 i_2)q}
-\Delta I_{10_{\tau_{p+1},\tau_p}}^{*(i_1 i_2)q}+
\frac{{\Delta}^3}{8}\Biggl[
\frac{2}{3\sqrt{5}}\zeta_0^{(i_2)}\zeta_2^{(i_1)}+\Biggr.
$$

\vspace{1mm}
$$
+\frac{1}{3}\zeta_0^{(i_1)}\zeta_0^{(i_2)}+
\sum_{i=0}^{q}\Biggl(
\frac{(i+1)(i+2)\zeta_{i+3}^{(i_2)}\zeta_{i}^{(i_1)}
-(i+2)(i+3)\zeta_{i}^{(i_2)}\zeta_{i+3}^{(i_1)}}
{\sqrt{(2i+1)(2i+7)}(2i+3)(2i+5)}+
\Biggr.
$$

\vspace{1mm}
$$
\Biggl.\Biggl.+
\frac{(i^2+3i-1)\zeta_{i+1}^{(i_2)}\zeta_{i}^{(i_1)}
-(i^2+i-3)\zeta_{i}^{(i_2)}\zeta_{i+1}^{(i_1)}}
{\sqrt{(2i+1)(2i+3)}(2i-1)(2i+5)}\Biggr)\Biggr],
$$

\vspace{9mm}

$$
I_{11_{\tau_{p+1},\tau_p}}^{*(i_1 i_2)q}
=-\frac{{\Delta}^2}{4}I_{00_{\tau_{p+1},\tau_p}}^{*(i_1 i_2)q}
-\frac{\Delta}{2}\left(
I_{10_{\tau_{p+1},\tau_p}}^{*(i_1 i_2)q}+
I_{01_{\tau_{p+1},\tau_p}}^{*(i_1 i_2)q}\right)+
\frac{{\Delta}^3}{8}\Biggl[
\frac{1}{3}\zeta_1^{(i_1)}\zeta_1^{(i_2)}+\Biggr.
$$

\vspace{1mm}
$$
+
\sum_{i=0}^{q}\Biggl(
\frac{(i+1)(i+3)\left(\zeta_{i+3}^{(i_2)}\zeta_{i}^{(i_1)}
-\zeta_{i}^{(i_2)}\zeta_{i+3}^{(i_1)}\right)}
{\sqrt{(2i+1)(2i+7)}(2i+3)(2i+5)}+
\Biggr.
$$

\vspace{2mm}
$$
\Biggl.\Biggl.+
\frac{(i+1)^2\left(\zeta_{i+1}^{(i_2)}\zeta_{i}^{(i_1)}
-\zeta_{i}^{(i_2)}\zeta_{i+1}^{(i_1)}\right)}
{\sqrt{(2i+1)(2i+3)}(2i-1)(2i+5)}\Biggr)\Biggr],
$$

\vspace{7mm}

$$
I_{{0001}_{\tau_{p+1},\tau_p}}^{*(i_1i_2i_3)q}
=
\sum_{j_1,j_2,j_3,j_4=0}^{q}
C_{j_4j_3 j_2 j_1}^{0001}
\zeta_{j_1}^{(i_1)}\zeta_{j_2}^{(i_2)}\zeta_{j_3}^{(i_3)}\zeta_{j_4}^{(i_4)},
$$

\vspace{4mm}
$$
I_{{0010}_{\tau_{p+1},\tau_p}}^{*(i_1i_2i_3)q}
=
\sum_{j_1,j_2,j_3,j_4=0}^{q}
C_{j_4j_3 j_2 j_1}^{0010}
\zeta_{j_1}^{(i_1)}\zeta_{j_2}^{(i_2)}\zeta_{j_3}^{(i_3)}\zeta_{j_4}^{(i_4)},
$$

\vspace{4mm}
$$
I_{{0100}_{\tau_{p+1},\tau_p}}^{*(i_1i_2i_3)q}
=
\sum_{j_1,j_2,j_3,j_4=0}^{q}
C_{j_4j_3 j_2 j_1}^{0100}
\zeta_{j_1}^{(i_1)}\zeta_{j_2}^{(i_2)}\zeta_{j_3}^{(i_3)}\zeta_{j_4}^{(i_4)},
$$

\vspace{4mm}
$$
I_{{1000}_{\tau_{p+1},\tau_p}}^{*(i_1i_2i_3)q}
=
\sum_{j_1,j_2,j_3,j_4=0}^{q}
C_{j_4j_3 j_2 j_1}^{1000}
\zeta_{j_1}^{(i_1)}\zeta_{j_2}^{(i_2)}\zeta_{j_3}^{(i_3)}\zeta_{j_4}^{(i_4)},
$$

\vspace{4mm}
$$
I_{{000000}_{\tau_{p+1},\tau_p}}^{*(i_1 i_2 i_3 i_4 i_5 i_6)q}=
\sum\limits_{j_1, j_2, j_3, j_4, j_5, j_6=0}^{q}
C_{j_5j_4 j_3 j_2 j_1}\zeta_{j_1}^{(i_1)}\zeta_{j_2}^{(i_2)}\zeta_{j_3}^{(i_3)}
\zeta_{j_4}^{(i_4)}\zeta_{j_5}^{(i_5)}\zeta_{j_5}^{(i_5)},
$$

\vspace{7mm}

\noindent
where formulas for the Fourier--Legendre coefficients 

\vspace{1mm}
$$
C_{j_3 j_2 j_1},\ C_{j_4 j_3 j_2 j_1},\ C_{j_3 j_2 j_1}^{001},\ 
C_{j_3 j_2 j_1}^{010},\ C_{j_3 j_2 j_1}^{100},\
C_{j_5 j_4 j_3 j_2 j_1},
C_{j_4 j_3 j_2 j_1}^{0001},\
C_{j_4 j_3 j_2 j_1}^{0010},\
$$

\vspace{1mm}
$$
C_{j_4 j_3 j_2 j_1}^{0100},\
C_{j_4 j_3 j_2 j_1}^{1000},\
C_{j_6 j_5 j_4 j_3 j_2 j_1}
$$

\vspace{4mm}
\noindent
can be found in Sect.~3.

On the basis of 
the presented 
approximations of 
iterated Stratonovich stochastic integrals we 
can see that increasing of multiplicities of these integrals 
leads to increasing 
of orders of smallness with respect to $\tau_{p+1}-\tau_p$
($\tau_{p+1}-\tau_p\ll 1$) in the mean-square sense 
for iterated Stratonovich stochastic integrals. This leads to a sharp decrease  
of member 
quantities
in the appro\-xi\-ma\-ti\-ons of iterated Stratonovich stochastic 
integrals (see the numbers $q$ in the approximations
of iterated Stratonovich stochastic integrals from this section), 
which are required for achieving the acceptable accuracy
of approximation.

From (\ref{qqqq1}) $(i_1\ne i_2)$ we have

$$
{\sf M}\left\{\left(I_{{00}_{\tau_{p+1},\tau_p}}^{*(i_1 i_2)}-
I_{{00}_{\tau_{p+1},\tau_p}}^{*(i_1 i_2)q}
\right)^2\right\}=\frac{\Delta^2}{2}
\sum\limits_{i=q+1}^{\infty}\frac{1}{4i^2-1}\le 
$$

\begin{equation}
\label{teac}
\le \frac{\Delta^2}{2}\int\limits_{q}^{\infty}
\frac{1}{4x^2-1}dx
=-\frac{\Delta^2}{8}{\rm ln}\left|
1-\frac{2}{2q+1}\right|\le C_1\frac{\Delta^2}{q},
\end{equation}

\vspace{3mm}
\noindent
where $C_1$ is a constant.

Since the value $\Delta=\tau_{p+1}-\tau_p$ plays the role of integration step 
in the numerical scheme (\ref{4.470}), 
then this value is a sufficiently small.

Keeping in mind this circumstance, it is easy to notice that there 
exists such a constant $C_2$ that

\vspace{-1mm}
\begin{equation}
\label{teac3}
{\sf M}\left\{\left(
I_{{l_1\ldots l_k}_{\hspace{0.2mm}\tau_{p+1},\tau_p}}^{*(i_1\ldots i_k)}-
I_{{l_1\ldots l_k}_{\hspace{0.2mm}\tau_{p+1},\tau_p}}^{*(i_1\ldots i_k)q}\right)^2\right\}
\le C_2 {\sf M}\left\{\left(I_{{00}_{\tau_{p+1},\tau_p}}^{*(i_1 i_2)}-
I_{{00}_{\tau_{p+1},\tau_p}}^{*(i_1 i_2)q}\right)^2\right\},
\end{equation}

\vspace{3mm}
\noindent
where 
$I_{{l_1\ldots l_k}_{\hspace{0.2mm}\tau_{p+1},\tau_p}}^{*(i_1\ldots i_k)q}$
is an approximation of the iterated Stratonovich stochastic integral 
$I_{{l_1\ldots l_k}_{\hspace{0.2mm}\tau_{p+1},\tau_p}}^{*(i_1\ldots i_k)}.$

From (\ref{teac}) and (\ref{teac3}) we finally obtain

\vspace{-1mm}
\begin{equation}
\label{teac4}
{\sf M}\left\{\left(
I_{{l_1\ldots l_k}_{\hspace{0.2mm}\tau_{p+1},\tau_p}}^{*(i_1\ldots i_k)}-
I_{{l_1\ldots l_k}_{\hspace{0.2mm}\tau_{p+1},
\tau_p}}^{*(i_1\ldots i_k)q}\right)^2\right\}
\le C \frac{\Delta^2}{q},
\end{equation}

\vspace{3mm}
\noindent
where constant $C$ does not depend on $\Delta$.

The same idea can be found in \cite{KlPl2} for the case 
of trigonometric functions.
Note that, in contrast to the estimate (\ref{teac4}), 
the constant $C$ in Theorems 7--9 does not depend on $q.$

Since 
$$
J^{*}[\psi^{(k)}]_{T,t}=J[\psi^{(k)}]_{T,t}\ \ \ \hbox{w.\ p.\ 1}
$$

\vspace{4mm}
\noindent
for pairwise different $i_1,\ldots,i_k=1,\ldots,m$, then we can write 
for pairwise different $i_1,\ldots,i_6=1,\ldots,m$ (see (\ref{qq1}))

\vspace{2mm}

$$
{\sf M}\left\{\left(
I_{{01}_{\tau_{p+1},\tau_p}}^{*(i_1i_2)}-
I_{{01}_{\tau_{p+1},\tau_p}}^{*(i_1i_2)q}\right)^2\right\}=
\frac{\Delta^{4}}{4}-\sum_{j_1,j_2=0}^{q}
\left(C_{j_2j_1}^{01}\right)^2,
$$

\vspace{3mm}
$$
{\sf M}\left\{\left(
I_{{10}_{\tau_{p+1},\tau_p}}^{*(i_1i_2)}-
I_{{10}_{\tau_{p+1},\tau_p}}^{*(i_1i_2)q}\right)^2\right\}=
\frac{\Delta^{4}}{12}-\sum_{j_1,j_2=0}^{q}
\left(C_{j_2j_1}^{10}\right)^2,
$$

\vspace{3mm}

$$
{\sf M}\left\{\left(
I_{{000}_{\tau_{p+1},\tau_p}}^{*(i_1i_2 i_3)}-
I_{{000}_{\tau_{p+1},\tau_p}}^{*(i_1i_2 i_3)q}\right)^2\right\}=
\frac{\Delta^{3}}{6}-\sum_{j_3,j_2,j_1=0}^{q}
C_{j_3j_2j_1}^2,
$$

\vspace{3mm}

$$
{\sf M}\left\{\left(
I_{{0000}_{\tau_{p+1},\tau_p}}^{*(i_1i_2 i_3 i_4)}-
I_{{0000}_{\tau_{p+1},\tau_p}}^{*(i_1i_2 i_3 i_4)q}\right)^2\right\}=
\frac{\Delta^{4}}{24}-\sum_{j_1,j_2,j_3,j_4=0}^{q}
C_{j_4j_3j_2j_1}^2,
$$

\vspace{3mm}

$$
{\sf M}\left\{\left(
I_{{100}_{\tau_{p+1},\tau_p}}^{*(i_1i_2 i_3)}-
I_{{100}_{\tau_{p+1},\tau_p}}^{*(i_1i_2 i_3)q}\right)^2\right\}=
\frac{\Delta^{5}}{60}-\sum_{j_1,j_2,j_3=0}^{q}
\left(C_{j_3j_2j_1}^{100}\right)^2,
$$

\vspace{3mm}

$$
{\sf M}\left\{\left(
I_{{010}_{\tau_{p+1},\tau_p}}^{*(i_1i_2 i_3)}-
I_{{010}_{\tau_{p+1},\tau_p}}^{*(i_1i_2 i_3)q}\right)^2\right\}=
\frac{\Delta^{5}}{20}-\sum_{j_1,j_2,j_3=0}^{q}
\left(C_{j_3j_2j_1}^{010}\right)^2,
$$

\vspace{3mm}

$$
{\sf M}\left\{\left(
I_{{001}_{\tau_{p+1},\tau_p}}^{*(i_1i_2 i_3)}-
I_{{001}_{\tau_{p+1},\tau_p}}^{*(i_1i_2 i_3)q}\right)^2\right\}=
\frac{\Delta^5}{10}-\sum_{j_1,j_2,j_3=0}^{q}
\left(C_{j_3j_2j_1}^{001}\right)^2,
$$

\vspace{3mm}

$$
{\sf M}\left\{\left(
I_{{00000}_{\tau_{p+1},\tau_p}}^{*(i_1 i_2 i_3 i_4 i_5)}-
I_{{00000}_{\tau_{p+1},\tau_p}}^{*(i_1 i_2 i_3 i_4 i_5)q}\right)^2\right\}=
\frac{\Delta^{5}}{120}-\sum_{j_1,j_2,j_3,j_4,j_5=0}^{q}
C_{j_5 i_4 i_3 i_2 j_1}^2,
$$

\vspace{3mm}

$$
{\sf M}\left\{\left(
I_{{20}_{\tau_{p+1},\tau_p}}^{*(i_1i_2)}-
I_{{20}_{\tau_{p+1},\tau_p}}^{*(i_1i_2)q}\right)^2\right\}=
\frac{\Delta^6}{30}-\sum_{j_2,j_1=0}^{q}
\left(C_{j_2j_1}^{20}\right)^2,
$$

\vspace{3mm}

$$
{\sf M}\left\{\left(
I_{{11}_{\tau_{p+1},\tau_p}}^{*(i_1i_2)}-
I_{{11}_{\tau_{p+1},\tau_p}}^{*(i_1i_2)q}\right)^2\right\}=
\frac{\Delta^6}{18}-\sum_{j_2,j_1=0}^{q}
\left(C_{j_2j_1}^{11}\right)^2,
$$

\vspace{3mm}

$$
{\sf M}\left\{\left(
I_{{02}_{\tau_{p+1},\tau_p}}^{*(i_1i_2)}-
I_{{02}_{\tau_{p+1},\tau_p}}^{*(i_1i_2)q}\right)^2\right\}=
\frac{\Delta^6}{6}-\sum_{j_2,j_1=0}^{q}
\left(C_{j_2j_1}^{02}\right)^2,
$$

\vspace{3mm}

$$
{\sf M}\left\{\left(
I_{{1000}_{\tau_{p+1},\tau_p}}^{*(i_1i_2 i_3i_4)}-
I_{{1000}_{\tau_{p+1},\tau_p}}^{*(i_1i_2 i_3i_4)q}\right)^2\right\}=
\frac{\Delta^{6}}{360}-\sum_{j_1,j_2,j_3, j_4=0}^{q}
\left(C_{j_4j_3j_2j_1}^{1000}\right)^2,
$$

\vspace{3mm}

$$
{\sf M}\left\{\left(
I_{{0100}_{\tau_{p+1},\tau_p}}^{*(i_1i_2 i_3i_4)}-
I_{{0100}_{\tau_{p+1},\tau_p}}^{*(i_1i_2 i_3i_4)q}\right)^2\right\}=
\frac{\Delta^{6}}{120}-\sum_{j_1,j_2,j_3, j_4=0}^{q}
\left(C_{j_4j_3j_2j_1}^{0100}\right)^2,
$$

\vspace{3mm}

$$
{\sf M}\left\{\left(
I_{{0010}_{\tau_{p+1},\tau_p}}^{*(i_1i_2 i_3i_4)}-
I_{{0010}_{\tau_{p+1},\tau_p}}^{*(i_1i_2 i_3 i_4)q}\right)^2\right\}=
\frac{\Delta^6}{60}-\sum_{j_1,j_2,j_3, j_4=0}^{q}
\left(C_{j_4j_3j_2j_1}^{0010}\right)^2,
$$

\vspace{3mm}

$$
{\sf M}\left\{\left(
I_{{0001}_{\tau_{p+1},\tau_p}}^{*(i_1i_2 i_3 i_4)}-
I_{{0001}_{\tau_{p+1},\tau_p}}^{*(i_1i_2 i_3 i_4)q}\right)^2\right\}=
\frac{\Delta^6}{36}-\sum_{j_1,j_2,j_3, j_4=0}^{q}
\left(C_{j_4j_3j_2j_1}^{0001}\right)^2,
$$

\vspace{3mm}

$$
{\sf M}\left\{\left(
I_{{000000}_{\tau_{p+1},\tau_p}}^{*(i_1 i_2 i_3 i_4 i_5 i_6)}-
I_{{000000}_{\tau_{p+1},\tau_p}}^{*(i_1 i_2 i_3 i_4 i_5 i_6)q}\right)^2\right\}=
\frac{\Delta^{6}}{720}-\sum_{j_1,j_2,j_3,j_4,j_5,j_6=0}^{q}
C_{j_6 j_5 j_4 j_3 j_2 j_1}^2.
$$

\vspace{6mm}

For example \cite{2006} (also see \cite{2017}-\cite{2013},
\cite{2018a}-\cite{2018aaa})

\vspace{3mm}

$$
{\sf M}\left\{\left(
I_{{000}_{\tau_{p+1}},\tau_p}^{*(i_1i_2 i_3)}-
I_{{000}_{\tau_{p+1}},\tau_p}^{*(i_1i_2 i_3)6}\right)^2\right\}=
\frac{\Delta^{3}}{6}-\sum_{j_3,j_2,j_1=0}^{6}
C_{j_3j_2j_1}^2
\approx
0.01956000\Delta^3,
$$

\vspace{3mm}

$$
{\sf M}\left\{\left(
I_{{100}_{\tau_{p+1},\tau_p}}^{*(i_1i_2 i_3)}-
I_{{100}_{\tau_{p+1},\tau_p}}^{*(i_1i_2 i_3)2}\right)^2\right\}=
\frac{\Delta^{5}}{60}-\sum_{j_1,j_2,j_3=0}^{2}
\left(C_{j_3j_2j_1}^{100}\right)^2
\approx
0.00815429\Delta^5,
$$

\vspace{3mm}

$$
{\sf M}\left\{\left(
I_{{010}_{\tau_{p+1},\tau_p}}^{*(i_1i_2 i_3)}-
I_{{010}_{\tau_{p+1},\tau_p}}^{*(i_1i_2 i_3)2}\right)^2\right\}=
\frac{\Delta^{5}}{20}-\sum_{j_1,j_2,j_3=0}^{2}
\left(C_{j_3j_2j_1}^{010}\right)^2
\approx
0.01739030\Delta^5,
$$

\vspace{3mm}

$$
{\sf M}\left\{\left(
I_{{001}_{\tau_{p+1},\tau_p}}^{*(i_1i_2 i_3)}-
I_{{001}_{\tau_{p+1},\tau_p}}^{*(i_1i_2 i_3)2}\right)^2\right\}=
\frac{\Delta^5}{10}-\sum_{j_1,j_2,j_3=0}^{2}
\left(C_{j_3j_2j_1}^{001}\right)^2
\approx 0.02528010\Delta^5,
$$

\vspace{3mm}

$$
{\sf M}\left\{\left(
I_{{0000}_{\tau_{p+1},\tau_p}}^{*(i_1i_2i_3 i_4)}-
I_{{0000}_{\tau_{p+1},\tau_p}}^{*(i_1i_2i_3 i_4)2}\right)^2\right\}=
\frac{\Delta^{4}}{24}-\sum_{j_1,j_2,j_3,j_4=0}^{2}
C_{j_4 j_3 j_2 j_1}^2\approx
0.02360840\Delta^4,
$$

\vspace{3mm}

$$
{\sf M}\left\{\left(
I_{{00000}_{\tau_{p+1},\tau_p}}^{*(i_1i_2i_3i_4 i_5)}-
I_{{00000}_{\tau_{p+1},\tau_p}}^{*(i_1i_2i_3i_4 i_5)1}\right)^2\right\}=
\frac{\Delta^5}{120}-\sum_{j_1,j_2,j_3,j_4,j_5=0}^{1}
C_{j_5 i_4 i_3 i_2 j_1}^2\approx
0.00759105\Delta^5.
$$

\vspace{6mm}

The theory presented in this article was realized \cite{Mikh-1}, \cite{Kuz-Kuz}
in the form of a software package in the Python programming language.
The mentioned software package implements the 
strong numerical methods with convergence orders 
0.5, 1.0, 1.5, 2.0, 2.5, and 3.0 for Ito
SDEs with multidimensional non-commutative noise
based on the unified Taylor--Ito
and Taylor--Stratonovich expansions.
At that for the numerical simulation of iterated 
Ito and Stratonovich stochastic integrals of multiplicities 
1 to 6 we applied the formulas 
based on multiple Fourier--Legendre series \cite{Mikh-1}, \cite{Kuz-Kuz}.
Moreover, we used \cite{Mikh-1}, \cite{Kuz-Kuz}
the database with 270,000 exactly calculated
Fourier--Legendre coefficients.

\end{document}